\documentclass[twoside,11pt]{article}

\usepackage{jmlr2e}
\usepackage{xcolor}

\usepackage{amsmath}
\usepackage{mathtools,subcaption,bbm}
\newtheorem{assumption}{Assumption}
\newtheorem{assumpB}{Assumption}

\newtheorem{comment}{Comment}

\newcommand{\R}{\mathbb{R}}
\newcommand{\E}{\mathbb{E}}
\newcommand{\Xout}{X^{\text{out}}}
\newcommand{\Xin}{X^{\text{in}}}
\newcommand{\constrtheta}{\bar{\theta}^\dagger}
\newcommand{\constrvartheta}{\bar{\vartheta}^\dagger}
\newcommand{\constrz}{ \bar{z}^\dagger}
\DeclareMathOperator*{\argmin}{arg\,min}

\jmlrheading{1}{2022}{1-??}{6/22}{01/23?}{SteinLeng}{Stefan Stein and Chenlei Leng}

\ShortHeadings{An Annotated Graph Model}{Stein and Leng}
\firstpageno{1}

\usepackage{lastpage}
\jmlrheading{24}{2023}{1-\pageref{LastPage}}{10/22; Revised
3/23}{4/23}{22-1138}{Stefan Stein and Chenlei Leng}
\ShortHeadings{An Annotated Graph Model}{Stein and Leng}

\begin{document}

\title{An Annotated Graph Model with Differential Degree Heterogeneity for Directed Networks}

\author{\name Stefan Stein \email s.stein@warwick.ac.uk
       \AND Chenlei Leng \email c.leng@warwick.ac.uk\\
       \addr University of Warwick\\
       Coventry, CV4 7AL, UK
    }

\editor{Pradeep Ravikumar}

\maketitle

\begin{abstract}
Directed networks are conveniently represented as graphs in which ordered edges encode interactions between vertices. Despite their wide availability, there is a shortage of statistical models amenable for inference, specially when contextual information and degree heterogeneity are present. This paper presents an annotated graph model with parameters explicitly accounting for these features. To overcome the curse of dimensionality due to modelling degree heterogeneity, we introduce a sparsity assumption  and propose a penalized likelihood approach with $\ell_1$-regularization for parameter estimation. We study the estimation and selection consistency of this approach under a sparse network assumption, and show that inference on the covariate parameter is straightforward, thus bypassing the need for the kind of debiasing commonly employed in $\ell_1$-penalized likelihood estimation. Simulation and data analysis corroborate our theoretical findings. 
\end{abstract}


\section{Introduction}

\begin{figure}[htbp]
	\begin{center}
		\includegraphics[width=3.5in]{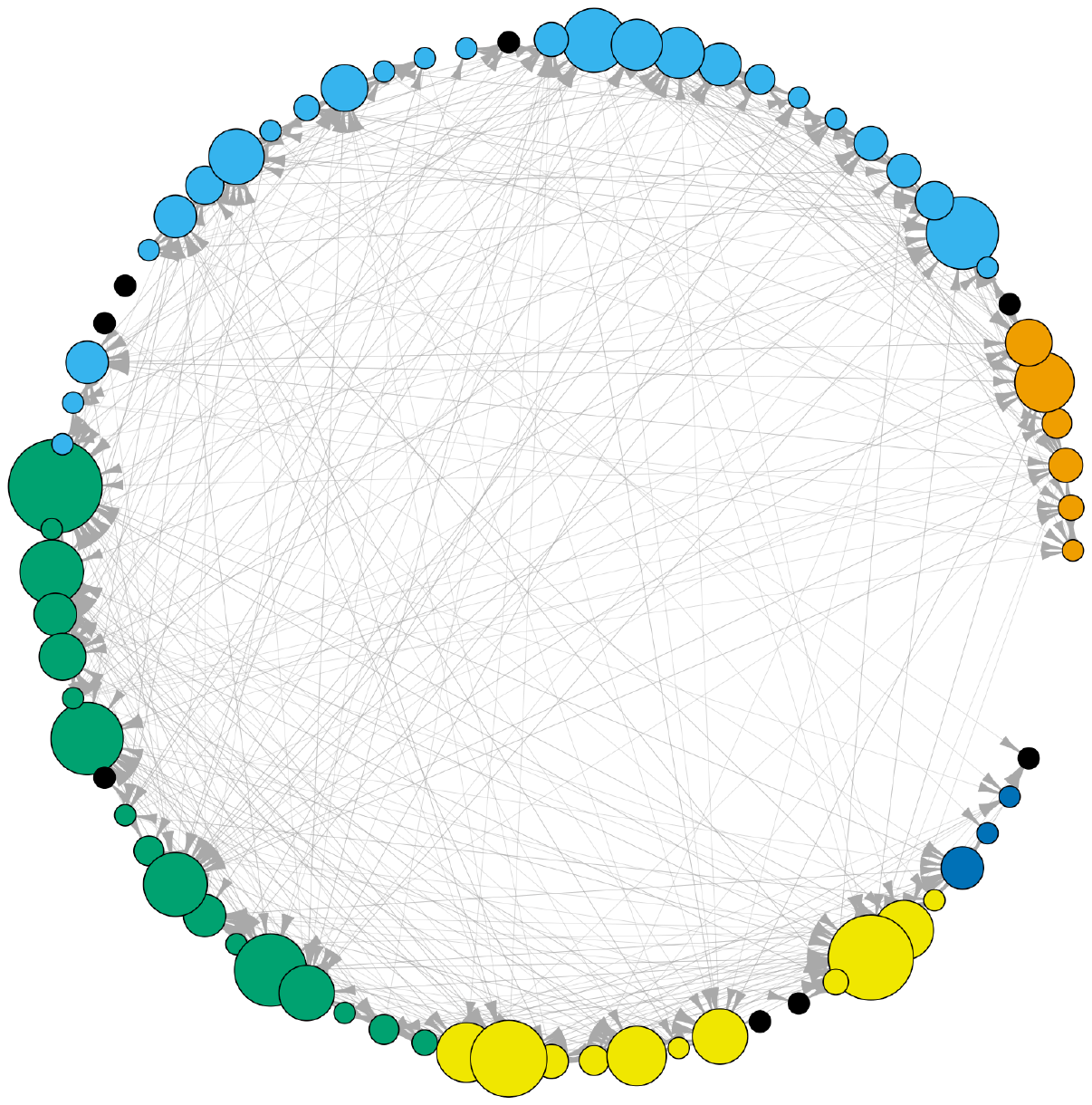}
		\caption{Lazega's lawyers friendship network. The size of the nodes corresponds to their in-degrees. For better visibility all nodes with an in-degree of five or less are plotted with the same size. The 71 lawyers are colour-coded by their age group: The lawyers aged 20-29 are represented in orange, those aged 30-39 in light-blue, those aged 40-49 in green, followed by the lawyers aged 50-59 in yellow and finally those lawyers aged 60 or older in dark-blue. The eight nodes in black correspond to lawyers with either zero in- or out-degree or both. 		}
		\label{fig0}
	\end{center}
\end{figure}

The need to examine inter-relationship of multiple entities in data is rapidly growing due to the increasing availability of datasets that can be conveniently represented as networks or graphs. This paper concerns a new random graph model for describing networks with directed edges.  
As a motivating example, Figure \ref{fig0} depicts the lawyer friendship data in \cite{lawyernetwork} in which 71 lawyers were asked to name their friends: 
An edge from node $i$ to node $j$ exists if and only if lawyer $i$ indicated in a survey that they socialized with lawyer $j$ outside work. 
Mathematically, we denote the adjacency matrix of a network on $n$ nodes as a binary matrix $A \in \R^{n \times n}$, where $n$ is the number of nodes and $A_{ij} = 1$ if $i$ points to $j$ and $0$ otherwise. 	 As is common for this type of data, the lawyer data is annotated with contextual information in the form of nodal covariates. These covariates include a lawyer's status (partner or associate), their gender (man or woman), which of three offices they worked in, the number of years they had spent with the firm, their age, their practice (litigation or corporate) and the law school they had visited (Harvard and Yale, UConn or other). The main interest here is to understand how links were formed, especially the effect of the covariates between pairs of nodes for forming ties.

The lawyer network features several stylized facts of a typical real-life network. First, it exhibits degree heterogeneity, the different tendency that the nodes in a network have in participating in network activities as can be seen from Figure \ref{fig0}. Second, the overall network is sparse, in that the observed number of ties does not scale proportionally to the total number of possible links. In Figure \ref{fig0}, the average in- (and out-) degree is 8.1, whereas the maximum possible value  is 70. Third, the contextual information in terms of the covariates has a role to play in determining how nodes are connected, as can be seen from the data analysis later. Here the covariates will be denoted as $Z_{ij} \in \R^p$ for an edge linking node $i$ to $j$. In case where we only have nodal covariates denoted as $X_i \in \R^p $ for the $i$th node, a common approach is to define $Z_{ij}$ as a function of $X_i$ and $X_j$ that measures their (dis)similarities.

This paper proposes a new annotated graph model that can effectively deal with the above features in directed networks.  This model postulates that links are independently made with the linking probability between node $i$ and $j$  as
\begin{equation}
	p_{ij}= P(A_{ij} = 1 \vert Z_{ij} )  = \frac{\exp(\alpha_i +\beta_j + \mu + \gamma^TZ_{ij})}{1 + \exp(\alpha_i + \beta_j + \mu + \gamma^TZ_{ij})},
	\label{Eq: srgm}
\end{equation}
where $\mu\in\R$, $\alpha_i\in\R$, $\beta_j\in\R$ and $\gamma\in\R^p$ are parameters. 	For identifiability, we assume $\min_i \{\alpha_i\} = \min_j \{\beta_j \} = 0$, because otherwise this model becomes trivially the model in \cite{Yan:etal:2019} by absorbing $\mu$ into $\alpha_i$ and $\beta_j$. Clearly, to fit a model with these many parameters, we will require the number of edges of a network to be relatively large which renders \eqref{Eq: srgm} less useful or not even applicable for the kind of networks encountered in practice. To reduce the dimensionality of the model, we assume that both $\alpha=(\alpha_1, ..., \alpha_n)^T$ and $\beta=(\beta_1, ..., \beta_n)^T$ are sparse. {\color{black} While it may seem appealing to impose $\min_{1 \le i \le n } \vert \alpha_i \vert = 0, \min_{1 \le j \le n } \vert \beta_j \vert = 0$ instead of $\min_i \{\alpha_i\} = \min_j \{\beta_j \} = 0$, restricting the degree heterogeneity parameters only in absolute value would result in an unidentifiable parameter.} 

In \eqref{Eq: srgm}, $\mu$ can be seen as a global density parameter that is allowed to diverge to $-\infty$ aiming to model sparse networks. The two node-specific parameters $\alpha_i$ and $\beta_i$ are used to explicitly capture out- and in-degree heterogeneity respectively. The sparsity assumption on $\alpha$ and $\beta$ introduces a notion of differential degree heterogeneity, in the sense that we only include them for nodes that are important. The effect of covariates is captured by $\gamma$. When a covariate encodes the similarity of a node attribute, a positive $\gamma$  implies homophily, the tendency of nodes similar in attributes to connect, which is a widely observed  phenomenon in real-life networks. 

Statistical analysis of random networks has attracted enormous research attention in recent years thanks to a deluge of network data \citep{Kolaczyk:2009,Fienberg:2012,Kolaczyk:2017}. 	Among others,	there are several major classes of statistical models that have been successfully applied to model degree heterogeneity, including the stochastic block model aiming to identify groups of nodes as communities \citep{Holland:etal:1983,Bickel:Chen:2009,karrer2011stochastic,rohe2011spectral,zhao2012consistency,lei2015consistency,gao2017achieving,amini2018semidefinite,Abbe:2018,zhang2021directed}, the $\beta$-model assigning individual parameters to different nodes \citep{Chatterjee:etal:2011,Yan:Xu:2013,Yan:etal:2016,Chen:etal:19}, and the exponential random graph models using network motifs as sufficient statistics \citep{Holland:Leinhardt:1981,frank1986markov,robins2007introduction}.

The increasing prevalence of covariates in network data calls for models that can effectively account for their effects.  For the stochastic block model, we point to \cite{binkiewicz2017covariate,zhang2016community,Huang:Feng:2018,Yan:Sarkar:2020}. For the $\beta$-model, see \cite{Graham:2017,Yan:etal:2019}, and \cite{Stein:2020}. Additional references include \cite{mama:2020} where a latent space model approach is investigated and 	\cite{zhang2021directed} that proposes a model-free approach to study the dependence of links on covariates. 

\cite{Yan:etal:2019} proposed a model similar to \eqref{Eq: srgm} with the critical difference that $\alpha$ and $\beta$ are dense. As a result, their model only handles relatively dense networks, and the inference on $\gamma$, the covariate parameter, warrants a bias correction step. Remarkable, the inference for this parameter in our model does not require debiasing, even when the network has vanishing link probabilities and hence are sparse. The techniques involved in deriving this result extends substantially many similar ones developed, for example in $M$-estimation \citep{Vaart:1998} and  LASSO theory \citep{Ravikumar:2010, vandegeer2011,vandegeer2014}, where the probabilities in similar models are typically assumed to be bounded away from zero.

Another model similar to \eqref{Eq: srgm} has been developed for undirected networks by \cite{Stein:2020}. The model in this paper is more complex due to the presence of two sets of heterogeneity parameters and thus more delicate analysis is needed \citep{Yan:etal:2016, Yan:etal:2019}. 
More importantly, this paper focuses more on variable selection while  \cite{Stein:2020} focused exclusively on estimation consistency. In particular, we place special focus on the study of the interplay of the rates of convergence for the network sparsity, the parameter sparsity and the penalty we use. Each of our main results requires these three quantities to be balanced in a suitable manner, made explicit in our \textit{balancing assumptions} (Assumptions \ref{Assum: model selection assumption}, \ref{Assum: rate of s^* directed}, \ref{Assum: new rate of s and rho_n directed}). These assumptions highlight  the effect of different sparsity regimes much more clearly than \cite{Stein:2020}.

{\color{black}It has been claimed by \cite{Barabasi:2016} that the degrees of real world networks often follow a power-law distribution. We show in the Appendix \ref{sec:power law} that the degree sequence of the \textit{sparse $\beta$-model} in \cite{Chen:etal:19} (as a representative of the family of models introduced in \cite{Chen:etal:19, Stein:2020} and this paper), follows a power-law distribution under the right assumptions.}

\subsection{Notation}
Denote $N = n(n-1)$ and $[n] \coloneqq \{1, \dots, n\}$. Denote $\R_+=[0, +\infty)$ as the non-negative real line. For a vector $v \in \R^n$, we use $S(v)=\{i: v_i \not=0 \}$ to denote its support and $\text{diag}(v) \in \R^{n \times n}$ the $n$-by-$n$ diagonal matrix with $v$ on the diagonal. 
Let $\Vert \, . \, \Vert_1, \Vert \, . \, \Vert_2. \Vert \, . \, \Vert_\infty$ denote the vector $\ell_1$-, $\ell_2$- and $\ell_\infty$-norm respectively and define $\| \, . \, \|_0$ denotes the $\ell_0$-``norm''. That is, $\| v \|_0 = | S(v) |$. For a vector $v \in \R^{N}$, we number its elements as $v = (v_{ij})_{i \neq j }$. 

For brevity, we write $\vartheta = (\alpha^T, \beta^T)^T$ and $\xi = (\mu, \gamma^T)^T \in \R^{p+1}$. Thus, $\theta = (\vartheta^T, \xi^T)^T $ 
with its true value denoted as $\theta_0 = (\vartheta^{T}_0, \xi_0^T)^T$. 
We write $S_0=S(\vartheta_0)$ and denote its cardinality by $s_0 = \vert S_0 \vert$. We write 
$S_{0,+} \coloneqq S_0 \cup \{2n+1, 2n+2, \dots 2n+1+p\}$ with cardinality $s_{0,+} = \vert S_{0,+} \vert = s_0 + p+1$ to refer to all active indices including those of $\mu$ and $\gamma$. Thus, $s_0$ and $s_{0,+}$ can be understood as the parameter sparsity of our model.  Let $S_\alpha =  \{i: \alpha_{0,i} > 0 \}, S_\beta = \{j : \beta_{0,j} > 0\}$ and $s_\alpha = \vert S_\alpha \vert, s_\beta = \vert S_\beta \vert$. 
When we want to make the dependence of the link probabilities given $Z_{ij}$ on different values of $\theta$ explicit, we write $p_{ij}(\theta) = \frac{\exp(\alpha_i + \beta_j + \mu + \gamma^T Z_{ij})}{1 + \exp(\beta_i + \beta_j + \mu + \gamma^T Z_{ij})}$.  	
Finally, we use $C$ for some generic, strictly positive constant that may change between displays.

{\color{black} For the convenience of the reader, we provide a (non-exhaustive) table of the most important quantities encountered in this paper in Table \ref{tab:notation}.

\begin{table}[!h]
	\centering
	\begin{tabular}{l|l}\hline
		Notation & Description \\\hline
		$\alpha, \beta \in \R^n$ & Degree heterogeneity parameters for incomingness and outgoingness \\
		$\mu \in \R$ & Global sparsity parameter, may diverge to $-\infty$ \\
		$\gamma \in \R^p$ & Covariate parameter, captures homophily, if covariates measure similarity\\
		$\theta \in \R^{2n+1+p}$ & Shorthand for $(\alpha^T, \beta^T, \mu, \gamma^T)^T$\\
		$\theta_0, \alpha_0, \beta_0, \mu_0, \gamma_0$ & The true parameter values \\
		$\vartheta \in \R^{2n}$ & Shorthand for $= (\alpha^T, \beta^T)^T$ \\
		$\xi \in \R^{p+1}$ & Shorthand for $(\mu, \gamma^T)^T$\\
		$S_0$ & The support of $\vartheta_0$ with cardinality $s_0 = \vert S_0 \vert$\\
		$S_{0,+}$ & $= S_0 \cup \{2n+1, 2n+2, \dots 2n+1+p\}$, with cardinality  $s_{0,+} = s_0 + p+1$ \\
		$\lambda$ & The penalty used in \eqref{Eq: Penalized llhd with covariates} \\
		$\bar{\theta} = (\bar \vartheta, \mu, \gamma)$ & $= \left( \frac{{1}}{\sqrt{n}} \vartheta, \mu, \gamma \right)$, the rescaled parameter values \\
		$\bar{\lambda}$ & $=  \sqrt{n}\lambda$, the rescaled penalty value \\
		$b_i, d_i$ & out- and in-degree of node $i$ respectively\\
		$d_+$ & $=\sum_{ i = 1}^n d_i$ \\
		$\mathcal{E}(\theta)$ & $= \frac{1}{N} \E[ \mathcal{L}(\theta) - \mathcal{L}(\theta_0) ]$, the excess risk at $\theta$\\
		$\lambda_{\min}, \lambda_{\max}$ & Minimum and maximum eigenvalue of $\frac{1}{N}\E[Z^T Z]$ respectively \\
		$c_{\min}$ & Constant independent of $n$ such that $\lambda_{\min} \ge c_{\min} > 0$ and $c_{\min} < 1/2$. \\
		\hline
	\end{tabular}
	\caption[Most important definitions]{\label{tab:notation overview} List of most important definitions and notations.}
	\label{tab:notation}
\end{table}

}

\section{Estimation}\label{Sec: Model}

A directed network on $n$ nodes is represented as a directed graph $G_n = (V, E)$, consisting of a node set $V$ with cardinality $n$ and an edge set $E \subseteq V \times V$. Without loss of generality, we assume $V = [n]$ and that $G_n$ is simple, having no self-loops nor multiple edges between any pair of nodes. Such a graph $G_n$ is represented as a binary adjacency matrix $A \in \R^{n \times n}$, where $A_{i,j} = 1$, if $(i,j) \in E$ and $A_{i,j} = 0$ otherwise. By assumption $A_{ii} = 0$ for all $i$.

Given $A$ and the covariates $\{Z_{ij} \}_{i \neq j}$,  the negative log-likelihood of the model \eqref{Eq: srgm} is 
\begin{equation*}
	\begin{split}
		\mathcal{L}(\alpha, \beta, \mu, \gamma) =& - \sum_{i = 1}^n \alpha_i b_i - \sum_{ i = 1}^n \beta_i d_i - d_+ \mu - \sum_{\substack{i,j = 1 \\ i \neq j}}^n (\gamma^TZ_{ij})A_{ij} \\
		& + \sum_{\substack{i,j = 1 \\ i \neq j}}^n \log(1 + \exp(\alpha_i + \beta_j + \mu + \gamma^TZ_{ij})),
	\end{split}
\end{equation*}
where $b_i = \sum_{j = 1, j \neq i}^n A_{ij}$ is the {out-degree} of vertex $i$ and $d_i = \sum_{j=1,  j \neq i }^n A_{ji}$ its in-degree. Write $d = (d_1, \dots, d_n)^T$ and $b = (b_1, \dots, b_n)^T$ as the corresponding degree sequences. Denote $d_+ \coloneqq \sum_{ i = 1}^n d_ i$ and $b_+ \coloneqq \sum_{i = 1}^n b_i$ for which we have $b_+=d_+$. As is common in the literature, we call a network sparse if $\mathbb{E}[d_+]\sim n^\kappa$ for some $\kappa \in (0,2)$, where $\mathbb{E}$ is the expectation with regard to the data generating process. A network is dense if $\mathbb{E}[d_+]\sim n^2$. 

Since $\alpha$ and $\beta$ are sparse, an estimate of $\theta = (\alpha^T, \beta^T, \mu, \gamma^T)^T$ can be obtained via the following penalized likelihood
\begin{equation}\label{Eq: Penalized llhd with covariates0}
	\argmin_{\theta \in \Theta} \frac{1}{N}\mathcal{L}(\alpha, \beta, \mu, \gamma) + \lambda (\Vert \alpha \Vert_1 + \Vert \beta \Vert_1),
\end{equation}
where $\lambda$ is a tuning parameter and $\Theta = \R^n_+ \times \R^n_+ \times \R \times \R^p$ is the parameter space. For simplicity, we have used the same amount of penalty on $\alpha$ and $\beta$ because $b_+=d_+$. The objective function in \eqref{Eq: Penalized llhd with covariates0} is similar to the penalized logistic regression with an $\ell_1$ penalty and thus can be easily solved. In this paper, we use the solver in the \textsf{R} package \texttt{glmnet} \citep{glmnet}. The similarity of our estimator to the LASSO estimator makes our estimation approach extremely scalable.

Since our focus is on sparse networks,  we assume the existence of a non-random sequence $ \rho_{n,0} \in (0, 1/2]$, allowing $\rho_{n,0} \rightarrow 0$ as $n \rightarrow \infty$, such that for all $i, j$, 
\begin{equation*}
	1 - \rho_{n,0} \ge p_{ij} \ge \rho_{n,0},
\end{equation*}
where $\rho_{n,0}$ is referred to as the\textit{ network sparsity parameter}. 
The above constraint is equivalent to 
\begin{equation*}
	\vert  \alpha_{0,i} + \beta_{0,j} + \mu_0 + \gamma_0^TZ_{ij} \vert \le - \text{logit}(\rho_{n,0}) \eqqcolon r_{n,0}, \quad \forall i, j,
\end{equation*}
where $\text{logit}(p) = \log(p / (1-p))$ for $p \in (0,1)$ and $r_{n,0} \ge 0$ since  $\rho_n \le 1/2$.	
This inequality can also be expressed in terms of the design matrix $D$ associated with the corresponding logistic regression problem, defined in \eqref{Eq: Def D} below, and is equivalent to
$
\Vert D\theta_0 \Vert_\infty \le r_{n,0}.
$
This motivates the following tweak to the estimation procedure in \eqref{Eq: Penalized llhd with covariates0}: Given a sufficiently large constant $r_n$,  we define the local parameter space
\begin{equation}\label{Eq: Def of Theta loc}
	\Theta_{\text{loc}} = \Theta_{\text{loc}}(r_n) \coloneqq \left\{ \theta \in \Theta : \Vert D\theta \Vert_\infty \le r_n  \right\},
\end{equation}
which is convex, 
and propose to estimate the parameters as
\begin{equation}\label{Eq: Penalized llhd with covariates}
	\hat{\theta} = (\hat{\alpha}^T, \hat{\beta}^T, \hat{\mu}, \hat{\gamma}^T)^T = \argmin_{\theta = (\alpha^T, \beta^T,\mu, \gamma^T)^T \in \Theta_{\text{loc}}} \frac{1}{N}\mathcal{L}(\alpha, \beta, \mu, \gamma) + \lambda (\Vert \alpha \Vert_1 + \Vert \beta \Vert_1),
\end{equation}
which is more amenable for theoretical analysis.

We now give an explicit form of the associated design matrix $D$. Since we have the presence/ absence of $N = n(n-1)$ directed edges and  $2n + 1 + p$ parameters,  $D$ has dimension $N \times (2n+1+p)$.  Define the {out-matrix} $\Xout \in \R^{N \times n}$ with rows $\Xout_{ij} \in \R^{1 \times n}, i  \neq j$, such that for each component $k = 1, \dots, n$, $\Xout_{ij, k} = 1$ if $k = i$ and zero otherwise. Likewise, define the {in-matrix} $\Xin \in \R^{N \times n}$ with rows $\Xin_{ij} \in \R^{1 \times n}, i  \neq j$, such that for each component $k = 1, \dots, n$, $\Xin_{ij, k} = 1$ if $k = j$ and zero otherwise. Let $Z = (Z_{ij}^T)_{i\neq j} \in \R^{N \times p}$ be the matrix of the covariate vectors written below each other. Then, the design matrix $D$ consists of four blocks, written next to each other:
\begin{equation}\label{Eq: Def D}
	D = \left[
	\begin{array}{c|c|c|c}
		X^{\text{out}} & X^{\text{in}} &\textbf{1} & Z
	\end{array}
	\right] \in \R^{N\times (2n+p+1)},
\end{equation}
where $\textbf{1} \in \R^N$ is a vector of all ones. We use the shorthand
$
X = [\Xout \mid \Xin] \in \R^{N \times 2n}.
$

The design matrix $D$ reveals an important property of our model \eqref{Eq: srgm}. While the columns of the the global parameters $\mu$ and $\gamma$ have non-zero entries in all $N \sim n^2$ rows of $D$, the local parameters $\alpha$ and  $\beta$ only have $n$ non-zero entries in their respective columns. Thus, the effective sample size for $\alpha$ and $\beta$ is only $n$, whereas it is $N$, that is, of order $n$ larger, for the global parameters $\mu$ and $\gamma$. This will also be reflected in the different rates of convergence we obtain in Theorem \ref{Thm: consistency directed} below.

A key quantity in the theory of high-dimensional statistics is the population Gram matrix $\Sigma$, closely linked to the Hessian of $\mathcal{L}$ and the precision matrix. Loosely speaking, it is given by
\[
\Sigma = \frac{1}{\text{sample size}} \E[D^TD].
\]
 Were we to naively ignore the differing sample sizes between local and global parameters and choose $\Sigma = 1/N \cdot \E[D^TD]$, our proofs would fail, due to the top-left corner of $\Sigma$ rapidly converging to the zero-matrix, making $\Sigma$ singular in the limit. In particular, the \textit{compatibility condition} (cf.~Section \ref{Sec: consistency}), crucial for proofs for LASSO-type problems, would not hold. We need to account for this fact and therefore propose to use a \textit{sample-size adjusted} Gram matrix. To that end, we introduce the matrix
\[
T = \begin{bmatrix}
	\sqrt{n-1} I_{2n} & 0 \\
	0 & \sqrt{N} I_{p+1}
\end{bmatrix},
\]
where $I_m$ is the $m \times m$ identity matrix and define the \textit{sample size adjusted Gram matrix} $\Sigma$ as
\begin{equation}\label{Eq: Gram matrix}
	\Sigma = T^{-1} \E[D^TD] T^{-1}.
\end{equation}
It will be convenient to cast problem \eqref{Eq: Penalized llhd with covariates} in terms of rescaled parameters $\bar{\theta}$ which adjust for the discrepancy in effective sample sizes.
This new formulation is equivalent to the one in (\ref{Eq: Penalized llhd with covariates}), but gives us a unified framework for treating convergence properties of our estimators. We will rely heavily on that rescaled version in our proofs. For now, we simply remark that for any parameter $\theta = (\vartheta^T, \mu, \gamma^T)^T \in \Theta$, we introduce the notation
\begin{equation}\label{Eq: bar theta directed}
	\bar{\theta} = (\bar \vartheta^T, \mu, \gamma^T)^T = \left( \frac{{1}}{\sqrt{n}} \vartheta^T, \mu, \gamma^T \right)^T
\end{equation}
and refer the reader to Section \ref{Sec: rescaled problem} in the appendix for a derivation and interpretation of this formulation.
Our original estimation problem \eqref{Eq: Penalized llhd with covariates} can then equivalently be rewritten in terms of these rescaled parameters, giving rise to a sample-size adjusted estimator
\[
\hat{\bar{ \theta}} = \left( \hat{\bar{\vartheta}}^T, \hat{\mu}, \hat{\gamma}^T \right)^T,
\]
that solves a problem similar to \eqref{Eq: Penalized llhd with covariates} with penalty parameter $\bar{\lambda} =  \sqrt{n}\lambda$ (see \eqref{Eq: Penalized llhd with covariates bar} in Section \ref{Sec: rescaled problem}). We denote the negative log-likelihood with respect to a rescaled parameter $\bar{\theta}$ by $\bar{\mathcal{L}}(\bar{\theta})$. Then, given a solution $\hat{ \bar \theta}$ for a given penalty parameter $\bar{\lambda}$ to this modified problem \eqref{Eq: Penalized llhd with covariates bar}, we can obtain a solution to our original problem (\ref{Eq: Penalized llhd with covariates}) with penalty parameter $\lambda = \bar{\lambda}/\sqrt{n}$, by setting
\begin{equation}\label{Eq: penalized llhd bar}
	(\hat{\vartheta}, \hat{\mu}, \hat{\gamma}) = \left( {\sqrt{n}}\hat{\bar{\vartheta}}, \hat{\mu}, \hat{\gamma} \right).
\end{equation}

\section{Theory}\label{Sec: Theory}
We outline the main assumptions first. 

\begin{assumption}\label{Assum: minimum EW directed}\label{Assum: maximum EW directed}
	The $Z_{ij}$'s are independent with $\E[Z_{ij}]=0$ and $|Z_{ij}|$ uniformly bounded. We also assume that  $\gamma_0$ lies in some compact, convex set $\Gamma \subset \R^p$ with a fixed $p$.
	Further assume that there are constants $C > c_{\min} > 0$ such that  the minimum eigenvalue $\lambda_{\min}$ and the maximum eigenvalue $\lambda_{\max}$ of $\frac{1}{N}\E[Z^T Z]$ fulfil $c_{\min} \le \lambda_{\min} \le \lambda_{\max} \le C$.
	Without loss of generality we assume $c_{\min} < 1/2$. 
\end{assumption}
As a result of Assumption \ref{Assum: minimum EW directed}, there exist constants $\kappa, c > 0$ such that $\vert Z_{ij}^T\gamma \vert \le \kappa$ for all $1 \le i \neq j \le n$ and $\vert Z_{ij,k} \vert \le c$ for all $1 \le i \neq j \le n, k = 1, \dots, p$.
\begin{assumption}\label{Assum: no approximation error}
	We assume that $\theta_0 \in \Theta_{\textup{loc}}$ or equivalently $r_n \ge r_{n,0}$. Therefore, without loss of generality we assume $r_n = r_{n,0}$ and consequently $\rho_n = \rho_{n,0}$. 
\end{assumption}

Assumption \ref{Assum: minimum EW directed} is standard.
Note that $Z_{ij}$'s are not necessarily  i.i.d., possibly having correlated entries and that $Z_{ij}$ can be asymmetric in that  $Z_{ij} \neq Z_{ji}$. We have chosen to focus on the random-design assumption which is somewhat more interesting than a fixed-design one. We have assumed a fixed $p$ and leave the study of diverging $p$ to future study. 
Assumption \ref{Assum: no approximation error} ensures no model misspecification.

\begin{assumpB}\label{Assum: model selection assumption}
	$\sqrt{n} s_+^2 \bar{\lambda} \rho_n^{-2} \rightarrow 0, n \rightarrow \infty$.
\end{assumpB}

For all of our theorems, striking the right balance between parameter sparsity $s_+$, network sparsity $\rho_n$ and penalty parameter $\bar{\lambda}$ is crucial. The restrictiveness of these balancing assumptions will depend on the complexity of the results being proven and we number them separately from the general assumptions as ``Assumption B$i$'', $i = 1,2,3$, to make their special standing explicit in our notation.
Our main result on model selection consistency, Theorem \ref{Thm: main result model selection}, is the most refined of our theorems and hence Assumption \ref{Assum: model selection assumption} is the strongest such balancing assumption. In particular, the weaker balancing assumptions required to establish parameter estimation consistency, Theorem \ref{Thm: consistency directed} (Assumption \ref{Assum: rate of s^* directed}), and asymptotic normality of the homophily parameter estimator $\hat{\gamma}$, Theorem \ref{Thm: inference directed} (Assumption \ref{Assum: new rate of s and rho_n directed}), follow from Assumptions \ref{Assum: model selection assumption} above and \ref{Assum: model selection lambda} below. Thus, our estimator $\hat{\theta}$ in \eqref{Eq: Penalized llhd with covariates} can simultaneously recover the correct support, consistently estimate the parameter values and produce an asymptotically normal estimate of $\gamma_0$.

\subsection{Model selection consistency}\label{Sec: model selection}

In this section we study under which conditions our estimator \eqref{Eq: Penalized llhd with covariates} identifies the correct subset of active variables $S_0$. Our main result, Theorem \ref{Thm: main result model selection}, is that under the appropriate conditions, our estimator $\hat{\theta}$ will correctly exclude all the truly inactive parameters and correctly include all those truly active parameters whose value exceeds a certain threshold. The latter minimal signal condition is typical for model selection \citep[e.g.]{Ravikumar:2010, Chen:etal:19}.

Recall that we use $S_0$ to refer to the active set of indices associated with $\vartheta_0 = (\alpha_0^T, \beta_0^T)^T$, whereas $S_{0,+} = S_0 \cup \{2n+1, \dots, 2n+1+p \}$. In the following derivations it will be crucial to distinguish the two correctly. We use $S_{0,+}^c$ to denote the complement of $S_{0,+}$ in $[2n+1+p]$, that is $S_{0,+}^c = [ 2n+1+p ] \backslash S_{0,+}$. Let $S_0^c$ refer to the complement of $S_0$ in $[2n]$ \textit{only}: $S_0^c = [2n] \backslash S_0$. While this may seem like a potential notational pitfall, this allows for much cleaner notation in our proofs.

We first state the main theorem of this section before giving more details on its derivation. Recall that $\bar{\lambda} = \sqrt{n}\lambda$ is the penalty parameter in the rescaled version of our problem \eqref{Eq: Penalized llhd with covariates}. Also notice that $\hat{S} \coloneqq \{i : \hat{\bar{\vartheta}}_i > 0\} =  \{i : \hat{{\vartheta}}_i > 0\}$, that is the estimators \eqref{Eq: Penalized llhd with covariates} and \eqref{Eq: penalized llhd bar} will always select the same active set of parameters.
\begin{assumption}\label{Assum: model selection lambda}
	$- \frac{N \bar{\lambda}^2}{18} + \log(n) \rightarrow - \infty, n \rightarrow \infty$.
\end{assumption}
{\color{black} Assumption \ref{Assum: model selection lambda} suggests we pick our penalty of order $\bar{\lambda} \asymp \sqrt{\frac{\log(n)}{N}}$. Thus, }
informally speaking, Assumption \ref{Assum: model selection lambda} requires that the rescaled penalty parameter $\bar{\lambda}$ must be at least of order $\sqrt{ { \log(\text{number variables} )} / { (\text{effective sample size} ) } }$, which is the typical rate for the penalty we would expect from classical LASSO literature \citep{vandegeer2011}.

\begin{theorem}\label{Thm: main result model selection}\label{Thm: model selection consistency}
	Under Assumptions \ref{Assum: minimum EW directed}, \ref{Assum: no approximation error}, \ref{Assum: model selection assumption} and \ref{Assum: model selection lambda}, and for $n$ sufficiently large, with probability approaching one, the penalized likelihood estimator $\hat{\theta}$ from \eqref{Eq: Penalized llhd with covariates}:
	\begin{enumerate}
		\item excludes all the truly inactive parameters: $\hat{S} \cap S^c = \emptyset$ and,
		\item with penalty of order $\bar{\lambda} \asymp \sqrt{\frac{\log(n)}{N}}$, it includes all those truly active parameters whose value is larger than $C \cdot \rho_n^{-1}\frac{\sqrt{\log(n)}}{\sqrt{n}}$:
		\[
		\left\{i :  \vartheta_{0,i} > C \cdot \rho_n^{-1}\frac{\sqrt{\log(n)}}{\sqrt{n}} \right\} \subseteq \hat{S},
		\]
		where the form of $C$ and the exact probability are given in the proof.
	\end{enumerate}
\end{theorem}

We have the following remarks.

\begin{remark}
	\begin{enumerate}
		\item Notice that Assumption \ref{Assum: model selection lambda} requires $\bar{\lambda} > 3\sqrt{2} \cdot \sqrt{\log(n)/N}$. This is the same regime as specified by Theorems \ref{Thm: consistency directed} and \ref{Thm: inference directed} which ensure consistent parameter estimation and asymptotic normality of $\hat{\gamma}$. Hence, consistent parameter estimation, inference on $\gamma$ and support recovery are all possible simultaneously.
		\item If we choose $\bar{\lambda} \asymp \sqrt{\frac{\log(n)}{N}}$ and $s_+$ is of lower order, such as growing logarithmically or constant, then, up to log-terms, Assumption \ref{Assum: model selection assumption} implies that we must have for the permissible network sparsity, $\rho_n = o(n^{-1/4})$.
	\end{enumerate}
\end{remark}

Our tool of choice for proving Theorem \ref{Thm: main result model selection} is a \textit{primal-dual witness construction}, similar to the one in \cite{Ravikumar:2010}. 
The idea is to construct a tuple $(\bar{\theta}^\dagger, \bar{z}^\dagger)$, such that $\bar{\theta}^\dagger$ solves the rescaled version of \eqref{Eq: Penalized llhd with covariates}, 
while also identifying the correct support $S_0$ and  $\bar{z}^\dagger$ is a solution to the Karush-Kuhn-Tucker conditions \eqref{Eq: 21} as outlined below.
In the construction of $(\bar{\theta}^\dagger, \bar{z}^\dagger)$, we make use of knowledge of the true active set $S_0$, which makes it infeasible to use in practice. However, by Lemma \ref{Lem: unique active set} below, if the construction succeeds -- we make precise what we mean by that below -- any solution to \eqref{Eq: Penalized llhd with covariates} must have the same support as $\bar{\theta}^\dagger$. In summary, if the construction succeeds, our estimator $\hat{{\theta}}$ must identify the correct support $S_0$, too. The bulk of the work in proving Theorem \ref{Thm: model selection consistency} is to show that the construction of $(\bar{\theta}^\dagger, \bar{z}^\dagger)$ will be successful with high probability for large $n$. 

It is important to point out that due to the mixture of deterministic and random columns in $D$ and the differing sample sizes between $\vartheta$ and $\xi$, the standard assumptions in \cite{Ravikumar:2010} imposed on the Hessian of $\mathcal{L}$ cannot simply be imposed in our model. Rather, a careful argument is needed to prove that analogous properties hold for sufficiently large $n$ with high probability. See Section \ref{Sec:auxiliary lemmas} in the appendix for details.

Our starting point for proving Theorem \ref{Thm: model selection consistency} are the Karush-Kuhn-Tucker conditions \citep[Chapter 5]{Bertsekas:1995}: Equation \eqref{Eq: Penalized llhd with covariates} and its rescaled version \eqref{Eq: Penalized llhd with covariates bar}  are a convex optimization problems. Hence, by subdifferential calculus, a vector $\bar{\theta}$ is a minimizer of \eqref{Eq: Penalized llhd with covariates bar} if and only if zero is contained in the subdifferential of $\frac{1}{N} \bar{\mathcal{L}} (\bar{\theta}) + \bar{\lambda} \Vert \bar{\vartheta} \Vert_1$ at $\bar{\theta}$.
That is, if and only if there is a vector $\bar{z} \in \R^{2n+1+p}$ such that
\begin{equation}\label{Eq: 21}
	0 = \frac{1}{N}\nabla \bar{\mathcal{L}}(\bar{\theta}) + \bar{\lambda} \bar{z}, 
\end{equation}
and
\begin{subequations}
	\begin{align}
		\bar{z}_i &= 1, \text{ if } \bar{\vartheta}_i > 0, i = 1, \dots, 2n, \label{Eq: 22a} \\
		\bar{z}_i &\in [-1, 1], \text{ if } \bar{\vartheta}_i = 0, i = 1, \dots, 2n, \label{Eq: 22b} \\
		\bar{z}_i &= 0, i = 2n+1, \dots, 2n+1+p \label{Eq: 22c}.
	\end{align}
\end{subequations}
We call such a pair $(\bar{\theta}, \bar{z}) \in \R^{2n+1+p} \times \R^{2n+1+p}$ \textit{primal-dual optimal} for the rescaled problem \eqref{Eq: Penalized llhd with covariates bar}. Note that in the first $2n$ components of $\nabla \bar{\mathcal{L}}$ we are taking the derivative with respect to $\bar{\vartheta}$ instead of $\vartheta$. This means we need to pay attention to additional $\sqrt{n}$-factors. For such a pair to identify the correct support $S_0$, it is sufficient for
\begin{subequations}
	\begin{align}
		&\bar{\theta}_i > 0, \text{ for all } i \in S_0, \text{ and} \label{Eq: correct inclusion}\\
		&\Vert \bar{z}_{S_{0,+}^c} \Vert_\infty < 1 \label{Eq: correct exclusion}
	\end{align}
\end{subequations}
to hold. Where \eqref{Eq: correct inclusion} ensures that all truly active indices are included and \eqref{Eq: correct exclusion} ensures that all truly inactive indices are excluded (due to \eqref{Eq: 22a}). We call \eqref{Eq: correct exclusion} the \textit{strict feasibility condition} as in \cite{Ravikumar:2010}.

We will proceed to construct a  pair $(\bar{\theta}^\dagger, \bar{z}^\dagger)$ that satisfies condition \eqref{Eq: 21}, \eqref{Eq: 22a} - \eqref{Eq: 22c} and \eqref{Eq: correct inclusion} - \eqref{Eq: correct exclusion} with high probability and for sufficiently large $n$. 
We say the construction \textit{succeeds}, if $(\bar{\theta}^\dagger, \bar{z}^\dagger)$ fulfils \eqref{Eq: 21} - \eqref{Eq: correct exclusion}, which in particular implies that $\bar{\theta}^\dagger$ identifies the correct support $S_0$ and also is a solution to \eqref{Eq: Penalized llhd with covariates bar}. 

By the following lemma, if the construction succeeds, any solution to \eqref{Eq: Penalized llhd with covariates bar} in the appendix must have the same support as $\bar{\theta}^\dagger$. Thus, if the construction succeeds, our estimator $\hat{\bar{\theta}}$ must identify the correct support $S_0$, too.
\begin{lemma}\label{Lem: unique active set}
	Suppose the construction $(\bar{\theta}^\dagger, \bar{z}^\dagger)$ fulfils equations \eqref{Eq: 21} and \eqref{Eq: 22a} - \eqref{Eq: 22c} and \eqref{Eq: correct exclusion}. Let $S^\dagger = \{i: \constrvartheta_i > 0\}$. Then,
	\[
	\hat{S} = S^\dagger.
	\]
	In particular, if $(\bar{\theta}^\dagger, \bar{z}^\dagger)$ additionally fulfils \eqref{Eq: correct inclusion}, then $S^\dagger = S_0$, and thus,
	$
	\hat{S} = S_0.
	$
\end{lemma}
We now give a detailed description of the primal-dual witness construction.

\noindent\textbf{Primal-dual witness construction.}
\begin{enumerate}
	\item Solve the restricted penalized likelihood problem
	\begin{equation}\label{Eq: restricted penalized llhd}
		\bar{\theta}^\dagger = (\bar{\vartheta}^{\dagger,T}, \mu^\dagger, \gamma^{\dagger,T})^T = \argmin \frac{1}{N}\bar{\mathcal{L}}(\bar{\theta}) + \bar{\lambda} \Vert \bar{\vartheta} \Vert_1,
	\end{equation}
	where the argmin is taken over all $\bar{\theta} = (\bar{\vartheta}^T, \mu, \gamma^T)^T \in \Theta_{\text{loc}}$ with support $S_{0,+}$, i.e.~$\bar{\theta}^\dagger_{S_{0,+}} = \bar{\theta}^\dagger$ or equivalently $\bar{\theta}^\dagger_{S_{0,+}^c} = 0$. Thus, by construction, $\bar{\theta}^\dagger$ correctly excludes all inactive indices.
	\item Since \eqref{Eq: restricted penalized llhd} is a convex problem, zero must be contained in its subdifferential at $\bar{\theta}^\dagger$. Thus, we
	set $\bar{z}^\dagger_{i} = 1$, if $\bar{\vartheta}^\dagger_i > 0$ such that \eqref{Eq: 22a} holds and $\bar{z}^\dagger_i = 0, i = 2n+1, \dots, 2n+1+p$, such that \eqref{Eq: 22c} holds.
	By subdifferential calculus we find $\bar{z}^\dagger_{i} \in [-1,1]$, for those $i \in S$ with  $\bar{\vartheta}^\dagger_i = 0$ (in case there are any), such that \eqref{Eq: 21} holds for those components in $S$.
	\item Plug $\bar{\theta}^\dagger$ and $\bar{z}^\dagger$ into \eqref{Eq: 21} and solve for the remaining components of $\bar{z}^\dagger$, such that \eqref{Eq: 21} holds for $(\bar{\theta}^\dagger, \bar{z}^\dagger)$.
\end{enumerate}
The challenge will be proving that \eqref{Eq: correct inclusion} and \eqref{Eq: correct exclusion} also hold, which together ensure that \eqref{Eq: 22b} holds, too.

\subsection{Consistency}\label{Sec: consistency}

In this section we will show that under assumptions similar to those of Theorem \ref{Thm: main result model selection}, our estimator $\hat{\theta}$ will also be consistent in terms of excess risk (cf.~\cite{GreenshteinRitov2004}, \cite{koltchinskii2011}) and $\ell_1$-error. 
To that end, define the excess risk for a parameter $\theta$ as
\begin{equation*}
	\mathcal{E}(\theta) \coloneqq \frac{1}{N} \E[ \mathcal{L}(\theta) - \mathcal{L}(\theta_0) ]. 
\end{equation*}
By construction, 
$
\theta_0 = \argmin_{\theta \in \Theta} \mathcal{E}(\theta) = \argmin_{\theta \in \Theta_{\text{loc}}(r_{n,0})} \mathcal{E}(\theta),
$
where the second equality follows from Assumption \ref{Assum: no approximation error}.

\textbf{A compatibility condition.}
A crucial identifiability assumption in LASSO theory is the so called \textit{compatibility condition} \citep{vandegeer2011, vandegeer2014}. It relates the quantities $\Vert (\hat{\theta} - \theta_0)_{S_{0,+}}\Vert_1$ and
\[
(\hat{\theta} - \theta_0) \Sigma (\hat \theta - \theta_0),
\]
in a suitable sense made precise below and is crucial for deriving consistency results. 
In our model, similar to the sparse $\beta$-model in \cite{Stein:2020}, the classical compatibility condition as for example defined for generalized linear models in \cite{vandegeer2014} does not hold. The reason for this is that $\vartheta$ and $(\mu, \gamma^T)^T$ have different effective sample sizes. Therefore, it is crucial that we use the sample size adjusted Gram matrix \eqref{Eq: Gram matrix}.
Using similar techniques as in \cite{Stein:2020}, we can now show that the sample size adjusted Gram matrix fulfils the compatibility condition.

\begin{proposition}[Compatibility condition]\label{Prop: population compatibility condition}
	Under Assumption \ref{Assum: minimum EW directed}, for $s_0 = o(\sqrt{n})$ and $n$ large enough,
	it holds for every ${\theta} \in \R^{2n+1+p}$ with $\Vert {\theta}_{S^{c}_{0,+}} \Vert_1 \le 3 \Vert {\theta}_{S_{0,+}} \Vert_1$, that
	\[
	\Vert {\theta}_{S_{0,+}} \Vert_1^2 \le \frac{2s_{0,+}}{c_{\min}} {\theta}^T \Sigma { \theta},
	\]
	{\color{black} where $\Sigma$ is the sample size adjusted Gram matrix defined in \eqref{Eq: Gram matrix}.}
\end{proposition}

Parameter estimation consistency is the most lenient of our theorems in terms of restrictions that we have to impose on the parameter sparsity $s_0$ and the network sparsity $\rho_n$. We may replace the stricter assumption \ref{Assum: model selection assumption} by the following.
\begin{assumpB}\label{Assum: rate of s^* directed}
	$\sqrt{n}s_0\rho_n^{-1}\bar{\lambda} \rightarrow 0, n \rightarrow \infty$.
\end{assumpB}

Theorem \ref{Thm: consistency directed} below suggests a choice of $\bar{\lambda} \asymp \sqrt{\log(n) /N}$. Under these conditions, Assumption \ref{Assum: rate of s^* directed} becomes $s_0\rho_n^{-1}\sqrt{\log(n)/n} \rightarrow 0$. That is, up to an additional factor $\rho_n^{-1}$, which is the price we have to pay for allowing our link probabilities to go to zero, the permissible sparsity for $\vartheta_0$ is the permissible sparsity in classical LASSO theory for an effective sample size of order $n$. This makes sense, considering the discussion of the differing effective sample sizes in Section \ref{Sec: Theory}.
Also, this choice of $\bar{\lambda}$ together with Assumption \ref{Assum: rate of s^* directed} implies $s_0 = o(\sqrt{n})$, as is required by Proposition \ref{Prop: population compatibility condition} and which thus is not a restriction.

\begin{theorem}\label{Thm: consistency directed}
	Let Assumptions \ref{Assum: minimum EW directed}, \ref{Assum: no approximation error} and \ref{Assum: rate of s^* directed} hold.
	Fix a confidence level $t$ and let
	\[
	a_n \coloneqq \sqrt{\frac{2\log(2(2n+p+1))}{N}} (1 \vee c).
	\]
	{\color{black}
	Choose $\bar{\lambda} =  \sqrt{n}\lambda$ such that
	\[
		\bar{\lambda} \ge 8 \cdot \left( 8a_n + 2 \sqrt{\frac{t}{N}( 11 (1 \vee (c^2p) ) + 16(1 \vee c) \sqrt{n} a_n  )} + \frac{4t(1 \vee c) \sqrt{n}}{3N} \right)
	\]
}
	Then, with probability at least $1 - \exp(-t)$ we have
	\begin{equation*}
		\mathcal{E}(\hat{\theta}) + \bar{ \lambda} \left( \frac{1}{\sqrt{n}} \Vert \hat{\vartheta} - \vartheta_0  \Vert_1 + \vert \hat{\mu} - \mu_0 \vert + \Vert \hat{\gamma} - \gamma_0 \Vert_1  \right) \le C  \frac{s_{0,+}\bar{ \lambda}^2}{\rho_{n,0}}.
	\end{equation*}
\end{theorem}

{\color{black} Theorem \ref{Thm: consistency directed} implies a lower bound on $\bar{\lambda}$  of order $\sqrt{\log(n)/N}$}, suggesting that we may choose $\bar{\lambda}$ of the same order. Thus, up to the additional factor $\rho_n^{-1}$, we obtain the classical LASSO rates of convergence for a parameter of effective sample size $N$ for $\mu$ and $\gamma$ and those for a parameter of effective sample size $n$ for $\alpha$ and $\beta$. 
If $s_0$ is a lower order term, such as growing logarithmically or constant, then up to log-factors Assumption \ref{Assum: rate of s^* directed} requires that $\rho_n$ tend to zero at rate at most as fast as $1 / \sqrt{n}$, which is faster and thus allows for sparser networks than what we obtained for model selection consistency in Theorem \ref{Thm: model selection consistency}. {\color{black}Under Assumption \ref{Assum: rate of s^* directed}, we can see that  the above average rate of convergence for estimating $\vartheta$ is slower than the rate for estimating $\gamma$. 
The rate of convergence for consistency of the estimator $\hat{{\theta}}$ is of the same order as the rate obtained for the analogous estimator for undirected networks in \cite{Stein:2020}. We refer to \cite{Stein:2020}, Section 2.2, for a comparison between $\ell_1$-penalized estimation procedures as discussed here and $\ell_0$-penalized estimators as advocated by \cite{Chen:etal:19}.}

\subsection{Inference}\label{Sec: inference}

Finally, we derive the limiting distribution of our estimator of the covariates weights, $\hat \gamma$. We will see that the same arguments used for deriving the limiting distribution for $\hat \gamma$ also work for $\hat \mu$ and as a by-product of our proofs we also obtain an analogous limiting result for $\hat \mu$. 

Our strategy will be inverting the  Karush-Kuhn-Tucker conditions, similar to \cite{vandegeer2014}. See Section \ref{Sec:goal and approach} in the appendix for details.
Denote by $H(\hat{\theta}) \coloneqq \left.H_{\xi \times \xi}(\theta)\right|_{\theta = \hat{\theta}}$, the Hessian of $\frac{1}{N}\mathcal{L}(\theta)$ with respect to {\color{black}$\xi = (\mu, \gamma^T)^T$} only, evaluated at $\hat{\theta}$.
\begin{comment}
	Consider the entries of $H(\hat{\theta})$:  For all $k,l = 1, \dots, (p+1)$,
	\begin{align*}
		H(\hat{\theta})_{k,l} = \frac{1}{N}\partial_{\xi_k, \xi_l}\mathcal{L}(\hat{\theta}) =  \frac{1}{N} \sum_{i \neq j} D_{ij,2n+k}D_{ij,2n+l} \cdot  p_{ij}(\hat{\theta})(1-p_{ij}(\hat{\theta})),
	\end{align*}
	where $D_{ij}^T$ is the $(i,j)$-th row of the design matrix $D$, i.e.~in particular $D_{ij,2n+k} = 1$ if $k=1$ and $D_{ij,2n+k} = Z_{ij,k-1}$ for $k=2, \dots, (p+1)$.
\end{comment}
Let $D_\xi = [\textbf{1}|Z]$ be the part of the design matrix $D$, {\color{black} defined in \eqref{Eq: Def D},} corresponding to $\xi$ with rows $D_{\xi, ij}^T = (1, Z_{ij}^T), i \neq j$. Also, let $\hat{W} = \text{diag}\left(\sqrt{p_{ij}(\hat{\theta})(1-p_{ij}(\hat{\theta}))}, i \neq j\right)$. It is then easy to see
\[
H(\hat{\theta}) = \frac{1}{N} D_\xi^T \hat{W}^2 D_\xi.
\]
Let $W_0 = \text{diag}(\sqrt{p_{ij}(\theta_0)(1-p_{ij}(\theta_0))}, i \neq j)$ and consider the corresponding population version: 
\[
\mathbb{E}[H(\theta_0)] = \frac{1}{N} \mathbb{E}[ D_\xi^T W_0^2 D_\xi].
\]
To be consistent with commonly used notation, call $\hat{\Sigma}_\xi = H(\hat \theta) = \frac{1}{N} D_\xi^T \hat{W}^2D_\xi$ and $\Sigma_\xi = \mathbb{E}[H(\theta_0)] = \frac{1}{N} \mathbb{E}[ D_\xi^T W_0^2 D_\xi]$ and
$
\hat \Theta_\xi \coloneqq \hat \Sigma_\xi^{-1},  \Theta_\xi \coloneqq \Sigma_\xi^{-1}.
$

For the proof of asymptotic normality, we need to invert $\hat \Sigma_\xi$ and $\Sigma_\xi$ and show that these inverses are close to each other in an appropriate sense. It is commonly assumed in LASSO theory (cf.~\cite{vandegeer2014}) that the minimum eigenvalues of these matrices stay bounded away from zero. In our case, however, such an assumption is invalid. Since we allow for the lower bound $\rho_{n,0}$ on the link probabilities to go to zero, any lower bound on the entries in $W_0$ will go to zero with $n$ and as a consequence our lower bound on the minimum eigenvalue of $\Sigma_\xi$ will tend to zero as $n$ goes to infinity as well. The best we can achieve is a strict positive definiteness of $\Sigma_\xi$ for finite $n$, but not uniformly in $n$. 
Since these lower bounds tend to zero with increasing $n$, a careful argument is needed and we need to impose a slightly stricter balancing assumption than Assumption \ref{Assum: rate of s^* directed}.

\begin{assumpB}\label{Assum: new rate of s and rho_n directed}
	$\sqrt{n}s_{0} \rho_n^{-2} \bar{\lambda} \rightarrow 0, n \rightarrow \infty$.
\end{assumpB}

\begin{theorem}\label{Thm: inference directed}
	Under Assumptions \ref{Assum: minimum EW directed}, \ref{Assum: no approximation error} and \ref{Assum: new rate of s and rho_n directed}, with $\lambda \asymp \sqrt{\log(n)/N}$ fulfilling the conditions of Theorem \ref{Thm: consistency directed}, we have for any $k = 1, \dots, p$, as $n \rightarrow \infty$,
	\begin{equation*}
		\sqrt{N}\frac{\hat \gamma_k - \gamma_{0,k}}{\sqrt{\hat \Theta_{\vartheta,k+1,k+1}}} \overset{d}{\longrightarrow} \mathcal{N}(0,1).
	\end{equation*}
	We also have for our estimator of the global sparsity parameter, $\hat \mu$, as $n \rightarrow \infty$,
	\begin{equation*}
		\sqrt{N}\frac{\hat \mu - \mu_0}{\sqrt{\hat \Theta_{\vartheta,1,1}}} \overset{d}{\longrightarrow} \mathcal{N}(0,1).
	\end{equation*} 
\end{theorem}
Contrary to what is commonly seen in the penalized likelihood literature \citep{zhang:zhang:2014, vandegeer2014}, no debiasing of $\hat{\gamma}$ and $\hat{\mu}$ is needed. The reason for this is that columns of $D$ pertaining to those parameters which are indeed biased, that is to $\vartheta$, and those pertaining to $\xi = (\mu, \gamma^T)^T$ become asymptotically orthogonal, meaning that the bias in $\hat{\xi}$ vanishes fast enough for the derivation of Theorem \ref{Thm: inference directed} to be possible.
For a lower order $s_0$, Assumption \ref{Assum: new rate of s and rho_n directed} essentially allows for the same level of network sparsity as Assumption \ref{Assum: model selection assumption}, up to lower order factors.

\section{Simulation}\label{Sec: Simulation}

In this section we demonstrate the effectiveness of our estimator \eqref{Eq: Penalized llhd with covariates} in performing simultaneous parameter estimation and model selection consistently. To this end, we test its performance on networks of varying sizes. Specifically, we let $n$ vary between $150$ and $800$ in steps of $50$ and choose the sparsity level $s_0$ to be close to $\sqrt{n}/2$. We let $s_0 = 6,  6, 6,  8,  8, 10, 10, 10, 10, 12, 12, 12, 12, 14$ and chose $s_\alpha = s_\beta = s_0 / 2$ in each case. We selected a heterogeneous configuration for the assignment of non-zero $\alpha$ and $\beta$ values. That is, we included dedicated `spreader' nodes, with large $\alpha$ and zero $\beta$ value as well as `attractor' nodes with large $\beta$ and zero $\alpha$ as well as some nodes with both active $\alpha$ and $\beta$. In detail, we let
\begin{align*}
	\alpha &= (2, 1.5, 1, 0.8, \dots, 0.8, 0, \dots, 0), \\
	\beta &= (0, \dots, 0, 2, 1.5, 1, 0.8, \dots, 0.8, 0, \dots, 0),
\end{align*}
where the number of entries with value $0.8$ was chosen to match the aforementioned sparsity level (zero for the first three values of $n$) and the number of leading zeros in $\beta$ was chosen such that there were exactly two nodes with both active $\alpha$ and $\beta$. We let the networks get progressively sparser and set $\mu = -1.2 \cdot \log(\log(n))$. In all cases we used $p = 2$, sampled the covariate values $Z_{ij,k}, k = 1,2, i \neq j$ from a centred $\text{Beta}(2,2)$ distribution, 
and set $\gamma = (1, 0.8)^T$.  
Our estimator requires us to choose a tuning parameter $\lambda$ and we explored the use of the Bayesian Information Criterion (BIC) as well as a heuristic based on our developed theory for model selection. While the former criterion is purely data-driven, the use of the latter is to ensure that our theoretical results are about right in terms of the rates. Specifically, our BIC is defined as
\[\text{BIC} = 2 \mathcal{L}(\hat \theta (\lambda)) + s(\lambda) \log(N)  \]
where $\hat \theta (\lambda)$ is the solution with $\lambda$ as the tuning parameter, and $s(\lambda)$ is the cardinality of its support.
On the other hand, our heuristic is motivated by the theory developed in the previous sections. We construct {\color{black} $\bar{ \lambda}$ as in Theorem \ref{Thm: consistency directed} based on confidence level $t=3$, choosing to drop the leading factor eight prescribed by Theorem \ref{Thm: consistency directed}.} 
It is known that in high-dimensional settings the penalty values prescribed by mathematical theory in practice tend to over-penalize the parameter values \citep{Yu:etal:2018}. Decreasing the penalty by removing that factor is thus in line with these empirical findings.

We drew $M = 500$ realizations for each value of $n$ and recorded the mean absolute error for estimation of $(\alpha^T, \beta^T)^T$, the absolute error for estimation of $\mu$ and the $\ell_1$-error for estimation of $\gamma$. We also constructed confidence intervals as prescribed by Theorem \ref{Thm: inference directed} and recorded the empirical coverage at the nominal $95\%$ level. Finally, we studied how well BIC and our heuristic did in terms of identifying the correct model.

\begin{figure}[t]
	\centering
	\begin{subfigure}{0.3\textwidth}
		\centering
		\includegraphics[width=1.8in]{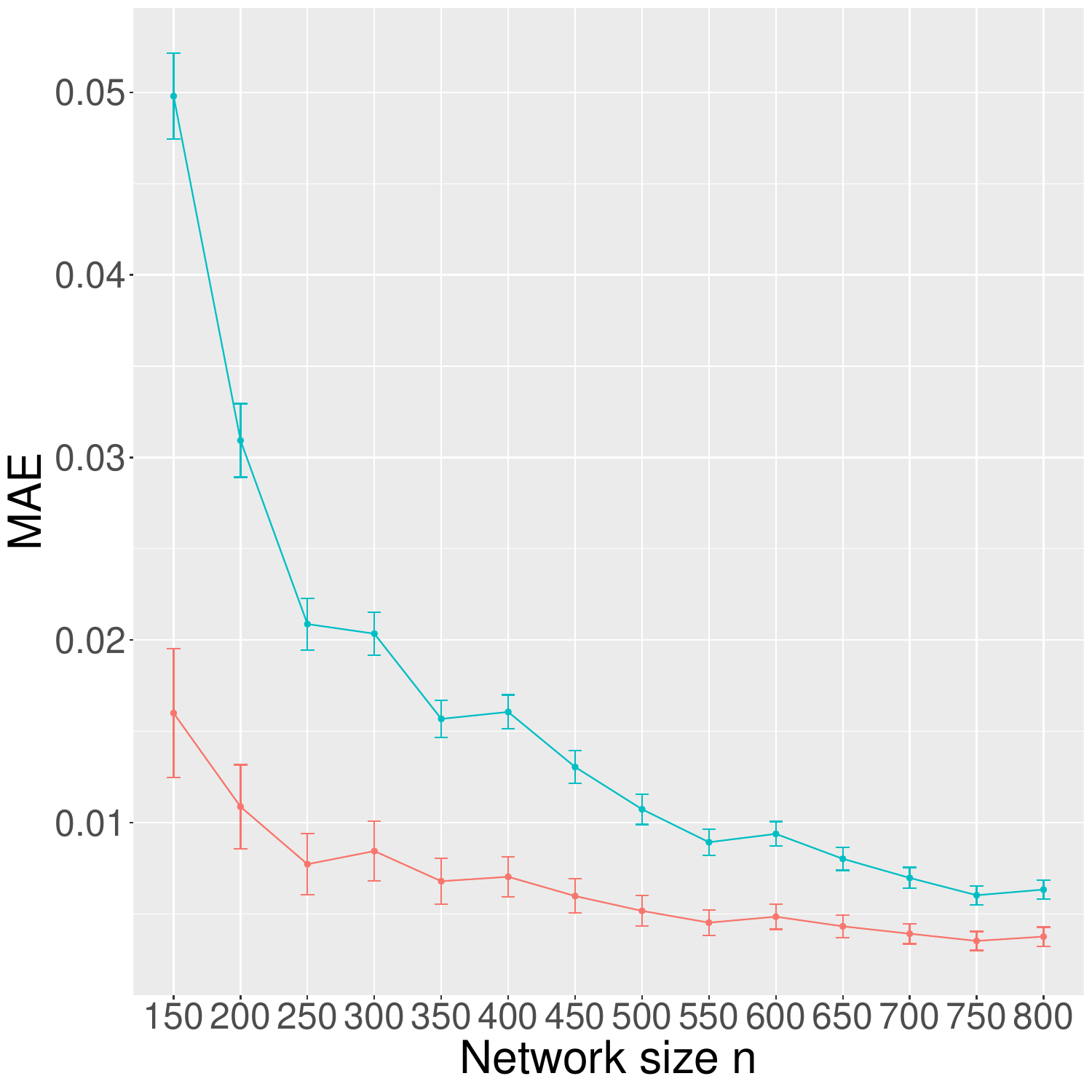}
		\caption{MAE for $\vartheta_0 = (\alpha_0^T, \beta_0^T)^T$.}
		\label{Fig: MAE vartheta}
	\end{subfigure}%
	\begin{subfigure}{0.3\textwidth}
		\centering
		\includegraphics[width=1.8in]{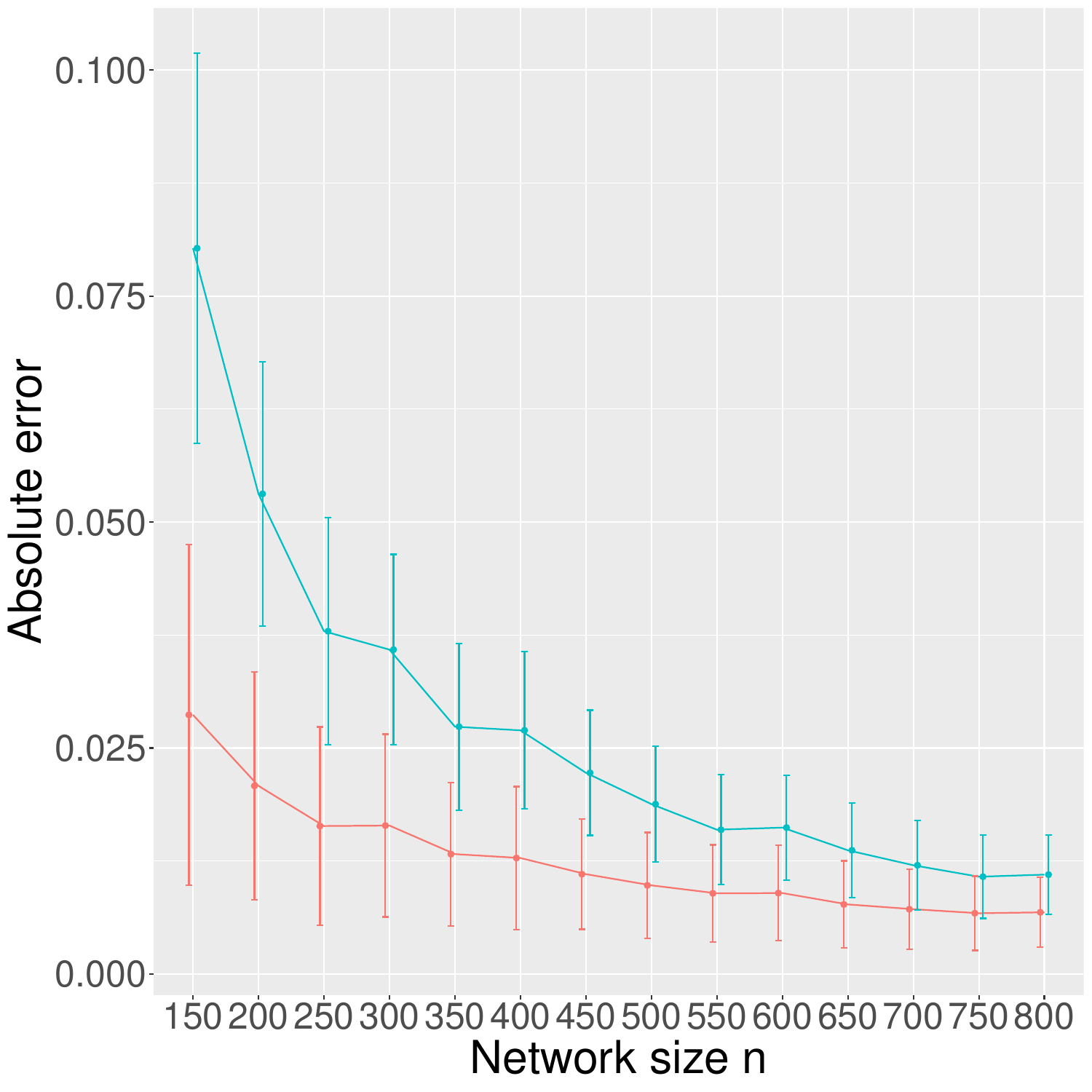}
		\caption{Absolute error for $\mu_0$}
		\label{Fig: mu abs directed}
	\end{subfigure}%
	\begin{subfigure}{0.3\textwidth}
		\centering
		\includegraphics[width=1.8in]{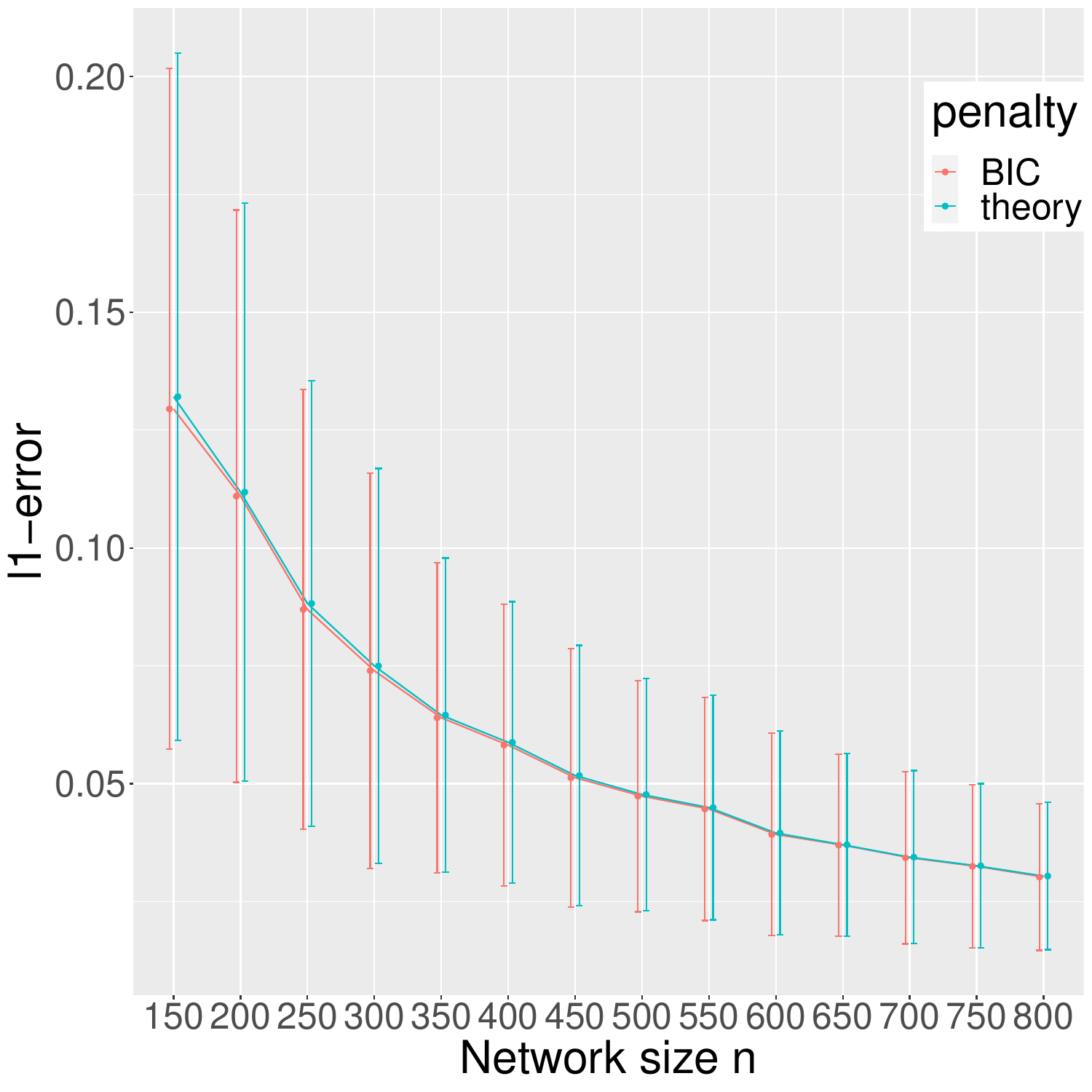}
		\caption{$\ell_1$-error for  $\gamma_0$.}
		\label{Fig: gamma l1 directed}
	\end{subfigure}
	\caption{Mean absolute error for $\vartheta_0 = (\alpha_0^T, \beta_0^T)$, absolute error for $\mu_0$ and $\ell_1$-error for $\gamma_0$ for varying $n$. The results for BIC are presented in red, the ones for our heuristic in green. The dots are the mean errors and the error bars are of length one standard deviation.}
\end{figure}

\noindent\textbf{Consistency.}
We display the various error statistics for estimation of $\vartheta_0 = (\alpha_0^T, \beta_0^T)^T, \mu_0$ and $\gamma_0$ in Figures \ref{Fig: MAE vartheta}, \ref{Fig: mu abs directed} and \ref{Fig: gamma l1 directed} respectively. We see that the error decreases with increasing network size for both model selection procedures. We see that especially for small $n$, BIC outperforms the heuristic for $\vartheta_0$ and $\mu_0$, while they both give essentially the same results for estimation of $\gamma_0$. The better performance of BIC is less prominent as $n$ increases. BIC selects the penalty in a purely data driven manner, which allows it to adapt to differing degrees of sparsity in the network, while for the heuristic the penalty value only depends on $n$ and $p$. This additional flexibility is what allows BIC to achieve lower error values.

\noindent\textbf{Asymptotic normality.}
We construct confidence intervals at the nominal 95\% level for our estimators of $\gamma_{0,1}$ and $\gamma_{0,2}$ as prescribed by Theorem \ref{Thm: inference directed}.
Table \ref{table: CV and CI} shows the results for $\gamma_{0,1}$ across three values of 
$n$. The results for other $n$ and $\gamma_{0,2}$ are similar and are omitted to save space.
The coverage is very close to the $95\%$-level across all network sizes, independent of which model selection criterion we use. This is to be expected, considering that there was hardly any difference for the estimation of $\gamma$ between our two model selection criteria. This empirically illustrates the validity of the asymptotic results derived in Theorem \ref{Thm: inference directed}. As expected, the median length of the confidence interval decreases with increasing network size.

\begin{table}[!htbp]
	\begin{center}
		\begin{tabular}{rrrrrrr}
			&$n$ & Coverage & CI && Coverage & CI\\\hline
			&& \multicolumn{2}{c}{{Heuristic $\lambda$}} && \multicolumn{2}{c}{{BIC}}\\\hline
			&200 &  0.952& 0.265&&0.962 & 0.266\\
			&400 & 0.950& 0.141&&0.964 & 0.141\\
			&800 & 0.946& 0.075&&0.952 & 0.075\\\hline
			
		\end{tabular}
	\end{center}
	\caption{Empirical coverage under nominal 95\% coverage and median lengths of confidence intervals (CIs). }
	\label{table: CV and CI}
\end{table}%

\begin{figure}[t]
	\centering
	\begin{subfigure}{0.45\textwidth}
		\centering
		\includegraphics[scale=0.4]{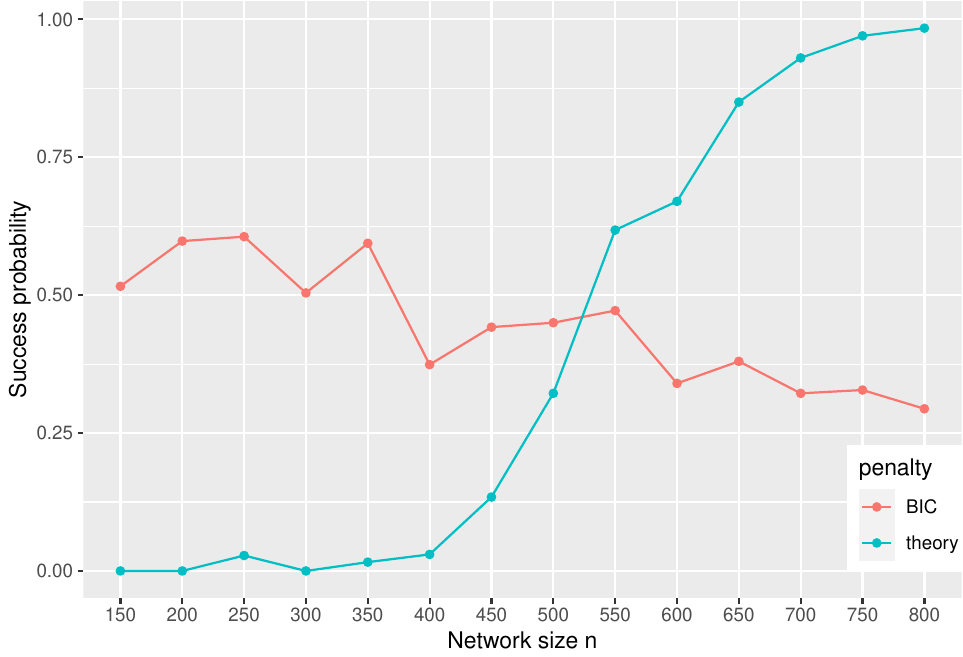}
		\caption{Probability of correct model selection.}
		\label{Fig: model selection}
	\end{subfigure}%
	\begin{subfigure}{0.45\textwidth}
		\centering
		\includegraphics[scale=0.4]{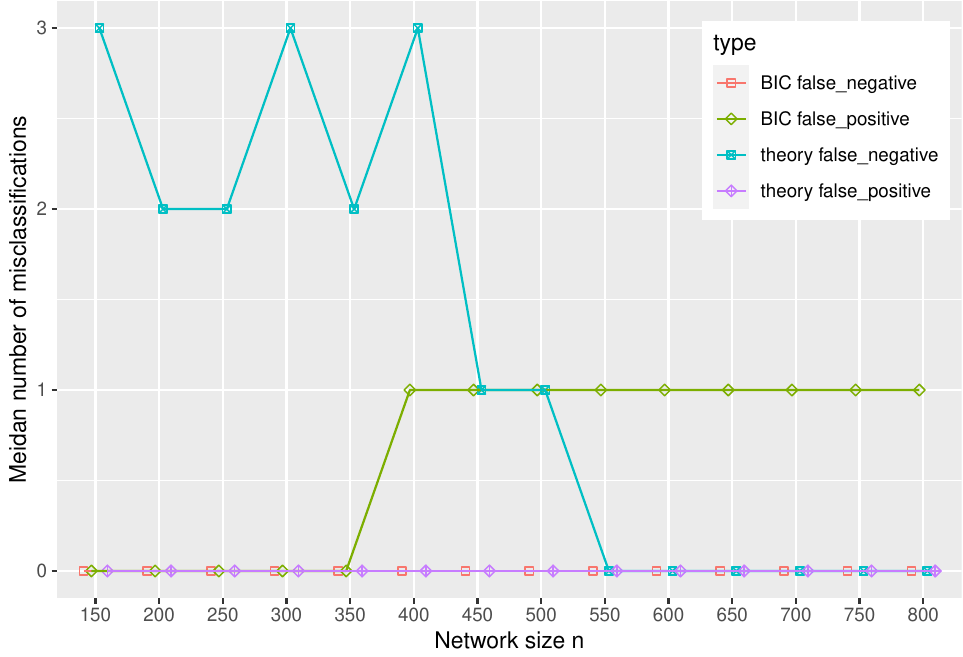}
		\caption{Median number of misclassifications.}
		\label{Fig: misclassifications}
	\end{subfigure}%
	\caption{(a): The empirical probability of selecting the correct subset of active indices. (b): The median number of misclassifications for each model selection procedure, split up into false positives and false negatives.}
\end{figure}

\noindent\textbf{Model selection.} Finally we compare model selection performance between BIC and our heuristic. Figure \ref{Fig: model selection} shows the empirical probability of selecting the correct model versus the various network sizes. We can see very clearly that asymptotically, as $n$ grows, our heuristic outperforms BIC, achieving correct model selection almost all the time.
Nonetheless, it is worth pointing out that even though BIC may not select the exact correct model, the number of misclassifications it does on average is not very large, as shown in Figure \ref{Fig: misclassifications}. Figure \ref{Fig: misclassifications} also shows that the heuristic, by virtue of selecting a larger penalty than BIC, will on average incur more false negatives for small $n$. On the other hand, as $n$ grows, BIC will incur false positives, resulting in the decreasing probability of selecting the correct subset.

\subsection{The lawyer data}
We return to our motivating example by comparing our estimates of the regression coefficients with those in \cite{Yan:etal:2019}. For the seven covariates in this dataset, we followed \cite{Yan:etal:2019} in using the absolute differences of the continuous variables and the indicators whether the categorical variables are equal as our covariates. {\color{black} The edge density in this network is $11.6\%$.}

To apply the model in \cite{Yan:etal:2019}, one needs to remove the eight nodes in black in Figure \ref{fig0} that have zero in-degree or out-degree.  Otherwise the maximum likelihood estimates would be $-\infty$ for $\alpha_i$ if node $i$ has no outgoing connections or for $\beta_i$ if the node has no incoming links. Another interesting aspect of the model in \cite{Yan:etal:2019} lies in the inference for the fixed-dimensional parameter $\gamma$. Because the rate of convergence of its estimate is slowed down by  those of the growing-dimensional heterogeneity parameters $\alpha$ and $\beta$, the estimator of $\gamma$  requires a bias correction to be asymptotically normal. In contrast, by making a sparsity assumption on $\alpha$ and $\beta$ in our model, we estimate the parameters via penalized likelihood and the inference of $\gamma$ is straightforward as seen in Theorem \ref{Thm: inference directed}. 

When the Bayesian information criterion is used to choose the tuning parameter in the penalized likelihood estimation, 
our model gives  7 nonzero $\alpha_i$'s and $7$ nonzero $\beta_i$'s. Four pairs of these nonzeros come from the same nodes. In Table \ref{tab:lawyer results} we present the estimated $\gamma$ and their standard errors when our model and the model in \cite{Yan:etal:2019} are fitted. We remark that since \cite{Yan:etal:2019}  removed eight nodes, akin to biased sampling, their estimates can be biased. In terms of the parameter estimates themselves, although generally similar, we can see a few differences. First we can see that the standard errors of our estimates are smaller than those in \cite{Yan:etal:2019}, reflecting that our estimates are based on a larger sample size (a network with 71 nodes compared to one with 63 nodes in the latter paper) with fewer parameters (22 versus 132). 
Second, the effect of age difference is not significant in our model while it is in the model in \cite{Yan:etal:2019}. To explore the age effect graphically, we colour-coded the lawyers by their age group in Figure \ref{fig0}. We can see that plenty of connections are made between age groups and "across the circle", i.e.~between lawyers with a large difference in age, suggesting that age may not have played an important role. Indeed, a third ($33.9\%$) of all friendships are formed between lawyers with an age difference of ten or more years. Third, we estimate the effect of attending the same law school as positive, implying that the lawyers tend to befriend those who graduated from the same school, while Yan et al.'s model states the opposite. The former conforms better to our intuition about social networks.

\begin{table}[!h]
	\centering
	\begin{tabular}{lrr|rr}\hline
		& \multicolumn{2}{c|}{This Paper} & \multicolumn{2}{c}{\cite{Yan:etal:2019}} \\
		Covariate & Estimate & SE & Estimate & SE \\\hline
		Same status & $1.52$ &$0.10$ &$1.76$&$ 0.16$\\
		Same gender & $0.44$& $0.09$& $0.96$&$ 0.14$\\
		Same office & $2.02$& $0.10$ & $3.23$&$ 0.18$\\
		Same practice & $0.58$&$ 0.09$& $1.11$&$0.12$\\
		Same law school & $0.29 $&$0.10$ &$-0.48$&$ 0.12$\\
		Difference in years with firm & $-0.01$ &$0.006$&$-0.064$&$ 0.014$\\
		Difference in age & $0.003$&$ 0.006$ & $-0.027$&$ 0.011$\\\hline
	\end{tabular}
	\caption[Estimated regression coefficients and standard errors for Lazega's lawyer friendship network (directed)]{\label{tab:lawyer results} Estimated regression coefficients and their standard errors (SE) for Lazega's lawyer friendship network.}
\end{table}

{\color{black}
\subsection{Link prediction for the lawyer data}

In this section we give a brief illustration of how our model can be used for link prediction. We remark that in general it is not possible to conduct vanilla train-test splits or cross-validation on network data, since randomly removing nodes or edges can destroy part of the network structure \citep{li2020network}. While some prior work on network cross-validation exists, notably the aforementioned paper, developing a rigorous cross-validation scheme for our model is beyond the scope of the current paper. Therefore, we present the results in this section as a guideline that the present model shows promising performance for link prediction, even when using vanilla cross-validation. We leave the rigorous mathematical treatment of this interesting problem for a future paper.

In the following, we mimic traditional K-fold cross-validation with K=10 for this dataset as follows. We randomly split the entries of the adjacency matrix $A$ of the lawyer data into 10 folds (ignoring the diagonal). This corresponds to randomly selecting node pairs, as advocated in \cite{li2020network}.
For each fold, we fit our model on the other 90\% of  of the entries of $A$ and use the fitted model to predict the values for the 10\% holdout fold. We run 10 repetitions of this modified 10-fold cross-validation and present the results in Figure \ref{Fig: AUC}. As we can see, even when using this vanilla approach, we achieve a decent performance in the high 80\% range usually, which suggests that with some further adjustments it might be possible to derive strong theoretical guarantees for link prediction for the present model. The key ingredient for why we believe this could be a fruitful research avenue is that link formation happens independently, conditional on the values of the covariates $Z$.

\begin{figure}[h]
	\begin{center}
		\includegraphics[scale=0.35]{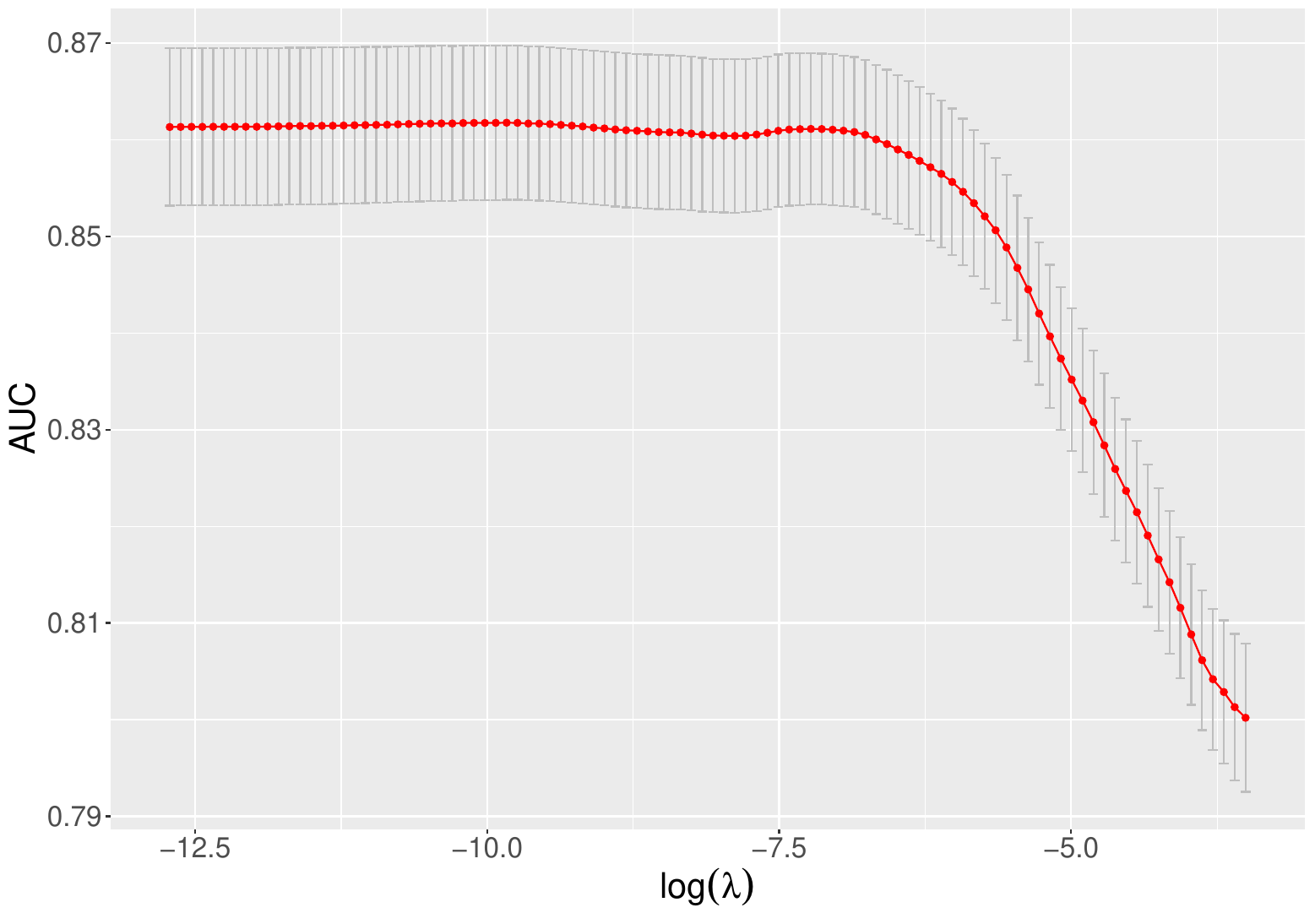}
		\caption{Area under the curve as a function of $\lambda$, averaged over 10 runs of 10-fold cross validations. Error bars are the length of one standard error, averaged over the 10 runs of 10-fold cross-validation.
		}
		\label{Fig: AUC}
	\end{center}
\end{figure}

\subsection{Sina Weibo data}

In our second case study we explore how well our estimation procedure scales to large networks. Towards this, we study the Sina Weibo data collected by \cite{CAI201832}, which was also analysed in \cite{Yan:etal:2019}. Sina Weibo is a Chinese social media platform similar to Twitter. In the original data set there are 4077 nodes representing MBA students and directed links represent who follows whom. Following \cite{Yan:etal:2019}, we focus on the largest strongly connected component consisting of 2242 nodes, leaving us with a very sparse network in which only $0.8\%$ of all possible edges are observed. The resulting in- and out-degree sequences have heavy tails, meaning this network also exhibits a high degree of degree heterogeneity, as illustrated in Figure \ref{Fig: weibo degree distribution}. The in-degrees range from 1 to 253, with the first quartile, median and third quartile equal to 4, 9 and 22 respectively. For the out-degrees we observe values between 1 and 715, with first quartile, median and third quartile equal to 2, 5 and 19 respectively.

\begin{figure}[t]
	\centering
	\begin{subfigure}{0.45\textwidth}
		\centering
		\includegraphics[scale=0.18]{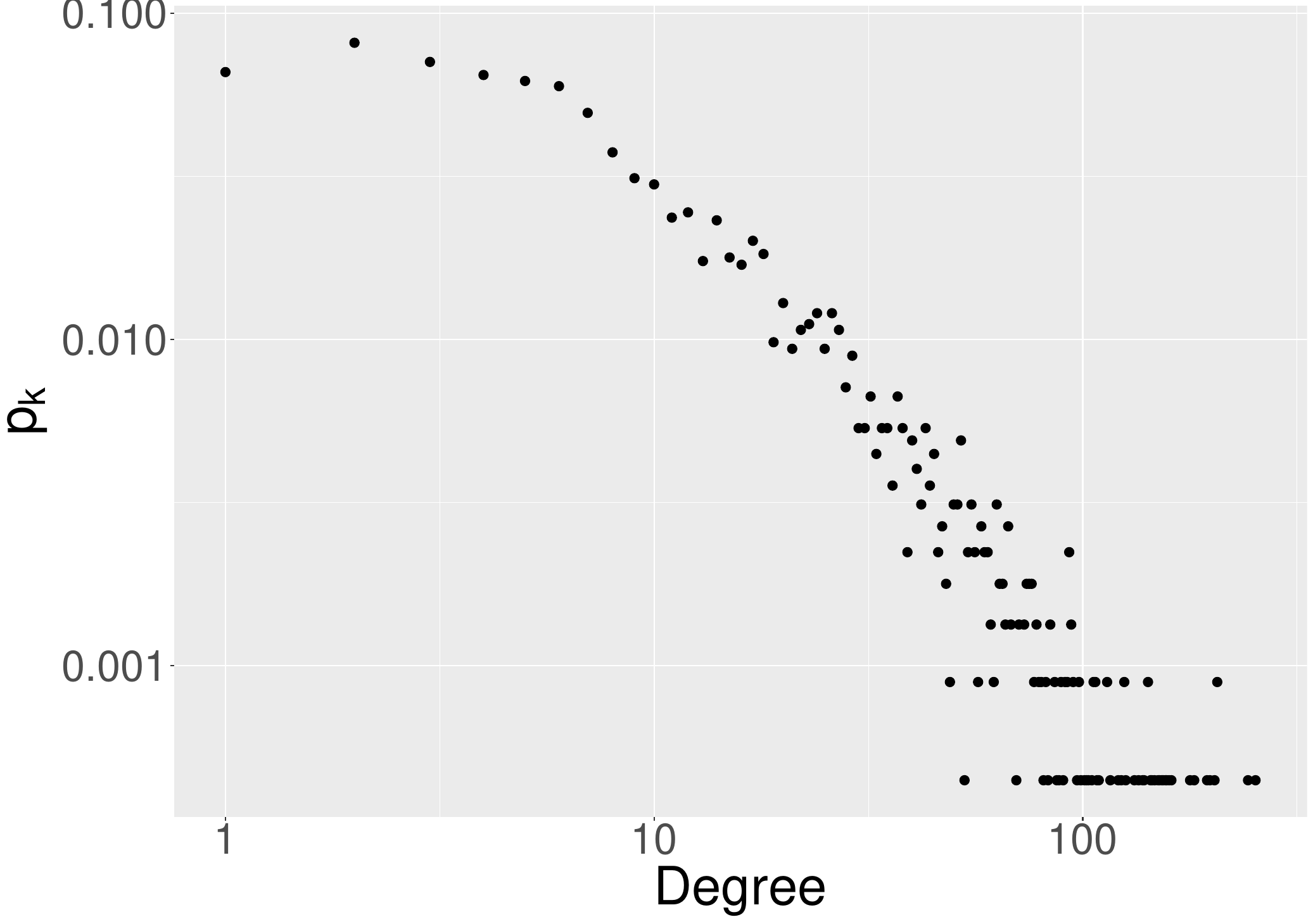}
		\label{Fig: weibo in degree sequence}
	\end{subfigure}%
	\begin{subfigure}{0.45\textwidth}
		\centering
		\includegraphics[scale=0.18]{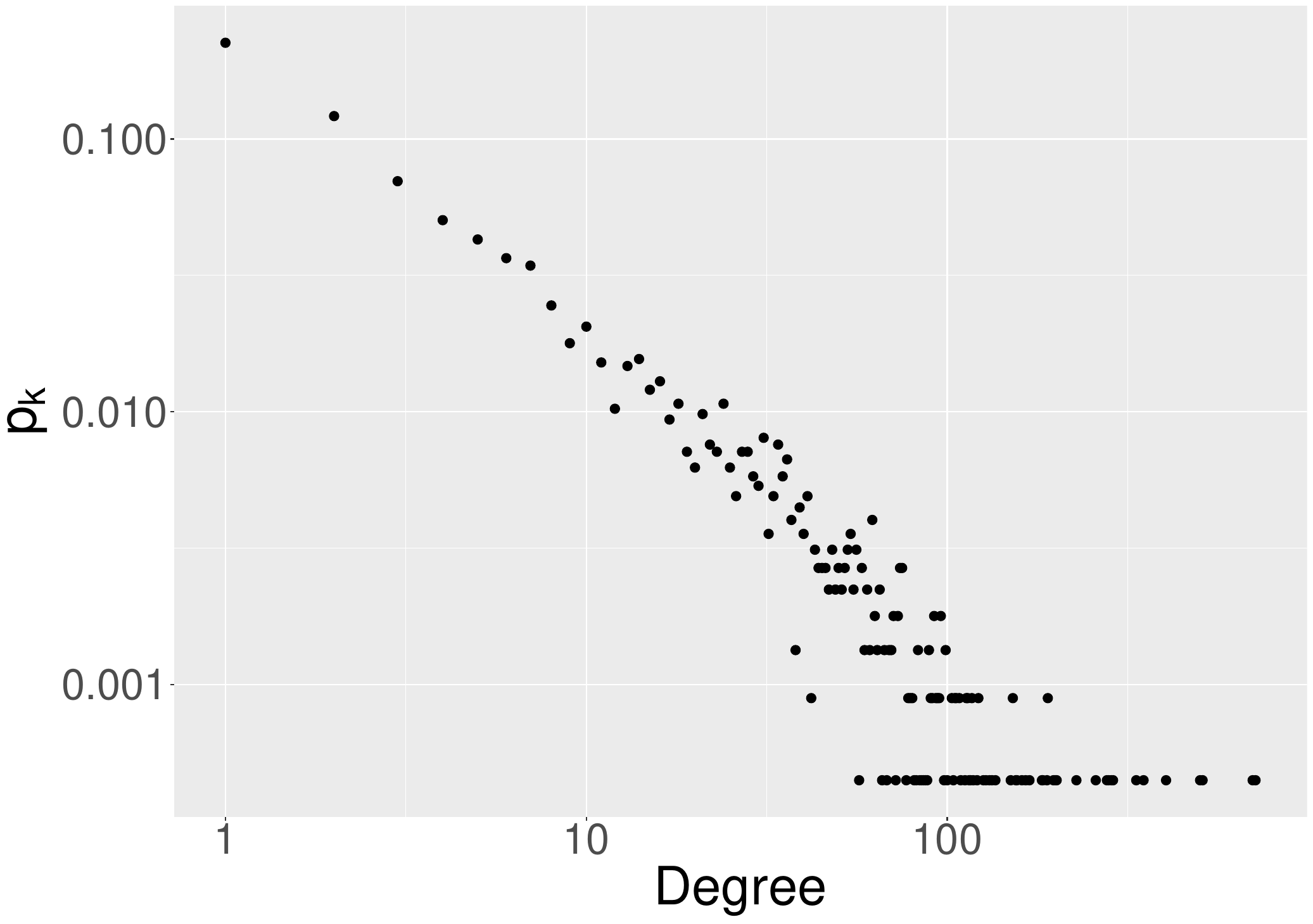}
		\label{Fig: weibo out degree sequence}
	\end{subfigure}%
	\caption{\label{Fig: Weibo degree distribution}(a): The in-degree distribution in the Sina Weibo dataset, log-log scale (b): The out-degree distribution in the Sina Weibo dataset, log-log scale. The y-axis is the empirical frequency with which a give in- / out- degree was observed.}\label{Fig: weibo degree distribution}
\end{figure}

For each node we observe the number of posts they have written, their tenure (measured in months since they joined the platform) and the number of characters in personal labels which were created by the users to describe their lifestyle. We used the absolute difference between these variables as covariates and standardized their values before fitting our model to the data using BIC. In Table \ref{tab:weibo results} we compare the estimates for $\gamma$ and their standard errors when our model is fitted with the results obtained by \cite{Yan:etal:2019}. While both papers estimate all covariates as significant with a negative sign, it is noteworthy that \cite{Yan:etal:2019} initially obtained positive covariate estimates for the difference in the number of posts and labels. Only after running their bias-correction procedure did they obtain the estimates presented here, which illustrates that indeed this ad-hoc procedure is necessary for their model to work properly. On the other hand, we obtain the the negative signs right off-the bat, while also giving us smaller standard errors.

\begin{table}[!h]
	\centering
	\begin{tabular}{lrr|rr}\hline
		& \multicolumn{2}{c|}{This Paper} & \multicolumn{2}{c}{\cite{Yan:etal:2019}} \\
		Covariate & Estimate & SE & Estimate & SE \\\hline
		Difference in posts & $-0.127$ &$0.006$ &$-0.391$&$ 0.018$\\
		Difference in tenure & $-0.050$& $0.005$& $-0.143$&$ 0.008$\\
		Difference in number of labels & $-0.080$& $0.005$&$-0.158$ & $0.008$\\
		\hline
	\end{tabular}
	\caption[Estimated regression coefficients and standard errors for Sina Weibo network]{\label{tab:weibo results} Estimated regression coefficients and their standard errors (SE) for the Sina Weibo network.}
\end{table}

}
\section{Discussion}\label{sec: conc}
We have assumed that links are formed independently between node pairs. This is a limitation because empirically reciprocity, a measure of the likelihood of vertices in a directed network to be mutually linked, may be present. In the motivating lawyers data for example, lawyer $j$ will be more  likely to call lawyer $i$ a friend if the converse is true. To address this layer of sophistication, the next natural step is to add a reciprocity parameter to the model. 
In this paper, we have focused on the inference of the covariate parameter. In some applications inference on $\alpha$ and $\beta$ may be of interest. Since their estimates are biased due to the shrinkage incurred by our $\ell_1$ penalty, this will require debiasing, possibly coupled with suitable balancing assumptions. We leave the exploration of these two interesting research questions for future work.


\acks{We thank the Action Editor and two reviewers for their very helpful comments that have led to a much improved paper.}



\appendix
\section{Appendix A} \label{SM}
We introduce the following additional notation.

For any $a, b \in \R$ we use the notation $a \wedge b = \min \{a, b\}$ and $a \vee b = \max \{a, b\}$.
For any subset $S \subset [n]$, we denote by $v_S$ the vector $v$ with components not belonging to $S$ set to zero. For a matrix $A \in \R^{d \times d}$ and a subset $S \subseteq [d]$, denote by $A_{S,S} \in R^{\vert S \vert \times \vert S \vert}$ the submatrix of $A$ obtained by only taking the rows and columns belonging to $S$. Denote by $A_{-,S} \in \R^{d \times \vert S \vert}$ the submatrix obtained by keeping all the rows and taking only those columns belonging to $S$ and by $A_{S,-} \in \R^{\vert S \vert \times d}$ the submatrix obtained by taking only those rows belonging to $S$ and keeping all the columns. For any square matrix $A$, we denote by $\text{maxeval}(A)$ its maximum eigenvalue and by $\text{mineval}(A)$ its minimum eigenvalue. 

\subsection{A compatibility condition}\label{Sec: proof compatibility condition}
In this section we first prove a \textit{sample compatibility condition} before providing a proof for the population compatibility condition in Proposition \ref{Prop: population compatibility condition}. That is, we first want to find a suitable relation between the quantities $\Vert \hat{\theta} - \theta_0\Vert_1$ and $(\hat{\theta} - \theta_0) \hat{\Sigma} (\hat \theta - \theta_0)$
where $\hat{\Sigma} = T^{-1}D^TDT^{-1}$ is the sample version of the sample size adjusted Gram matrix $\Sigma$.

We make this mathematically precise now:	
For a general matrix $A \in \R^{(2n+1+p) \times (2n+1+p)}$ we say the \textit{compatibility condition holds}, if $A$ has the following property: There is a constant $b$ independent of $n$ such that for every $\theta \in \R^{2n+1+p}$ with $\Vert \theta_{S^{c}_{0,+}} \Vert_1 \le 3 \Vert \theta_{S_{0,+}} \Vert_1$ it holds that
\[
\Vert \theta_{S_{0,+}} \Vert_1^2 \le \frac{s_{0,+}}{b} \theta^T A \theta.
\]

	Notice that the compatibility condition is clearly equivalent to the condition that
\[
\kappa^2(A, s_0) \coloneqq \min_{\substack{\theta \in \R^{2n+1+p} \backslash\{0\} \\ \Vert \theta_{S^{c}_{0,+}} \Vert_1 \le 3 \Vert \theta_{S_{0,+}} \Vert_1}} \frac{\theta^TA\theta}{\frac{1}{s_{0,+}} \Vert \theta_{S_{0,+}} \Vert_1^2}
\]
stays bounded away from zero.

We first show that the compatibility condition holds for the matrix
\begin{equation*}
	\Sigma_A \coloneqq \begin{bmatrix}
		I_{2n}& \textbf{0} & \textbf{0} \\
		\textbf{0} & 1 & \textbf{0} \\
		\textbf{0} & \textbf{0} &  \mathbb{E} [ Z^TZ/N]
	\end{bmatrix} \in \R^{(2n+1+p) \times (2n+1+p)},
\end{equation*}
where $I_{2n}$ is the $(2n) \times (2n)$ identity matrix. 

Recall that by assumption \ref{Assum: minimum EW directed}, the minimum eigenvalue $\lambda_{\min} = \lambda_{\min}(n)$ of $\frac{1}{N} \mathbb{E} [ Z^TZ]$ stays uniformly bounded away from zero. That is, there is a finite constant $c_{\min} > 0$ independent of $n$, such that $\lambda_{\min} > c_{\min} > 0$ for all $n$.
Then, clearly, for any $\theta = (\vartheta^T, \mu, \gamma^T)^T$,
\[
\theta^T \Sigma_A \theta = \Vert \vartheta \Vert_2^2 + \mu^2 + \gamma^T \frac{1}{N} \mathbb{E} [ Z^TZ] \gamma \ge \Vert \vartheta \Vert_2^2 + \mu^2 + c_{\min} \Vert \gamma \Vert_2^2 \ge (1 \wedge c_{\min}) \Vert \theta \Vert_2^2.
\]
Thus, $\Sigma_A$ is strictly positive definite. Furthermore, by Cauchy-Schwarz' inequality, for any $\theta \in \R^{2n+1+p}$ with $\Vert \theta_{S^c_{0,+}} \Vert_1 \le 3 \Vert \theta_{S_{0,+}} \Vert_1$,
\[
\frac{1}{s_{0,+}} \Vert \theta_{S_{0,+}} \Vert_1^2 \le \Vert \theta_{S_{0,+}} \Vert_2^2 \le \Vert \theta \Vert_2^2.
\] 
Thus,
\[
\kappa^2(\Sigma_A, s_0) = \min_{\substack{\theta \in \R^{2n+1+p} \backslash\{0\} \\ \Vert \theta_{S^{c}_{0,+}} \Vert_1 \le 3 \Vert \theta_{S_{0,+}} \Vert_1}} \frac{\theta^T\Sigma_A\theta}{\frac{1}{s_{0,+}} \Vert \theta_{S_{0,+}} \Vert_1^2} \ge \frac{(1 \wedge c_{\min})
	\Vert \theta \Vert_2^2}{\Vert \theta \Vert_2^2} > 0.
\]
We conclude that the compatibility condition holds for $\Sigma_A$.
Now, we need to show that with high probability $\kappa(\hat{\Sigma}, s_0) \ge \kappa(\Sigma_A, s_0)$, which would imply that the compatibility condition holds with high probability for $\hat{\Sigma}$. To that end, we have the following auxiliary lemma found in \cite{kock_tang_2019}. For completeness, we give the short proof of it. The notation is adapted to our setting.

\begin{lemma}[Lemma 6 in \cite{kock_tang_2019}]\label{Lem: Lemma 6 in KockTang}
	Let $A$ and $B$ be two positive semi-definite $(2n+1+p) \times (2n+1+p)$ matrices and $\delta = \max_{ij} \vert A_{ij} - B_{ij} \vert$. For any set $S_0 \subset \{1, \dots, 2n\}$ with cardinality $s_0$, one has
	\[
	\kappa^2(B, s_0) \ge \kappa^2(A, s_0) - 16\delta(s_0 + p + 1).
	\]
\end{lemma}

\begin{proof}
	Denote by $S_{0,+} = S_0 \cup \{2n+1, \dots, 2n+1+p\}$ and $s_{0,+} = s_0 + (1+p)$.
	Let $\theta \in \R^{2n+1+p}\backslash\{0\}$, with $\Vert \theta_{S^{c}_{0,+}} \Vert_1 \le 3\Vert \theta_{S_{0,+}} \Vert_1$. Then,
	\begin{align*}
		\vert \theta^TA\theta - \theta^TB\theta \vert &= \vert \theta^T(A-B) \theta \vert \le \Vert \theta \Vert_1 \Vert (A-B) \theta \Vert_\infty \le \delta \Vert \theta \Vert_1^2 \\
		&= \delta (\Vert \theta_{S_{0,+}} \Vert_1 + \Vert \theta_{S^{c}_{0,+}} \Vert_1)^2 \le \delta (\Vert \theta_{S_{0,+}} \Vert_1 + 3 \Vert \theta_{S_{0,+}} \Vert_1)^2 \\
		&\le 16\delta\Vert \theta_{S_{0,+}} \Vert_1^2.
	\end{align*}
	Hence, $\theta^TB\theta \ge \theta^TA\theta - 16\delta\Vert \theta_{S_{0,+}} \Vert_1^2$ and thus
	\begin{align*}
		\frac{\theta^TB\theta }{\frac{1}{s_{0,+}}\Vert \theta_{S_{0,+}} \Vert_1^2} \ge \frac{\theta^TA\theta }{\frac{1}{s_{0,+}}\Vert \theta_{S_{0,+}} \Vert_1^2} - 16\delta s_{0,+} \ge \kappa^2(A, s_0) - 16\delta s_{0,+}.
	\end{align*}
	Minimizing the left-hand side over all $\theta \neq 0$ with $\Vert \theta_{S^{c}_{0,+}} \Vert_1 \le 3 \Vert \theta_{S_{0,+}} \Vert_1$ proves the claim.
\end{proof}
This shows that to control $\kappa^2(\hat{\Sigma}, s_0)$, we need to control the maximum element-wise distance between $\hat{\Sigma}$ and $\Sigma_A$: $\max_{ij} \vert \hat{\Sigma}_{ij} - \Sigma_{A,ij} \vert$.
Introduce the set
\[
\mathcal{J} = \left\{  \max_{ij} \vert \hat{\Sigma}_{ij} - \Sigma_{A,ij} \vert \le \frac{c_{\text{min}}}{32s_{0,+}}  \right\}.
\]
On the set $\mathcal{J}$, by Lemma \ref{Lem: Lemma 6 in KockTang}, we have 
$
\kappa^2(\hat{\Sigma}, s_0) \ge \kappa(\Sigma_A, s_0) - \frac{c_{\text{min}}}{2} \ge \frac{c_{\text{min}}}{2} > 0
$
and thus the compatibility condition holds for $\hat{\Sigma}$ on $\mathcal{J}$.

\begin{lemma}\label{Lem: Lemma 7 in KockTang}
	If $s_0 = o(\sqrt{n})$, for $n$ large enough, with $\delta = \frac{c_{\text{min}}}{32s_{0,+}}$ and $\tilde{ c} = c^2 \vee (2c^4)$, where $c > 0$ is the universal constant such that $\vert Z_{k,ij} \vert \le c$ for all $k,i,j$, we have
	\begin{align*}
		P(\mathcal{J}) = P\left(  \max_{ij} \vert \hat{\Sigma}_{ij} - \Sigma_{A,ij} \vert \le \frac{c_{\text{min}}}{32s_{0,+}}  \right) &\ge 1 - p(p+3) \exp\left( - N \frac{c_{\min}^2}{2048 s^{2}_{0,+}\tilde{ c}}  \right).
	\end{align*}
\end{lemma}
\begin{proof}
	To make referencing of sections of $\hat{\Sigma}$ easier, we number its blocks as follows
	\begin{equation*}
		\hat{\Sigma} = T^{-1} \begin{bmatrix}
			\underbrace{ X^TX }_{\textcircled{1}} & \underbrace{ X^T\textbf{1} }_{\textcircled{2}} & \underbrace{ X^TZ }_{\textcircled{3}} \\
			\underbrace{\textbf{1}^TX}_{\textcircled{4}} & \underbrace{\textbf{1}^T\textbf{1} }_{\textcircled{5}}& \underbrace{\textbf{1}^TZ}_{\textcircled{6}} \\
			\underbrace{ Z^TX}_{\textcircled{7}} & \underbrace{ Z^T\textbf{1} }_{\textcircled{8}}& \underbrace{Z^TZ}_{\textcircled{9}}
		\end{bmatrix}T^{-1}.
	\end{equation*}
	For block $\textcircled{1}$, i.e.~$i,j = 1, \dots, 2n$, notice that $(\Xout)^T\Xout = (\Xin)^T\Xin = (n-1) I_n$ and $(\Xout)^T\Xin$ is a matrix with zero on the diagonal and ones everywhere else. Therefore, we have either $\hat{\Sigma}_{ij} = \Sigma_{A,ij}$ or
	\[
	\vert \hat{\Sigma}_{ij} - \Sigma_{A,ij} \vert = \frac{1}{n-1} < \frac{c_{\text{min}}}{32s_{0,+}},
	\]
	for $n$ large enough, since $s_{0,+} = o(\sqrt{n})$. Blocks $\textcircled{2}$ and $\textcircled{4}$ are a $2n$ dimensional column and row vector respectively in which each entry is equal to $n-1$. Thus, for $i,j$ corresponding to these blocks,
	\[
	\vert\hat{\Sigma}_{ij} - \Sigma_{A,ij} \vert = \frac{n-1}{\sqrt{(n-1)N}} = \frac{1}{\sqrt{n}} \le  \frac{c_{\text{min}}}{32s_{0,+}},
	\]
	for $n$ large enough, since $s_{0,+} = o(\sqrt{n})$. For $i,j$ corresponding to blocks $\textcircled{3}$ and $\textcircled{7}$, we have
	\[
	\vert \hat{\Sigma}_{ij} - \Sigma_{A,ij} \vert = \frac{c}{\sqrt{n}} < \frac{c_{\text{min}}}{32s_{0,+}},
	\]
	for $n$ large enough. Block $\textcircled{5}$ is a single real number and equal for $\hat{\Sigma}$ and $\Sigma_A$.
	
	The only cases left to consider are those entries corresponding to blocks \textcircled{6}, \textcircled{8} and \textcircled{9}.
	For the blocks \textcircled{6} and \textcircled{8}, that is for $i = 2n+1, j = 2n+2, \dots, 2n+1+p$ and $i = 2n+2, \dots, 2n+1+p, j = 2n+1$, $\hat{\Sigma}_{ij} - \Sigma_{A,ij} = \hat{\Sigma}_{ij}$ is the scaled sum of all the entries of some column $Z_k$ of the matrix $Z$ for an appropriate $k$. That is, there is a $1 \le k \le p$ such that
	\[
	\hat{\Sigma}_{ij} - \Sigma_{A,ij} = \frac{1}{N} Z_k^T\textbf{1} = \frac{1}{N} \sum_{ s \neq t} Z_{k,st}.
	\]
	Note, that thus by model assumption $\mathbb{E}[\hat{\Sigma}_{ij} - \Sigma_{A,ij}] = 0$. We know that for each $k,s,t: Z_{k,st} \in [-c,c]$. Hence, by Hoeffding's inequality, for all $\delta > 0$,
	\[
	P\left( \vert 	\hat{\Sigma}_{ij} - \Sigma_{A,ij}  \vert \ge \delta  \right) = P\left( \left\vert \sum_{ s \neq t} Z_{k,st} \right\vert \ge N \delta  \right) \le 2 \exp\left( - \frac{2N^2\delta^2}{\sum_{i\neq j} (2c)^2}  \right) = 2\exp\left( - N\frac{\delta^2}{2c^2}  \right).
	\]
	For block \textcircled{9}, that is for $i,j = 2n+2, \dots, 2n+1+p$, a typical element has the form
	\[
	\hat{\Sigma}_{ij} - \Sigma_{A,ij} = \frac{1}{N} \sum_{ s \neq t} \left\{ Z_{k,st}Z_{l,st} - \mathbb{E}[Z_{k,st}Z_{l,st}] \right\},
	\]
	for appropriate $k,l$. In other words, $\hat{\Sigma}_{ij} - \Sigma_{A,ij}$ is the inner product of two columns of $Z$, minus their expectation, scaled by $1/N$. Since $Z_{k,st}Z_{l,st} \in [-c^2,c^2]$ for all $k,l,s,t$, we have that for all $k,l,s,t$: $Z_{k,st}Z_{l,st} - \mathbb{E}[Z_{k,st}Z_{l,st}] \in [-2c^2,2c^2]$. Thus, by Hoeffding's inequality, for all $\delta > 0$,
	\[
	P\left( \vert \hat{\Sigma}_{ij} - \Sigma_{A,ij}  \vert \ge \delta  \right) = P\left( \left\vert \sum_{ s \neq t} \{ Z_{k,st}Z_{l,st} - \mathbb{E}[Z_{k,st}Z_{l,st}]\}  \right\vert \ge N \delta \right) \le 2 \exp\left( - N \frac{\delta^2}{8c^4}  \right).
	\]
	Thus, with $\tilde{ c} = c^2 \vee (2c^4)$, we have for any entry in blocks \textcircled{6}, \textcircled{8}, \textcircled{9}, that for any $\delta > 0$,
	\[
	P\left( \vert \hat{\Sigma}_{ij} - \Sigma_{A,ij}  \vert \ge \delta \right) \le 2 \exp\left( - N \frac{\delta^2}{2\tilde{ c}}  \right).
	\]
	Choosing $\delta = \frac{c_{\min}}{32s_{0,+}}$, by the exposition above we know that all entries in blocks \textcircled{1} - \textcircled{5} and \textcircled{7} are bounded by $\delta$ for $n \gg 0$. Also, because block \textcircled{6} is the transpose of block \textcircled{8}, it is sufficient to control one of them. By symmetry of block \textcircled{9} it suffices to control the upper triangular half, including the diagonal, of block \textcircled{9}. Thus, we only need to control the entries $\hat{\Sigma}_{ij} - \Sigma_{A,ij}$ for $i,j$ in the following index set
	\begin{align*}
		\mathcal{A} &= \{ i,j : i,j \text{ belong to block \textcircled{8} or the upper triangular half or the diagonal of block \textcircled{9}} \} \\
		&=	\{ (i,j) \in \{ n+2, \dots, n+1+p\} \times \{n+1\}  \} \cup  \{ i \le j: i,j = n+2, \dots, n+1+p  \}.
	\end{align*}
	Keep in mind that block \textcircled{8} has $p$ elements, while the upper triangular part of block \textcircled{9} plus its diagonal has $\binom{p}{2} + p = \binom{p+1}{2}$ elements.
	Thus, for $n \gg 0$,
	\begin{align*}
		P(\mathcal{J}^c) &= P\left( \max_{ij} \vert \hat{\Sigma}_{ij} - \Sigma_{A,ij} \vert \ge \frac{c_{\text{min}}}{32s_{0,+}}  \right)\\
		& \le \sum_{ i,j \in \mathcal{A}} P\left(  \vert \hat{\Sigma}_{ij} - \Sigma_{A,ij} \vert \ge \frac{c_{\text{min}}}{32s_{0,+}}   \right) \\
		& \le 2p \exp\left( - N\frac{\delta^2}{2c^2}  \right) + 2\binom{p+1}{2}\exp\left( - N \frac{\delta^2}{8c^4}  \right) \\
		&\le 2\left(p+\binom{p+1}{2}\right) \exp\left( - N \frac{\delta^2}{2\tilde{ c}}  \right) \\
		&= p(p+3) \exp\left( - N \frac{\delta^2}{2\tilde{ c}}  \right).
	\end{align*}
	This proves the claim.
\end{proof}

We summarize these results in the following proposition
\begin{proposition}\label{Prop: compatibility condition Sigma}
	Under Assumption \ref{Assum: minimum EW directed}, for $s_0 = o(\sqrt{n})$ and $n$ large enough, with $\tilde{ c} = c^2 \vee (2c^4)$, where $c > 0$ is the universal constant such that $\vert Z_{k,ij} \vert \le c$ for all $k,i,j$: With probability at least
	\[
	1 - p(p+3) \exp\left( - N \frac{c_{\min}^2}{2048 s^{2}_{0,+}\tilde{ c}}  \right)
	\]
	it holds that for every ${\theta} \in \R^{2n+1+p}$ with $\Vert {\theta}_{S^{c}_{0,+}} \Vert_1 \le 3 \Vert {\theta}_{S_{0,+}} \Vert_1$,
	\[
	\Vert {\theta}_{S_{0,+}} \Vert_1^2 \le \frac{2s_{0,+}}{c_{\min}} {\theta}^T \hat{\Sigma} { \theta}.
	\]
\end{proposition}

\begin{proof}
	This follows from Lemma \ref{Lem: Lemma 7 in KockTang}.
\end{proof}

\begin{proof}[Proof of Proposition \ref{Prop: population compatibility condition}]
	To prove that the compatibility condition holds for the population sample-size adjusted Gram matrix $\Sigma$ we may follow the same steps as in the proof of Proposition \ref{Prop: compatibility condition Sigma}: Number the blocks of $\Sigma$ as \textcircled{1} - \textcircled{9} as we did for $\hat{\Sigma}$. $\Sigma$ and $\Sigma_A$ are equal on blocks \textcircled{3}, \textcircled{5}, \textcircled{6}, \textcircled{7}, \textcircled{8} and \textcircled{9}. For blocks \textcircled{1}, \textcircled{2} and \textcircled{4} we use the exact same arguments as in the proof of Proposition \ref{Prop: compatibility condition Sigma} to find that for $n$ sufficiently large, almost surely,
	\[
	\max_{ij} \vert \Sigma_{ij} - \Sigma_{A,ij} \vert \le \frac{c_{\min}}{32s_{0,+}}.
	\]
	The claim follows from Lemma \ref{Lem: Lemma 6 in KockTang}.
\end{proof}

\subsection{A rescaled estimation problem}\label{Sec: rescaled problem}

We now formally introduce the notion of sample-size adjusted parameters $\bar{\theta}$.
Precisely, define the \textit{sample size adjusted design matrix} $\bar D$ as
\[
\bar D = \left[
\begin{array}{c|c|c}
	\bar X &\textbf{1} & Z
\end{array}
\right] \in \R^{N\times (2n+p+1)},
\]
where
\[
\bar X = \left[ \begin{array}{c|c}
	\bar X^{\text{out}} & \bar X^{\text{in}}
\end{array} \right] =  \left[ \begin{array}{c|c}
	\sqrt{n} X^{\text{out}} & \sqrt{n} X^{\text{in}}
\end{array} \right],
\]
is blowing up the entries in $D$ belonging to $\vartheta$.  Recall that for any parameter $\theta = (\vartheta^T, \mu, \gamma^T)^T \in \Theta$, we use
\begin{equation*}
	\bar{\theta} = (\bar \vartheta, \mu, \gamma) = \left( \frac{{1}}{\sqrt{n}} \vartheta, \mu, \gamma \right)
\end{equation*}
to refer to its sample-size adjusted version.
In particular we use the notation $\bar{\theta}_0 = (\bar{\vartheta}_0^T, \mu_0, \gamma_0^T)^T$, to denote the re-parametrized true parameter value.
The blow-up factor $\sqrt{n}$ was chosen precisely such that we can now reformulate our problem as a problem in which each parameter effectively has sample size $N$ in the sense that
\[
\Sigma = \frac{1}{N} \E[\bar{D}^T\bar{D}].
\]
Our original penalized likelihood problem can be rewritten as
\begin{align}\label{Eq: Penalized llhd with covariates bar}
	\begin{split}
		\hat{\bar{ \theta}} = (\hat{\bar{\vartheta}} , \hat{\mu}, \hat{\gamma}) = \argmin_{\substack{ \bar \vartheta= (\bar \alpha^T, \bar \beta^T)^T, \\ \mu, \gamma}}
		&\frac{1}{N} \Bigg( - \sum_{i = 1}^n  \sqrt{n}\bar{\alpha_i} b_i - \sum_{i = 1}^n  \sqrt{n}\bar{\beta_i} d_i - d_+\mu - \sum_{i \neq j} (Z_{ij}^T\gamma)A_{ij} \\
		&+ \sum_{i \neq j} \log\left(1 + \exp\left( \sqrt{n}\bar{\alpha_i} + \sqrt{n}\bar{\beta_j} + \mu + Z_{ij}^T\gamma\right)\right) \Bigg) \\
		&+ \bar{\lambda} \Vert \bar{\vartheta} \Vert_1,
	\end{split}
\end{align}
where $\bar{\lambda} =  \sqrt{n}\lambda$ and the argmin is taken over $\bar{\Theta}_{\text{loc}} = \{\bar{\theta} : \theta \in \Theta, \Vert \bar{ D}\bar{\theta} \Vert_\infty \le r_n \}$. Note that by the same arguments as before, $\bar{\Theta}_{\text{loc}}$ is convex.
Then, given a solution $\hat{ \bar \theta}$ for a given penalty parameter $\bar{\lambda}$ to this modified problem (\ref{Eq: Penalized llhd with covariates bar}), we can obtain a solution to our original problem (\ref{Eq: Penalized llhd with covariates}) with penalty parameter $\lambda = \bar{\lambda}/\sqrt{n}$, by setting
\[
(\hat{\vartheta}, \hat{\mu}, \hat{\gamma}) = \left( {\sqrt{n}}\hat{\bar{\vartheta}}, \hat{\mu}, \hat{\gamma} \right).
\]
Note that for any $\theta \in \Theta$, $D\theta = \bar{ D}\bar{\theta}$, and hence the bound $r_n$ is the same in the definitions of $\Theta_{\text{loc}}$ and $\bar{ \Theta}_{\text{loc}}$. Note also that $\theta \in \Theta_{\text{loc}}$ if and only if $\bar{\theta} \in \bar{ \Theta}_{\text{loc}}$.
For any $\bar{\theta} = (\bar{\vartheta}^T, \mu, \gamma)^T$, denote the negative log-likelihood function corresponding to the rescaled problem (\ref{Eq: Penalized llhd with covariates bar}) as $\bar{\mathcal{L}}(\bar{\theta}) $.
Then, clearly $\bar{\mathcal{L}}(\bar{\theta}) = \mathcal{L}(\theta)$ and $\mathbb{E}[\bar{\mathcal{L}}(\bar{\theta})] = \mathbb{E}[\mathcal{L}(\theta)]$.
Thus, $\bar{\theta}_0$ satisfies that
$	\bar{\theta}_0 = \argmin_{\theta \in \bar{\Theta}} \E[\bar{\mathcal{L}}(\bar{\theta})].$
We define the excess risk for a sample-size adjusted parameter $\bar{\theta}$ as
\begin{equation*}
 \bar{\mathcal{E}}(\bar{\theta}) = \frac{1}{N} \E[ \bar{\mathcal{L}}(\bar{\theta}) - \bar{\mathcal{L}}(\bar{\theta}_0) ]
\end{equation*}
By construction, $\bar{\mathcal{E}} (\bar{\theta}) = \mathcal{E}(\theta)$.

\subsection{A basic Inequality}

A key result in the consistency proofs in classical LASSO settings is the so called \textit{basic inequality} (cf. \cite{vandegeer2011}, Chapter 6). 
Let $P_n$ denote the empirical measure with respect to our observations $(A_{ij}, Z_{ij})$, that is, for any suitable function $g$,
$$
P_ng \coloneqq \frac{1}{N} \sum_{i \neq j} g(A_{ij}, Z_{ij}).
$$
In particular, if we let for each $\theta \in \Theta$, $l_\theta(A_{ij}, Z_{ij}) = -A_{ij} (\alpha_i + \beta_j + \mu + \gamma^TZ_{ij}) + \log(1 + \exp(\alpha_i + \beta_j + \mu + \gamma^TZ_{ij}))$, then
$
P_nl_{\theta} =  \mathcal{L}(\theta)/N.
$
Similarly, we define the theoretical risk as $P = \mathbb{E}P_n$. In particular,
$$
Pl_{\theta} = \mathbb{E}P_n l_\theta = \frac{1}{N} \mathbb{E}[\mathcal{L}(\theta)],
$$
where we suppress the dependence of the theoretical risk on $n$ in our notation. 
We may write the excess risk as
\begin{equation*}
	\mathcal{E}(\theta) \coloneqq P(l_{\theta} - l_{\theta_0}).
\end{equation*}
We define the \textit{empirical process} as
\[
\left\{ v_n(\theta) = (P_n - P)l_{\theta} : \theta \in \Theta \right\}.
\]

\begin{lemma}[Basic Inequality]\label{Lem: basic inequality}
	For any $\theta = (\beta^T, \mu, \gamma^T)^T \in \Theta_{\textup{loc}}$ it holds
	\[
	\mathcal{E}(\hat{\theta}) + \lambda \Vert \hat{\beta} \Vert_1 \le - [v_n(\hat{\theta}) - v_n(\theta)] + \mathcal{E}(\theta) + \lambda \Vert \beta \Vert_1.
	\]
\end{lemma}

\begin{proof}
	By plugging in the definitions and rearranging, we see that the above equation is equivalent to
	\[
	\frac{1}{N} \mathcal{L}(\hat{\theta}) + \lambda \Vert \hat{\beta} \Vert_1
	\le \frac{1}{N} \mathcal{L}(\theta) + \lambda \Vert \beta \Vert_1,
	\]
	which is true by definition of $\hat{\theta}$.
\end{proof}

Notice that since the basic inequality in Lemma \ref{Lem: basic inequality} only relies on the argmin property of the estimator $\hat \theta$, an analogous result follows line by line for the rescaled parameter $\hat{ \bar \theta}$. Writing
\[
\bar{v}_n(\bar{\theta}) \coloneqq \frac{1}{N} (\bar{\mathcal{L}}(\bar{\theta}) - \mathbb{E}[\bar{\mathcal{L}}(\bar{\theta})] ) = v_n(\theta).
\]
for the rescaled empirical process, we have the following.

\begin{lemma}\label{Lem: basic inequality bar directed}
	For any $\bar{\theta} \in \bar{\Theta}_{\textup{loc}}$ it holds
	\[
	\bar{\mathcal{E}}(\hat{\bar{\theta}}) + \bar{\lambda} \Vert \hat{\bar{\vartheta}} \Vert_1 \le - [\bar{v}_n(\hat{\bar{\theta}}) - \bar{v}_n(\bar{\theta})] + \bar{\mathcal{E}}(\bar{\theta}) + \bar{\lambda} \Vert \bar{\vartheta} \Vert_1.
	\]
\end{lemma}

\begin{remark}
	For any $0 < t < 1$ and $\theta \in \Theta_{\text{loc}}$, let $\tilde{\theta} = t\hat{\theta} + (1 - t)\theta$. Since $\Gamma$ is convex, $\tilde{\theta} \in \Theta_{\text{loc}}$ and since $\theta \rightarrow l_{\theta}$ and $\Vert \, . \, \Vert_1$ are convex functions, we can replace $\hat{\theta}$ by $\tilde{\theta}$ in the basic inequality and still obtain the same result. Plugging in the definitions, we see that the basic inequality is equivalent to the following:
	\begin{gather*}
		\mathcal{E}(\tilde{\theta}) + \lambda \Vert \tilde{\beta} \Vert_1 \le - [v_n(\tilde{\theta}) - v_n(\theta)] + \lambda \Vert \beta \Vert_1 + \mathcal{E}(\theta)\\
		\iff \frac{1}{N} \mathcal{L}(\tilde{\theta}) + \lambda \Vert \tilde{\beta} \Vert_1
		\le \frac{1}{N} \mathcal{L}(\theta) + \lambda \Vert \beta \Vert_1
	\end{gather*}
	and by convexity
	\[
	\frac{1}{N} \mathcal{L}(\tilde{\theta}) + \lambda \Vert \tilde{\beta} \Vert_1
	\le  \frac{1}{N} t \mathcal{L}(\hat{\theta}) + \frac{1}{N} (1-t) \mathcal{L}(\theta) + t \lambda \Vert \hat{\beta} \Vert_1 + (1-t) \lambda \Vert \beta \Vert_1 \le \frac{1}{N} \mathcal{L}(\theta) + \lambda \Vert \beta \Vert_1,
	\]
	where the last inequality follows by definition of $\hat{\theta}$. In particular, for any $M > 0$, choosing
	\[
	t = \frac{M}{M + \Vert \hat{\theta} - \theta \Vert_1},
	\]
	gives $\Vert \tilde{\theta} - \theta \Vert_1 \le M$. The completely analogous result holds for $\bar{\theta}$.
\end{remark}

\subsection{Two norms and one function space}\label{Sec: two norms}

To give us a more compact way of writing, for any $\bar{\theta} \in \Theta$ we introduce functions $f_{\bar{\theta}}: \R^{2n+1+p} \rightarrow \R, f_{\bar{\theta}} (v) = v^T\bar{\theta}$ and denote the function space of all such $f_{\bar{\theta}}$ by $\bar{\mathbb{F}} \coloneqq \{ f_{\bar{\theta}} : \bar{\theta} \in \Theta \}$. We endow  $\bar{\mathbb{F}}$ with two norms as follows: 
\newline
Denote the law of the rows of $\bar{D}$ on $\R^{2n+1+p}$,~i.e. the probability measure induced by $(\bar{X}_{ij}^T, 1, Z_{ij}^T)^T, i \neq j$, by $\bar{Q}$. That is, for a measurable set $A = A_1 \times A_2 \subset \R^{2n+1} \times \R^p$, 
\[
\bar{Q}(A) = \frac{1}{N} \sum_{ i \neq j} P( \bar{D}_{ij} \in A) = \frac{1}{N} \sum_{ i \neq j} {\delta}_{ij}(A_1) \cdot P(Z_{ij} \in A_2),
\]
where ${\delta}_{ij}(A_1) = 1$ if $(\bar{X}_{ij}^T,1)^T \in A_1$ and zero otherwise, is the Dirac-measure. We are interested in the $L_2$ and $L_\infty$ norm on $\bar{\mathbb{F}}$ with respect to the measure $\bar{Q}$ on $\R^{2n+1}\times\R^p$. Denote the $L_2(\bar{Q})$-norm of $f \in \bar{\mathbb{F}}$ simply by $\Vert \, . \, \Vert_{\bar{Q}}$ and let $\E_Z$ be the expectation with respect to $Z$:
\[
\Vert f \Vert^2_{\bar{Q}} \coloneqq \Vert f \Vert_{L_2(\bar{Q} )}^2 = \int_{\R^{2n+1}\times\R^p} f(v)^2 \bar{Q}(dv) = \frac{1}{N}\sum_{ i \neq j} \E_Z[ f((\bar{X}_{ij}^T, 1, Z_{ij}^T)^T)^2]
\]
and define the $L_\infty(\bar{Q})$-norm as usual as the $\bar{Q}$-a.s. smallest upper bound of $f$:
\begin{align*}
	\Vert f \Vert_{\bar{Q}, \infty} &= \inf\{ C \ge 0 : \vert f(v) \vert \le C \text{ for } \bar{Q} \text{-almost every } v \in \R^{2n+1+p} \}.
\end{align*}
Notice in particular, that for any $f_{\bar \theta} \in \bar{\mathbb{F}}, \bar{\theta}\in \bar{\Theta}_{\text{loc}}$: $\Vert f_{\bar{\theta }}\Vert_\infty \le \sup_{Z_{ij}}\Vert \bar{D}\bar{\theta} \Vert_\infty \le r_n$.

We make the analogous definitions for the unscaled design matrix. Let $Q$ denote the probability measure induced by the rows of $D$. Since $\bar D\bar \theta = D \theta$, for any $\theta$ with rescaled version $\bar \theta$, we have
\[
\Vert f_{\bar \theta} \Vert_{L_2({ \bar Q} )} = \Vert f_\theta \Vert_{L_2({Q} )}, \quad \Vert f_{\bar \theta} \Vert_{\bar{Q}, \infty} = \Vert f_\theta \Vert_{ Q,\infty}.
\] 
We want to apply the compatibility condition to vectors of the form ${\theta} = {\theta}_1 - {\theta}_2, {\theta}_1, {\theta}_2 \in {\Theta}_{\text{loc}}$. 

Notice, that we have the following relation between the $L_2(Q)$-norm and the sample size adjusted Gram matrix $\Sigma$: For any $\theta$ we have
\begin{equation}\label{Eq: Q bar norm identity}
	\Vert f_{{\theta}} \Vert_{{Q}}^2 = \E_Z \left[\frac{1}{N}\sum_{ i \neq j}  ( D_{ij}^T \theta)^2\right] = \bar \theta^T \Sigma \bar \theta.
\end{equation}
We have the following corollary which follows immediately from Proposition \ref{Prop: population compatibility condition} (see e.g.~\cite{vandegeer2011}, section 6.12 for a general treatment).

\begin{corollary}\label{Lem: compatibility condition}
	Under assumption \ref{Assum: minimum EW directed}, for $s_0 = o(\sqrt{n})$ and $n$ large enough and with $\tilde{c} = c^2 \vee (2c^4)$, where $c > 0$ is the universal constant such that $\vert Z_{k,ij} \vert \le c$ for all $k,i,j$, it holds that for every $\bar{\theta} = \bar{\theta}_1 - \bar{\theta}_2, \bar{\theta}_1, \bar{\theta}_2 \in \bar{\Theta}_{\textup{loc}}$ with $\Vert \bar{\theta}_{S^{c}_{0,+}} \Vert_1 \le 3 \Vert \bar{\theta}_{S_{0,+}} \Vert_1$,
	\[
	\Vert \bar{\theta}_{S_{0,+}} \Vert_1^2 \le\frac{s_{0,+}}{C} \Vert f_{{\theta}_1} - f_{{\theta}_2} \Vert_{{Q}}^2,
	\]
	where $C = c_{\min}/2$.
\end{corollary}

\begin{proof}
	By Proposition \ref{Prop: population compatibility condition},
	\[
	\Vert \bar{\theta}_{S_{0,+}} \Vert_1^2 \le \frac{2s_{0,+}}{c_{\min}} \bar{\theta} \Sigma \bar{\theta}.
	\]
	The claim follows from \eqref{Eq: Q bar norm identity} and the fact that $\theta \mapsto f_{\theta}$ is linear.
\end{proof}

\subsection{Lower quadratic margin for \texorpdfstring{$\mathcal{E}$}{E}}\label{Sec: quadratic margin covariates}

In this section we will derive a lower quadratic bound on the excess risk $\mathcal{E}(\theta)$ if the parameter $\theta$ is close to the truth $\theta_0$. This is a necessary property for the proof to come and is referred to as the \textit{margin condition} in classical LASSO theory (cf. \cite{vandegeer2011}). 

The proof mainly relies on a second order Taylor expansion of the function $l_\theta$ of introduced in section \ref{Sec: Theory}.
Given a fixed $\theta$, we treat $l_\theta$ as a function in $\theta^Tx$ and define new functions $l_{ij}: \R \rightarrow \R, i \neq j,$
\[
l_{ij}(a) = \mathbb{E}[l_\theta(A_{ij}, a) \vert Z_{ij}] = -p_{ij}a + \log(1+\exp(a)),
\]
where $p_{ij} = P(A_{ij} = 1 \vert Z_{ij})$ and by slight abuse of notation we use $l_\theta(A_{ij}, a) \coloneqq -A_{ij} a + \log(1 + \exp(a))$. Taking derivations, it is easy to see that
\[
f_{\theta_0} ((X_{ij}^T, 1, Z_{ij}^T)^T) \in \arg\min_a l_{ij}(a).
\]
All $l_{ij}$ are clearly twice continuously differentiable with derivative
\[
\frac{\partial^2}{\partial a^2}l_{ij}(a) = \frac{\exp(a)}{(1 + \exp(a))^2} > 0, \forall a \in \R.
\]
Using a second order Taylor expansion around $a_0 = f_0((X_{ij}^T, 1, Z_{ij}^T)^T)$ we get
\[
l_{ij}(a) = l_{ij}(a_0) + l'(a_0) (a - a_0) + \frac{l''(\bar{a})}{2}(a - a_0)^2 =  l_{ij}(a_0) + \frac{l''(\bar{a})}{2}(a - a_0)^2,
\]
with an $\bar{a}$ between $a$ and $a_0$. Note that $\vert a_0 \vert \le r_n$. Then, for any $a$ with $\vert a \vert \le r_n$, we must have that for any intermediate point $\bar{a}$ between $a_0$ and $a$ it also holds that $\vert \bar{a} \vert \le r_n$. 
Also note that $\frac{\exp(a)}{(1 + \exp(a))^2}$ is symmetric and monotone decreasing for $a \ge 0$. 
Thus, for any $a$ with $\vert a \vert \le r_n$,
\begin{align}\label{Eq: derivation lower bound}
	\begin{split}
		l_{ij}(a) - l_{ij}(a_0) &= \frac{\exp(\bar{a})}{(1 + \exp(\bar{a}))^2} \frac{(a - a_0)^2}{2} \\
		&= \frac{\exp(\vert\bar{a}\vert)}{(1 + \exp(\vert\bar{a}\vert))^2} \frac{(a - a_0)^2}{2}, \quad \text{by symmetry} \\
		&\ge \frac{\exp( r_n )}{(1 + \exp(r_n))^2} \frac{(a - a_0)^2}{2}.
	\end{split}
\end{align}
In particular, if we pick any $\theta$ and let $a = f_\theta((X_{ij}^T, 1, Z_{ij}^T)^T)$, we have
\begin{align*}
	l_{ij}(f_{\theta} ((X_{ij}^T, 1, Z_{ij}^T)^T)) &- l_{ij}(f_0((X_{ij}^T, 1, Z_{ij}^T)^T)) \\
	&\ge  \frac{\exp(r_{n})}{(1 + \exp( r_{n}))^2} \frac{(f_\theta((X_{ij}^T, 1, Z_{ij}^T)^T) - f_0((X_{ij}^T, 1, Z_{ij}^T)^T))^2}{2}.
\end{align*}
Let
\begin{equation}\label{Eq: Def K_n}
	K_n = \frac{2(1 + \exp(r_{n} ))^2}{\exp(r_{n} )}.
\end{equation}
Define a subset $\mathbb{F}_{\text{local}} \subset \mathbb{F}$ as $\mathbb{F}_{\text{local}} = \{ f_\theta : \theta \in \Theta_{\text{loc}}  \}$. Now, for all $f_\theta \in \mathbb{F}_{\text{local}}$:
\begin{align*}
	\mathcal{E}(\theta)
	&= \frac{1}{N} \sum_{ i \neq j} \mathbb{E} [ l_{\theta}(A_{ij}, D_{ij}) - l_{\theta_0}(A_{ij},  D_{ij}) ] \\
	&= \frac{1}{N} \sum_{ i \neq j} \E[( l_{ij}(f_\theta ( D_{ij}) - l_{ij}(f_0( D_{ij})))] \\
	&\ge \frac{1}{K_n} \cdot \frac{1}{N} (\theta - \theta_0)^T \E_Z[D^TD] (\theta - \theta_0) \\
	&=  \frac{1}{K_n} \cdot \Vert f_\theta - f_0 \Vert_{{Q}}^2.
\end{align*}
Thus, we have obtained a lower bound for the excess risk given by the quadratic function $G_n(\Vert f_{\theta} - f_0 \Vert)$ where $G_n(u) = 1/K_n \cdot u^2$.
Recall that the convex conjugate of a strictly convex function $G$ on $[0, \infty)$ with $G(0) = 0$ is defined as the function
\[
H(v) = \sup_{u} \{ uv - G(u) \}, \quad v > 0,
\]
and in particular, if $G(u) = cu^2$ for a positive constant $c$, we have $H(v) = v^2/(4c)$. 
Hence, the convex conjugate of $G_n$ is
\begin{equation*}
	H_n(v) = \frac{v^2 K_n}{4}.
\end{equation*}
Keep in mind that by definition for any $u,v$
\[
uv \le G(u) + H(v).
\]

\subsection{Consistency on a special set}

In this section we will show that the penalized likelihood estimator is consistent. We will first define a set $\mathcal{I}$ and show that consistency holds on $\mathcal{I}$. It will then suffice to show that the probability of $\mathcal{I}$ tends to one as well. The proof follows in spirit \cite{vandegeer2011}, Theorem 6.4.

We define some objects that we will need for the proof of consistency. We want to use the quadratic margin condition derived in section \ref{Sec: quadratic margin covariates}. 
Recall that the quadratic margin condition holds for any $\theta \in {\Theta}_{\text{loc}}$.
Define
\[
\epsilon^* = H_n\left(   \frac{4\sqrt{2}\sqrt{s_{0,+}}\bar{\lambda}}{\sqrt{c_{\min}}}\right).
\]
Recall the definition of $\bar \theta$ in equation \eqref{Eq: bar theta directed} and let for any $M > 0$
\[
Z_M \coloneqq \sup_{\substack{{\theta} \in {\Theta}_{\text{loc}}, \\  \Vert \bar{\theta} - \bar{\theta}_0 \Vert_1 \le M}} \vert {v}_n({\theta}) - {v}_n({\theta}_0) \vert,
\]
where ${v}_n$ denotes the empirical process. 
The set over which we are maximizing in the definition of $Z_M$ can be expressed in terms of parameters $\theta$ on the original scale as
\[
\left\{ \theta = (\vartheta^T, \mu, \gamma^T)^T \in \Theta_{\text{loc}}: \frac{1}{\sqrt{n}} \Vert \vartheta - \vartheta_0 \Vert_1 + \vert \mu - \mu_0 \vert + \Vert \gamma - \gamma_0 \Vert_1 \le M \right\}.
\]
Set
\[
M^* \coloneqq \epsilon^* / \lambda_0,
\]
where $\lambda_0$ is a lower bound on $\bar{\lambda}$ that will be made precise in the proof showing that $\mathcal{I}$ has large probability. Define
\begin{equation}\label{Def: I covariates}
	\mathcal{I} \coloneqq \{ Z_{M^*} \le \lambda_0 M^* \} = \{ Z_{M^*} \le \epsilon^* \}.
\end{equation}

\begin{theorem}\label{Thm: Consistency covariates}
	Assume that assumptions \ref{Assum: minimum EW directed} and \ref{Assum: rate of s^* directed} hold and that
	$
	\bar{\lambda} \ge 8 \lambda_0.
	$
	Then, on the set $\mathcal{I}$, we have
	\begin{align*}
		\mathcal{E}(\hat{\theta}) + \bar{ \lambda} \left( \frac{{1}}{\sqrt{n}} \Vert \hat{\vartheta} - \vartheta_0  \Vert_1 + \vert \hat{\mu} - \mu_0 \vert + \Vert \hat{\gamma} - \gamma_0 \Vert_1  \right) \le 4\epsilon^*
		= 4 H_n\left(   \frac{4\sqrt{2}\sqrt{s_{0,+}}\bar{ \lambda}}{\sqrt{c_{\min}}}\right).
	\end{align*}
\end{theorem}

\begin{proof}[Proof of Theorem \ref{Thm: Consistency covariates}]
	We assume that we are on the set $\mathcal{I}$ throughout. Set
	\[
	t = \frac{M^*}{M^* + \Vert \hat{\bar{ \theta}} - \bar{ \theta}_0 \Vert_1}
	\]
	and $\tilde{\theta} = (\tilde{\vartheta}^T, \tilde{\mu}, \tilde{\gamma}^T)^T = t \hat{\bar{ \theta}} + (1 - t)\bar{\theta}_0$. Then, 
	\[
	\Vert \tilde{\theta} - \bar{\theta}_0 \Vert_1 = t \Vert \hat{\bar{ \theta}} - \bar{\theta}_0 \Vert \le M^*.
	\]
	Since $\hat{\bar{ \theta}}, \bar{ \theta}_0 \in \bar{\Theta}_{\text{loc}}$ and by the convexity of $\bar{\Theta}_{\text{loc}}$, $\tilde{ \theta} \in \bar{\Theta}_{\text{loc}}$, and by the remark after Lemma \ref{Lem: basic inequality bar directed}, the basic inequality holds for $\tilde{\theta}$. Also, recall that $\bar{ \mathcal{E}}(\bar{ \theta}_0) = 0$:
	\begin{align*}
		\bar{ \mathcal{E}}(\tilde{ \theta}) + \bar{ \lambda} \Vert \tilde{\vartheta} \Vert_1 &\le - (\bar{v}_n(\tilde{ \theta}) - \bar{v}_n(\bar{\theta}_0) ) + \bar{ \mathcal{E}}(\bar{ \theta}_0) + \bar{ \lambda} \Vert \bar{\vartheta}_0 \Vert_1 \\
		&\le Z_{M^*} + \bar{ \lambda} \Vert \bar{\vartheta}_0 \Vert_1 \\
		&\le \epsilon^* + \bar{ \lambda} \Vert \bar{\vartheta}_0 \Vert_1.
	\end{align*}
	From now on write $\tilde{\mathcal{E}} = \bar{\mathcal{E}}(\tilde{\theta})$.
	Note, that $\Vert \tilde{\vartheta} \Vert_1 = \Vert \tilde{\vartheta}_{S^{c}_0} \Vert_1 + \Vert \tilde{\vartheta}_{S_0} \Vert_1$ and thus, by the triangle inequality,
	\begin{align}\label{Eq: 6.29 covariates}
		\begin{split}
			\tilde{\mathcal{E}} + \bar{\lambda} \Vert \tilde{\vartheta}_{S^{c}_0} \Vert_1 &\le \epsilon^* + \bar{ \lambda} (\Vert \bar{\vartheta}_0 \Vert_1 - \Vert \tilde{\vartheta}_{S_0} \Vert_1) \\
			&\le \epsilon^* + \bar{ \lambda} (\Vert \bar{\vartheta}_0 - \tilde{\vartheta}_{S_0} \Vert_1) \\
			&\le \epsilon^* + \bar{ \lambda} (\Vert \bar{\vartheta}_0 - \tilde{\vartheta}_{S_0} \Vert_1 + \Vert (\mu_0, \gamma^{T}_0)^T - (\tilde{\mu}, \tilde{\gamma}^T)^T \Vert_1)  \\
			&= \epsilon^* + \bar{ \lambda} \Vert (\tilde{\theta} - \bar{ \theta}_0)_{S_{0,+}} \Vert_1.
		\end{split}
	\end{align}
	\,
	\newline
	\newline
	\textbf{Case i)} If $\bar{ \lambda} \Vert (\tilde{\theta} - \bar{\theta}_0)_{S_{0,+}} \Vert_1 \ge \epsilon^*$, then
	\begin{equation}\label{Eq: 6.30 covariates}
		\bar{ \lambda} \Vert \tilde{\vartheta}_{S^{c}_0} \Vert_1 \le \tilde{\mathcal{E}} + \bar{ \lambda} \Vert \tilde{\vartheta}_{S^{c}_0} \Vert_1 \le 2\bar{ \lambda} \Vert (\tilde{\theta} - \bar{\theta}_0)_{S_{0,+}} \Vert_1.
	\end{equation}
	Since $\Vert (\tilde{\theta} - \bar{\theta}_0)_{S^{c}_{0,+}} \Vert_1 = \Vert \tilde{\vartheta}_{S^{c}_0} \Vert_1$, we may thus apply the compatibility condition corollary \ref{Lem: compatibility condition} (note that $\bar{\vartheta}_0 = \bar{\vartheta}_{0,S_0}$) to obtain
	\[
	\Vert (\tilde{\theta} - \bar{ \theta}_0)_{S_{0,+}} \Vert_1 \le \sqrt{2} \cdot \frac{\sqrt{s_{0,+}}}{\sqrt{c_{\min}}} \Vert f_{\tilde{ \theta}} - f_{\bar{ \theta}_0} \Vert_{\bar{Q}},
	\]
	where we have used that $\theta \mapsto f_\theta$ is linear and hence $f_{\tilde{ \theta} - \bar{\theta}_0 } = f_{\tilde{ \theta}} - f_{\bar{ \theta}_0}$.
	Observe that
	\begin{equation}\label{Eq: theta norm decomposition}
		\Vert \tilde{\theta} - \theta_0\Vert_1 = \Vert \tilde{\vartheta}_{S^{c}_0} \Vert_1 + \Vert (\tilde{\theta} - \theta_0)_{S_{0,+}} \Vert_1.
	\end{equation}
	Hence, 
	\begin{align*}
		\tilde{\mathcal{E}} + \bar{ \lambda} \Vert \tilde{\theta} - \bar{ \theta}_0 \Vert_1  &= \tilde{\mathcal{E}} + \bar{ \lambda} ( \Vert \tilde{\vartheta}_{S^{c}_0} \Vert_1 + \Vert (\tilde{\theta} - \bar{\theta}_0)_{S_{0,+}} \Vert_1)\\
		&\le \epsilon^* + 2\bar{ \lambda} \Vert (\tilde{\theta} - \bar{\theta}_0)_{S_{0,+}} \Vert_1 \\
		&\le \epsilon^* + 2\sqrt{2} \bar{ \lambda}\frac{\sqrt{s_{0,+}}}{\sqrt{c_{\min}}} \Vert f_{\tilde{\theta}} - f_{\bar{ \theta}_0} \Vert_{\bar{Q}}.
	\end{align*}
	Recall that for a convex function $G$ and its convex conjugate $H$ we have $uv \le G(u) + H(v)$. Thus, we obtain
	\begin{align*}
		2\sqrt{2}\bar{\lambda}\frac{\sqrt{s_{0,+}}}{\sqrt{c_{\min}}} \Vert f_{\tilde{\theta}} - f_{\bar{ \theta}_0} \Vert_{\bar{Q}} &= 4\sqrt{2} \bar{\lambda} \frac{\sqrt{s_{0,+}}}{\sqrt{c_{\min}}} \frac{\Vert f_{\tilde{\theta}} - f_{\bar{ \theta}_0} \Vert_{\bar{Q}}}{2} \\
		&\le H_n\left(4\sqrt{2} \bar{\lambda} \frac{\sqrt{s_{0,+}}}{\sqrt{c_{\min}}}\right) + G_n\left( \frac{\Vert f_{\tilde{ \theta}} - f_{\bar{ \theta}_0} \Vert_{\bar{Q}} }{2}\right) \\
		&\overset{G_n \text{ convex}}{\le}H_n\left(4\sqrt{2} \bar{ \lambda} \frac{\sqrt{s_{0,+}}}{\sqrt{c_{\min}}}\right)+ \frac{G_n(\Vert f_{\tilde{\theta}} - f_{\bar{ \theta}_0} \Vert_{\bar{Q}})}{2}  \\
		&\overset{\text{margin condition}}{\le} H_n\left(4\sqrt{2} \bar{ \lambda} \frac{\sqrt{s_{0,+}}}{\sqrt{c_{\min}}}\right)+ \frac{\tilde{\mathcal{E}}}{2}.
	\end{align*}
	It follows
	\begin{equation*}
		\tilde{\mathcal{E}} + \bar{ \lambda} \Vert\tilde{\theta} - \bar{\theta}_0\Vert_1 \le \epsilon^* + H_n\left(4\sqrt{2} \bar{\lambda} \frac{\sqrt{s_{0,+}}}{\sqrt{c_{\min}}}\right)+ \frac{\tilde{\mathcal{E}}}{2} = 2 \epsilon^* + \frac{\tilde{\mathcal{E}}}{2}
	\end{equation*}
	and therefore
	\begin{equation}\label{Eq: 6.31 covariates}
		\frac{\tilde{\mathcal{E}}}{2} +\bar{\lambda} \Vert\tilde{\theta} - \bar{\theta}_0\Vert_1 \le 2 \epsilon^*.
	\end{equation}
	Finally, this gives
	\[
	\Vert\tilde{\theta} - \bar{\theta}_0\Vert_1 \le \frac{2 \epsilon^*}{\bar{\lambda}} = \frac{2 \lambda_0 M^*}{\bar{\lambda}} \underbrace{\le}_{\bar{ \lambda} \ge 4 \lambda_0} \frac{M^*}{2}.
	\]
	From this, by using the definition of $\tilde{\theta}$, we obtain
	\begin{align*}
		\Vert \tilde{\theta} - \bar{\theta}_0 \Vert_1 = t \Vert \hat{\bar{\theta}} - \bar{\theta}_0 \Vert_1 = \frac{M^*}{M^* + \Vert \hat{\bar{\theta}} - \bar{\theta}_0 \Vert_1} \Vert \hat{\bar{\theta}} - \bar{\theta}_0 \Vert_1 \le \frac{M^*}{2}.
	\end{align*}
	Rearranging gives
	\[
	\Vert \hat{\bar{\theta}} - \bar{\theta}_0 \Vert_1 \le M^*.
	\]
	\,
	\newline
	\textbf{Case ii)} If $\bar{\lambda} \Vert (\bar{\theta}_0 - \tilde{\theta})_{S_{0,+}} \Vert_1 \le \epsilon^*$, then from (\ref{Eq: 6.29 covariates})
	\[
	\tilde{\mathcal{E}} + \bar{\lambda} \Vert \tilde{\vartheta}_{S^{c}_0} \Vert_1 \le 2\epsilon^*.
	\]
	Using once more (\ref{Eq: theta norm decomposition}), we get
	\begin{equation}\label{Eq: 6.32 covariates}
		\tilde{\mathcal{E}} + \bar{\lambda} \Vert \tilde{\theta} - \bar{\theta}_0 \Vert_1 = \tilde{\mathcal{E}} + \bar{\lambda} \Vert \tilde{\vartheta}_{S^{c}_0} \Vert_1 + \bar{\lambda}\Vert (\tilde{\theta} - \bar{\theta}_0)_{S_{0,+}} \Vert_1 \le 3\epsilon^*.
	\end{equation}
	Thus,
	\[
	\Vert \tilde{\theta} - \bar{\theta}_0 \Vert_1 \le 3 \frac{\epsilon^*}{\bar{\lambda}} = 3 \frac{\lambda_0}{\bar{\lambda}}M^* \le \frac{M^*}{2}
	\]
	by choice of $\lambda \ge 6 \lambda_0$. Again, plugging in the definition of $\tilde{\theta}$, we obtain
	\[
	\Vert \hat{\bar{\theta}} - \bar{\theta}_0 \Vert_1 \le M^*.
	\]
	Hence, in either case we have $\Vert \hat{\bar{\theta}} - \bar{\theta}_0 \Vert_1 \le M^*$. That means, we can repeat the above steps with $\hat{\bar{\theta}}$ instead of $\tilde{ \theta}$: Writing $\hat{\mathcal{E}} \coloneqq \bar{\mathcal{E}}(\hat{\bar{\theta}})$, following the same reasoning as above we arrive once more at (\ref{Eq: 6.29 covariates}):
	\begin{equation*}
		\hat{\mathcal{E}} + \bar{\lambda} \Vert \hat{\bar{\vartheta}}_{S^{c}_0} \Vert_1 \le \epsilon^* + \bar{\lambda} \Vert \bar{\vartheta}^* - \hat{\bar{\vartheta}}_{S_0} \Vert_1 \le 2\epsilon^* + \bar{\lambda}\Vert (\hat{\bar{\theta}} - \bar{\theta}_0)_{S_{0,+}} \Vert_1.
	\end{equation*}
	From this, in \textbf{case i)} we obtain (\ref{Eq: 6.30 covariates}) which allows us to use the compatibility assumption to arrive at (\ref{Eq: 6.31 covariates}):
	\[
	\frac{\hat{\mathcal{E}}}{2} + \bar{\lambda} \Vert \hat{\bar{\theta}} - \bar{\theta}_0 \Vert_1 \le 2 \epsilon^*,
	\]
	resulting in
	\[
	\hat{\mathcal{E}} + \bar{\lambda} \Vert \hat{\bar{\theta}} - \bar{\theta}_0 \Vert_1 \le 4 \epsilon^*.
	\]
	In \textbf{case ii)} on the other hand, we arrive directly at (\ref{Eq: 6.32 covariates}), and hence
	\[
	\hat{\mathcal{E}} + \bar{\lambda} \Vert \hat{\bar{\theta}} - \bar{\theta}_0 \Vert_1 \le 3 \epsilon^*.
	\]
	Plugging in the definitions of $\hat{\bar{\theta}}$ and $\bar{\theta}_0$ and using the fact that $\hat{\mathcal{E}} = \bar{\mathcal{E}}(\hat{\bar{\theta}}) = \mathcal{E}(\hat{\theta})$ proves the claim.
\end{proof}

\subsection{Controlling the special set \texorpdfstring{$\mathcal{I}$}{I}}\label{Sec: controlling I directed}

We now show that $\mathcal{I}$ has probability tending to one. Recall some results on concentration inequalities.

\subsection*{Concentration inequalities}
We first recall some probability inequalities that we will need. This is based on chapter 14 in \cite{vandegeer2011}. Throughout let $Z_1, \dots, Z_n$ be a sequence of independent random variables in some space $\mathcal{Z}$ and $\mathcal{G}$ be a class of real valued functions on $\mathcal{Z}$.

\begin{definition}
	A \textit{Rademacher sequence} is a sequence $\epsilon_1, \dots, \epsilon_n$ of i.i.d. random variables with $P(\epsilon_i = 1) = P(\epsilon_i = -1) = 1/2$ for all $i$.
\end{definition}
\begin{theorem}[Symmetrization Theorem as in \cite{vandervaartwellner1996}, abridged]\label{Thm: symmetrization theorem}
	Let $\epsilon_1, \dots, \epsilon_n$ be a Rademacher sequence independent of $Z_1, \dots, Z_n$. Then
	\[
	\mathbb{E}\left( \sup_{g \in \mathcal{G}} \left| \sum_{ i = 1}^n \{ g(Z_i) - \mathbb{E}[g(Z_i)] \} \right| \right) \le 2 \mathbb{E}\left( \sup_{g \in \mathcal{G}} \left| \sum_{ i = 1}^n \epsilon_i g(Z_i) \right|  \right).
	\]
\end{theorem}

\begin{theorem}[Contraction Theorem as in \cite{Ledoux:Talagrand:1991}]\label{Thm: contraction theorem Ledoux Talagrand}
	Let $z_1, \dots, z_n$ be non-random elements of $\mathcal{Z}$ and let $\mathcal{F}$ be a class of real-valued functions on $\mathcal{Z}$. Consider Lipschitz functions $g_i: \R \rightarrow \R$ with Lipschitz constant $L =1$, i.e. for all $i$
	\[
	\vert g_i(s) - g_i(s') \vert \le \vert s - s' \vert, \forall s,s' \in \R.
	\]
	Let $\epsilon_1, \dots, \epsilon_n$ be a Rademacher sequence. Then for any function $f^*: \mathcal{Z} \rightarrow \R$ we have
	\[
	\mathbb{E}\left( \sup_{f \in \mathcal{F}} \left| \sum_{ i = 1}^n \epsilon_i \{ g_i(f(z_i)) - g_i(f^*(z_i)) \} \right| \right) \le 2 \mathbb{E}\left( \sup_{f \in \mathcal{F}} \left| \sum_{ i = 1}^n \epsilon_i \{ f(z_i) - f^*(z_i) \} \right| \right).
	\]
\end{theorem}

The last theorem we need is a concentration inequality due to \cite{Bousquet:2002}. We give a version as presented in \cite{vandegeer2008}.

\begin{theorem}[Bousequet's concentration theorem]\label{Thm: concentration theorem Bousquet}
	Suppose $Z_1, \dots, Z_n$ and all $g \in \mathcal{G}$ satisfy the following conditions for some real valued constants $\eta_n$ and $\tau_n$
	\[
	\Vert g \Vert_{\infty} \le \eta_n, \; \forall g \in \mathcal{G}
	\]
	and
	\[
	\frac{1}{n} \sum_{ i = 1}^n \textnormal{Var}(g(Z_i)) \le \tau_n^2, \; \forall g \in \mathcal{G}.
	\]
	Define
	\[
	\textbf{Z} \coloneqq \sup_{g \in \mathcal{G}} \left| \frac{1}{n}\sum_{ i = 1}^n g(Z_i) - \mathbb{E}[g(Z_i)]  \right|.
	\]
	Then for any $z>0$
	\[
	P\left(\textbf{Z} \ge \mathbb{E}[\textbf{Z}] + z \sqrt{2(\tau_n^2 + 2 \eta_n\mathbb{E}[\textbf{Z}])} + \frac{2z^2\eta_n}{3} \right) \le \exp(-nz^2).
	\]
\end{theorem}

\begin{remark}
	Looking at the original paper of \cite{Bousquet:2002}, their result looks quite different at first. To see that the above falls into their framework, set the variables in \cite{Bousquet:2002} as follows
	\begin{align*}
		f(Z_i) &= (g(Z_i) - \mathbb{E}[g(Z_i)])/(2\eta_n), &\tilde{Z}_k = \sup_f \vert \sum_{ i \neq k} f(Z_i) \vert, \\
		f_k &= \arg\sup_f \vert \sum_{ i \neq k} f(Z_i) \vert, &\tilde{Z}_k' = \vert \sum_{ i = 1}^n f_k(Z_i) \vert - \tilde{Z}_k \\
		\tilde{Z} &= \frac{2\eta_n}{n}\textbf{Z}.
	\end{align*}
	Now apply Theorem 2.1 in \cite{Bousquet:2002}, choosing for their $(Z, Z_1, \dots, Z_n)$ the above defined $(\tilde{Z}, \tilde{Z}_1, \dots, \tilde{Z}_n)$, for their $(Z_1', \dots, Z_n')$ the above defined $(\tilde{Z}_1', \dots, \tilde{Z}_n')$ and setting $u = 1$ and $\sigma^2 = \frac{\tau_n^2}{4\eta_n^2}$ in their theorem: The result is exactly Theorem \ref{Thm: concentration theorem Bousquet} above.
\end{remark}
Finally we have a Lemma derived from Hoeffding's inequality. The proof can be found in \cite{vandegeer2011}, Lemma 14.14 (here we use the special case of their Lemma for $m=1$).
\begin{lemma}\label{Lem: 14.14}
	Let $\mathcal{G} = \{g_1, \dots, g_p \}$ be a set of real valued functions on $\mathcal{Z}$ satisfying for all $i= 1, \dots, n$ and all $j=1, \dots, p$
	\[
	\mathbb{E}[g_j(Z_i)] = 0, \; \vert g_j(Z_i) \vert \le c_{ij}
	\]
	for some positive constants $c_{ij}$.
	Then
	\[
	\mathbb{E}\left[ \max_{1 \le j \le p} \left| \sum_{ i = 1}^n g_j(Z_i) \right| \right] \le \left[ 2\log(2p)\right]^{1/2} \max_{1 \le j \le p} \left[ \sum_{ i = 1}^n c_{ij}^2   \right]^{1/2}.
	\]
\end{lemma}

\subsection*{The expectation of \texorpdfstring{$Z_M$}{Z\_M}}

Recall the definition of $Z_M$ 
\[
Z_M \coloneqq \sup_{\substack{\bar{\theta} \in \bar{\Theta}_{\text{loc}}, \\  \Vert \bar{\theta} - \bar{\theta}_0 \Vert_1 \le M}} \vert \bar{v}_n(\bar{\theta}) - \bar{v}_n(\bar{\theta}_0) \vert,
\]
where $\bar{v}_n$ denotes the re-parametrized empirical process. 
Recall, that there is a constant $c \in \R$ such that uniformly $\vert Z_{ij,k} \vert \le c, 1 \le i \neq j \le n, k = 1, \dots, p$.

\begin{lemma}\label{Lem: 14.20}\label{Lem: 14.20 covariates} For any $M > 0$ we have
	\[
	\mathbb{E}[Z_M] \le 8 M(1 \vee c) \sqrt{\frac{2\log(2(2n+p+1))}{N}}.
	\]
\end{lemma}

\begin{proof}
	Let $\epsilon_{ij}, i \neq j,$ be a Rademacher sequence independent of $A_{ij}, Z_{ij}, i \neq j$. We first want to use the symmetrization theorem \ref{Thm: symmetrization theorem}: For the random variables $Z_1, \dots, Z_n$ we choose $T_{ij} = (A_{ij}, \bar{X}_{ij}^T,1, Z_{ij}^T)^T \in \{0,1\} \times \R^{2n+1+p}$. 
	For any $\bar{\theta} \in \bar{\Theta}_{\text{loc}}$ we consider the functions
	\[
	g_{\bar{\theta}}(T_{ij}) = \frac{1}{N} \left\{-A_{ij} \bar D_{ij}^T(\bar \theta - \bar \theta_0) + \log(1 + \exp(\bar D_{ij}^T\bar \theta)) - \log(1 + \exp(\bar D_{ij}^T \bar \theta_0)) \right\}
	\]
	and the function set $\mathcal{G} = \mathcal{G}(M) \coloneqq \{ g_{\bar{\theta}} : \bar{\theta} \in \bar{ \Theta}_{\text{loc}}, \Vert \bar{\theta} - \bar{\theta}_0 \Vert_1 \le M \}$.
	Note, that
	\[
	\bar{v}_n(\bar{\theta}) - \bar{v}_n(\bar{\theta}_0) = \sum_{ i \neq j} \{g_{\bar{\theta}}(T_{ij}) - \mathbb{E}[g_{\bar{\theta}}(T_{ij}) ]\}.
	\]	
	Then, the symmetrization theorem gives us
	\begin{align*}
		\E[Z_M ] &= \mathbb{E}\left[ \sup_{g_{\bar{\theta}} \in \mathcal{G}} \left\vert \sum_{ i \neq j} g_{\bar{\theta}}(T_{ij}) - \mathbb{E}[g_{\bar{\theta}}(T_{ij}) ] \right\vert  \right] \\
		&\le 2 \mathbb{E}\left[ \sup_{g_{\bar{\theta}} \in \mathcal{G}} \left\vert \sum_{ i \neq j} \epsilon_{ij} g_{\bar{\theta}}(T_{ij}) \right\vert   \right].
	\end{align*}
	Next, we want to apply the contraction Theorem \ref{Thm: contraction theorem Ledoux Talagrand}. Denote $T = (T_{ij})_{i\neq j}$ and let $\mathbb{E}_{T}$ be the conditional expectation given $T$. We need the conditional expectation at this point, because Theorem \ref{Thm: contraction theorem Ledoux Talagrand} requires non-random arguments in the functions. This does not hinder us, as later we will simply take iterated expectations, cancelling out the conditional expectation, see below. For the functions $g_i$ in Theorem \ref{Thm: contraction theorem Ledoux Talagrand} we choose
	\[
	g_{ij}(x) = \frac{1}{2} \{-A_{ij}x + \log(1 + \exp(x))\}
	\]
	Note, that $ \log(1 + \exp(x))$ has derivative bounded by one and thus is Lipschitz continuous with constant one by the Mean Value Theorem. Thus, all $g_{ij}$ are also Lipschitz continuous with constant $1$:
	\[
	\vert g_{ij}(x) - g_{ij}(x') \vert \le \frac{1}{2} \{ \vert A_{ij} (x - x') \vert  + \vert  \log(1 + \exp(x)) -  \log(1 + \exp(x')) \vert \} \le \vert x - x' \vert.
	\]
	For the function class $\mathcal{F}$ in Theorem $\ref{Thm: contraction theorem Ledoux Talagrand}$ we choose $\mathcal{F} = \mathcal{F}_M \coloneqq \{ f_{\bar{\theta}} : \bar{\theta} \in \bar{\Theta}_{\text{loc}}, \Vert \bar{\theta} - \bar{\theta}_0 \Vert_1 \le M \}$ and pick $f^* = f_{\bar{\theta}_0}$. Then, by Theorem \ref{Thm: contraction theorem Ledoux Talagrand}
	\begin{align*}
		\mathbb{E}_T &\left[ \sup_{\substack{\bar{\theta} \in \bar{\Theta}_{\text{loc}}, \\ \Vert \bar{\theta} -\bar{\theta}_0 \Vert_1 \le M}} \left\vert \frac{1}{N} \sum_{ i \neq j} \epsilon_{ij}(g_{ij}(f_{\bar{\theta}}((\bar{X}_{ij}^T, 1, Z_{ij}^T)^T) )- g_{ij}(f_{\bar{\theta}_0}((\bar{X}_{ij}^T, 1, Z_{ij}^T)^T))) \right\vert  \right] \\
		&\le 2 \mathbb{E}_T \left[ \sup_{\substack{\bar{\theta} \in \bar{\Theta}_{\text{loc}}, \\ \Vert \bar{\theta} -\bar{\theta}_0 \Vert_1 \le M}} \left\vert \frac{1}{N} \sum_{ i \neq j} \epsilon_{ij}(f_{\bar{\theta}}((\bar{X}_{ij}^T, 1, Z_{ij}^T)^T) - f_{\bar{\theta}_0}((\bar{X}_{ij}^T, 1, Z_{ij}^T)^T)) \right\vert  \right].
	\end{align*}
	Recall that we can express the functions $f_{\bar{\theta}} = f_{\bar \alpha, \bar{\beta}, \mu, \gamma}$ as
	\[
	f_{\bar \alpha, \bar{\beta}, \mu, \gamma} (\, . \,)= \sum_{i = 1}^n \bar{\alpha}_i e_i(\, . \,) + \sum_{i = n+1}^{2n} \bar{\beta}_{i-n} e_i(\, . \,) +  \mu e_{2n+1}(\, . \,) + \sum_{ i = 1}^p \gamma_ie_{2n+1+i}(\,.\,),
	\]
	where $e_i(\,.\,)$ is the projection on the $i$-th coordinate. Consider any $\bar{\theta} \in \bar{\Theta}_{\text{loc}}$ with $\Vert \bar{\theta} - \bar{\theta}_0 \Vert_1 \le M$. For the sake of a compact representation we use our shorthand notation $\bar{ \theta} = (\bar{\theta}_i)_{i=1}^{2n+1+p}$ where the components $\theta_i$ are defined in the canonical way and we also simply write $e_k(\bar X_{ij}, 1, Z_{ij})$ for the projection of the 
	the vector $(\bar X_{ij}^T, 1, Z_{ij}^T)^T \in \R^{2n+p+1}$ to its $k$-th component, i.e. instead of $e_k((\bar X_{ij}^T, 1, Z_{ij}^T)^T)$. Then,
	\begin{align*}
		&\left\vert \frac{1}{N} \sum_{ i \neq j} \epsilon_{ij}(f_{\bar{\theta}}((\bar{X}_{ij}^T, 1, Z_{ij}^T)^T) - f_{\bar{\theta}_0}((\bar{X}_{ij}^T, 1, Z_{ij}^T)^T)) \right\vert \\
		&= \left\vert \frac{1}{N} \sum_{ i \neq j} \epsilon_{ij} \left( \sum_{k = 1}^{2n+p+1} (\bar{\theta}_k - \bar{\theta}_{0,k}) e_k(\bar{X}_{ij}, 1, Z_{ij})\right)\right\vert \\
		&\le \frac{1}{N} \sum_{k = 1}^{2n+p+1} \left\lbrace \vert \bar{\theta}_k - \bar{\theta}_{0,k} \vert
		\max_{1 \le l \le 2n+p+1} \left| \sum_{ i \neq j} \epsilon_{ij} e_l(\bar{X}_{ij}, 1, Z_{ij}) \right| \right\rbrace \\
		&\le M \max_{1 \le l \le 2n+p+1} \left|\frac{1}{N} \sum_{ i \neq j} \epsilon_{ij} e_l(\bar{X}_{ij}, 1, Z_{ij}) \right|.
	\end{align*}
	Note, that the last expression no longer depends on $\bar{\theta}$. To bind the right hand side in the last expression we use Lemma \ref{Lem: 14.14}: In the language of the Lemma, choose $Z_1, \dots, Z_n$ as $T_{ij} = (\epsilon_{ij}, \bar{X}_{ij}^T, 1, Z_{ij}^T)^T$. We choose for the $p$ in the formulation of the Lemma $2n+p+1$ and pick for our functions
	\[
	g_{k}(T_{ij}) = \frac{1}{N}\epsilon_{ij}e_k(\bar{X}_{ij}, 1, Z_{ij}), k=1,\dots, 2n+p+1.
	\]
	Note, that then $\mathbb{E}[g_{k}(T_{ij})] = 0$. We want to employ Lemma \ref{Lem: 14.14} which requires us to bound $\vert g_k(T_{ij}) \vert \le c_{ij,k}$ for all $i \neq j$ and $k = 1, \dots, n+1+p$.
	
	For any fixed $1 \le k \le n$ we have
	\[
	\vert g_k (T_{ij}) \vert \le 
	\begin{cases}
		\frac{\sqrt{n}}{N} = \frac{1}{(n-1)\sqrt{n}}, & i \text{ or } j = k \\
		0, & \text{otherwise}.
	\end{cases}
	\]
	Note that the first case occurs exactly $(n-1)$ times for each $k$. Thus, for any $k \le 2n$,
	\[
	\sum_{ i \neq j} c_{ij,k}^2 = \left( \frac{1}{(n-1)\sqrt{n}}  \right)^2 (n-1) = \frac{1}{N}.
	\]
	If $k = 2n+1$, $\vert g_k (T_{ij}) \vert = 1/N$ and hence
	\[
	\sum_{ i \neq j} c_{ij,2n+1}^2 = \frac{1}{N}.
	\]
	Finally, if $k > 2n+1$, $\vert g_k(T_{ij}) \vert \le c /N$ and therefore,
	\[
	\sum_{ i \neq j} c_{ij, k}^2 \le \frac{c^2}{N}.
	\]
	In total, this means
	\[
	\max_{1 \le k \le 2n+1+p} \sum_{ i \neq j} c_{ij,k}^2 \le \frac{1 \vee c^2}{N}.
	\]
	Therefore, an application of Lemma \ref{Lem: 14.14} results in
	\begin{align*}
		\mathbb{E}\left[ \max_{1 \le l \le 2n+p+1} \left|\frac{1}{N} \sum_{ i \neq j} \epsilon_{ij} e_l(\bar{X}_{ij}, Z_{ij}) \right| \right] &\le \sqrt{2\log(2(2n+1+p))} \max_{1 \le k \le 2n+1+p} \left[\sum_{ i \neq j} c_{ij,k}^2\right]^{1/2} \\
		&\le \sqrt{2\log(2(2n+1+p))} \sqrt{\frac{1 \vee c^2}{N}} \\
		&= \sqrt{\frac{2\log(2(2n+1+p))}{N}} (1 \vee c).
	\end{align*}
	Putting everything together, we obtain
	\begin{align*}
		\mathbb{E}[Z_M] &\le 2 \mathbb{E}\left[ \sup_{\substack{\bar{\theta} \in \bar{\Theta}_{\text{loc}}, \\ \Vert \bar{\theta} -\bar{\theta}_0 \Vert_1 \le M}}  \left\vert \frac{1}{N} \sum_{ i \neq j} \epsilon_{ij}(-A_{ij}(f_{\bar{\theta}}(\bar{X}_{ij}, 1,Z_{ij}) - f_{\bar{\theta}_0}  (\bar{X}_{ij},1, Z_{ij}))) \right\vert  \right] \\
		&= 2 \mathbb{E}\left[ \mathbb{E}_T \left[ \sup_{\substack{\bar{\theta} \in \bar{\Theta}_{\text{loc}}, \\ \Vert \bar{\theta} -\bar{\theta}_0 \Vert_1 \le M}}  \left\vert \frac{1}{N} \sum_{ i \neq j} \epsilon_{ij}(-A_{ij}(f_{\bar{\theta}}(\bar{X}_{ij},1, Z_{ij}) - f_{\bar{\theta}_0}  (\bar{X}_{ij},1, Z_{ij}))) \right\vert  \right]\right] \\
		&\le 8 \mathbb{E}\left[ \mathbb{E}_T \left[ \sup_{\substack{\bar{\theta} \in \bar{\Theta}_{\text{loc}}, \\ \Vert \bar{\theta} -\bar{\theta}_0 \Vert_1 \le M}}  \left\vert \frac{1}{N} \sum_{ i \neq j} \epsilon_{ij}(f_{\bar{\theta}}(\bar{X}_{ij},1, Z_{ij}) - f_{\bar{\theta}_0}  (\bar{X}_{ij},1, Z_{ij})) \right\vert  \right]\right] \\
		&\le 8M \mathbb{E} \left[ \mathbb{E}_T \left[ \max_{1 \le l \le 2n+p+1} \left|\frac{1}{N} \sum_{ i \neq j} \epsilon_{ij} e_l(\bar{X}_{ij}, 1,Z_{ij}) \right| \right]\right] \\
		&\le 8M\sqrt{\frac{2\log(2(2n+1+p))}{N}} (1 \vee c).
	\end{align*}
	This concludes the proof.
\end{proof}

We now want to show that $Z_M$ does not deviate too far from its expectation. The proof relies on the concentration theorem due to Bousquet, Theorem \ref{Thm: concentration theorem Bousquet}.

\begin{corollary}\label{Kor: Probability bound Z_M}\label{Kor: Probability bound Z_M covariates}
	Pick any confidence level $t > 0$. 
	Let
	\[
	a_n \coloneqq \sqrt{\frac{2\log(2(2n+p+1))}{N}} (1 \vee c)
	\]
	and choose $\lambda_0 = \lambda_0(t,n)$ as
	\[
	\lambda_0 = 8a_n + 2 \sqrt{\frac{t}{N}( 11 (1 \vee (c^2p) ) + 16(1 \vee c) \sqrt{n} a_n  )} + \frac{4t(1 \vee c) \sqrt{n}}{3N}
	\]
	Then, we have the inequality
	\[
	P\left( Z_M \ge M\lambda_0  \right) \le \exp(-t).
	\]
\end{corollary}

\begin{proof}
	We want to apply Bousquet's concentration Theorem \ref{Thm: concentration theorem Bousquet}. For the random variables $Z_i$ in the formulation of the theorem we choose once more $T_{ij} = (A_{ij}, \bar{X}_{ij},1, Z_{ij}), i\neq j,$ and as functions we consider
	\begin{align*}
		g_{\bar{\theta}}(T_{ij}) &= -A_{ij} \bar D_{ij}^T(\bar \theta - \bar \theta_0) + \log(1 + \exp(\bar D_{ij}^T\bar \theta)) - \log(1 + \exp(\bar D_{ij}^T \bar \theta_0)), \\
		\mathcal{G} &= \mathcal{G}_M \coloneqq \{ g_{\bar{\theta}} : \bar{\theta} \in \bar{\Theta}_{\text{loc}}, \Vert \bar{\theta} -\bar{\theta}_0 \Vert_1 \le M \}.
	\end{align*}
	Then, we have
	\[
	Z_M = \sup_{g_{\bar{\theta}} \in \mathcal{G}} \frac{1}{N} \left\vert \sum_{ i \neq j} \{g_{\bar{\theta}}(T_{ij}) - \mathbb{E}[g_{\bar{\theta}}(T_{ij}) ] \} \right\vert.
	\]
	To apply Theorem \ref{Thm: concentration theorem Bousquet}, we need to bound the infinity norm of $g_{\bar{\theta}}$. Recall that we denote the distribution of $[\bar{X} | 1| Z]$ by $\bar{Q}$ and the infinity norm is defined as the $\bar{Q}$-almost sure smallest upper bound on the value of $g_{\bar{\theta}}$. We have for any $g_{\bar{\theta}} \in \mathcal{G}$, using the Lipschitz continuity of $\log(1 + \exp(x))$:
	\begin{align*}
		\vert g_{\bar{\theta}}(T_{ij}) \vert &\le \vert \bar D_{ij}^T(\bar \theta - \bar \theta_0) \vert + \vert  \log(1 + \exp(\bar D_{ij}^T\bar \theta)) - \log(1 + \exp(\bar D_{ij}^T \bar \theta_0)) \vert \\
		&\le 2 \vert \bar D_{ij}^T(\bar \theta - \bar \theta_0) \vert \\
		&\le 2 \Vert \vartheta - \vartheta_0 \Vert_1 + \vert \mu - \mu_0 \vert + c \Vert \gamma - \gamma_0 \Vert_1.
		\shortintertext{Thus,}
		\Vert g_{\bar{\theta}} \Vert_{\infty} &\le 2 \Vert \vartheta - \vartheta_0 \Vert_1 + \vert \mu - \mu_0 \vert + c \Vert \gamma - \gamma_0 \Vert_1 \\
		&\le 2(1 \vee c) \Vert \theta - \theta_0 \Vert_1 \\
		&\le 2(1 \vee c) \sqrt{n} M \eqqcolon \eta_n.
	\end{align*}
	For the last inequality we used that for any $\theta$ with $\Vert \bar{\theta} - \bar{\theta}_0 \Vert_1 \le M$ it follows that $\Vert \theta - \theta_0 \Vert_1 \le \sqrt{n}M$, which is possibly a very generous upper bound. This does not matter, however, as the term associated with the above bound will be negligible, as we shall see.
	
	The second requirement of Theorem \ref{Thm: concentration theorem Bousquet} is that the average variance of $g_{\bar{\theta}}(T_{ij})$ has to be uniformly bounded. To that end we calculate
	\begin{align*}
		\frac{1}{N} \sum_{ i \neq j}\text{Var}(g_{\bar{\theta}}(T_{ij})) &= \frac{1}{N} \sum_{ i \neq j} \text{Var}(-A_{ij}D_{ij}^T(\theta - \theta_0)) \\
		& + \frac{1}{N} \sum_{ i \neq j}\text{Var}(\log(1 + \exp(\bar D_{ij}^T\bar \theta)) - \log(1 + \exp(\bar D_{ij}^T \bar \theta_0))) \\
		& +\frac{2}{N} \sum_{ i \neq j} \text{Cov}(-A_{ij}D_{ij}^T(\theta - \theta_0), \log(1 + \exp(\bar D_{ij}^T\bar \theta)) - \log(1 + \exp(\bar D_{ij}^T \bar \theta_0))).
	\end{align*}
	Let us look at these terms in term. For the first term, we obtain
	\begin{align*}
		\frac{1}{N} \sum_{ i \neq j} \text{Var}(-A_{ij}D_{ij}^T(\theta - \theta_0)) &\le \frac{1}{N} \sum_{ i \neq j} \E[(-A_{ij}D_{ij}^T(\theta - \theta_0))^2] \le \E \left[ \frac{1}{N} \sum_{ i \neq j} (D_{ij}^T(\theta - \theta_0))^2 \right].
	\end{align*}
	For the second term we get
	\begin{align*}
		\frac{1}{N} \sum_{ i \neq j}\text{Var}(&\log(1 + \exp(\bar D_{ij}^T\bar \theta)) - \log(1 + \exp(\bar D_{ij}^T \bar \theta_0))) \\
		&\le \frac{1}{N} \sum_{ i \neq j} \E[(\log(1 + \exp(\bar D_{ij}^T\bar \theta)) - \log(1 + \exp(\bar D_{ij}^T \bar \theta_0)))^2]\\
		&\le  \E \left[ \frac{1}{N} \sum_{ i \neq j} (D_{ij}^T(\theta - \theta_0))^2 \right].
	\end{align*}
	The last term decomposes as
	\begin{align*}
		\frac{2}{N} \sum_{ i \neq j} & \text{Cov}(-A_{ij}D_{ij}^T(\theta - \theta_0), \log(1 + \exp(\bar D_{ij}^T\bar \theta)) - \log(1 + \exp(\bar D_{ij}^T \bar \theta_0))) \\
		&= \frac{2}{N} \sum_{ i \neq j}  \E[-A_{ij}D_{ij}^T(\theta - \theta_0) \cdot (\log(1 + \exp(\bar D_{ij}^T\bar \theta)) - \log(1 + \exp(\bar D_{ij}^T \bar \theta_0)))] \\
		&\quad - \frac{2}{N} \sum_{ i \neq j}  \E[-A_{ij}D_{ij}^T(\theta - \theta_0)] \cdot \E[ \log(1 + \exp(\bar D_{ij}^T\bar \theta)) - \log(1 + \exp(\bar D_{ij}^T \bar \theta_0))]
	\end{align*}
	For the first term in that decomposition we have
	\begin{align*}
		\frac{2}{N} &\sum_{ i \neq j}  \left\vert\E[-A_{ij}D_{ij}^T(\theta - \theta_0) \cdot (\log(1 + \exp(\bar D_{ij}^T\bar \theta)) - \log(1 + \exp(\bar D_{ij}^T \bar \theta_0)))] \right\vert \\
		&\le \frac{2}{N} \sum_{ i \neq j}  \E[\vert D_{ij}^T(\theta - \theta_0)\vert \cdot \vert\log(1 + \exp(\bar D_{ij}^T\bar \theta)) - \log(1 + \exp(\bar D_{ij}^T \bar \theta_0))\vert] \\
		&\le \frac{2}{N} \sum_{ i \neq j}  \E[\vert D_{ij}^T(\theta - \theta_0)\vert^2] 
	\end{align*}
	and for the second term using the same arguments, we get
	\[
	\frac{2}{N} \sum_{ i \neq j}  \E[-A_{ij}D_{ij}^T(\theta - \theta_0)] \cdot \E[ \log(1 + \exp(\bar D_{ij}^T\bar \theta)) - \log(1 + \exp(\bar D_{ij}^T \bar \theta_0))] \le \frac{2}{N} \sum_{ i \neq j}  \E[\vert D_{ij}^T(\theta - \theta_0)\vert]^2. 
	\]
	Meaning that in total
	\begin{align*}
		\frac{2}{N} \sum_{ i \neq j} &\left\vert \text{Cov}(-A_{ij}D_{ij}^T(\theta - \theta_0), \log(1 + \exp(\bar D_{ij}^T\bar \theta)) - \log(1 + \exp(\bar D_{ij}^T \bar \theta_0))) \right\vert \\
		&\le \frac{2}{N} \sum_{ i \neq j}  \E[\vert D_{ij}^T(\theta - \theta_0)\vert^2] + \frac{2}{N} \sum_{ i \neq j}  \E[\vert D_{ij}^T(\theta - \theta_0)\vert]^2.
	\end{align*}
	In total, we thus get
	\begin{align}\label{Eq: decomposition of variance directed}
		\frac{1}{N} \sum_{ i \neq j}\text{Var}(g_{\bar{\theta}}(T_{ij})) \le 4 \cdot \E \left[ \frac{1}{N} \sum_{ i \neq j} (D_{ij}^T(\theta - \theta_0))^2 \right] + \frac{2}{N} \sum_{ i \neq j}  \E[\vert D_{ij}^T(\theta - \theta_0)\vert]^2.
	\end{align}
	Furthermore,
	\begin{align*}
		\frac{1}{N} \sum_{ i \neq j} (D_{ij}^T(\theta - \theta_0))^2 &= \frac{1}{N} \sum_{ i \neq j} (\alpha_i + \beta_j + \mu - \alpha_{0,i} - \beta_{0,j} - \mu_0 + (\gamma - \gamma_0)^TZ_{ij})^2 \\
		&\overset{\text{Cauchy-Schwarz}}{\le} \frac{4}{N} \sum_{ i \neq j} \left\{ (\alpha_i - \alpha_{0,i})^2 + (\beta_j - \beta_{0,j})^2 + (\mu - \mu_0)^2 + ((\gamma - \gamma_0)^TZ_{ij})^2\right\}.
	\end{align*}
	Recall that for any $x \in \R^p, \Vert x \Vert_2 \le \Vert x \Vert_1 \le \sqrt{p}\Vert x \Vert_2$ and note that
	\[
	\vert(\gamma - \gamma_0)^TZ_{ij}\vert \le c \Vert \gamma - \gamma_0 \Vert_1 \le c\sqrt{p}\Vert \gamma  - \gamma_0 \Vert_2.
	\]	
	Then, from the above
	\begin{align}\label{Eq: bound on sum D_ij theta -theta*}
		\begin{split}
			\frac{1}{N} \sum_{ i \neq j} (D_{ij}^T(\theta - \theta_0))^2 &\le \frac{4}{N} \sum_{ i \neq j} \left\{ (\alpha_i - \alpha_{0,i})^2 + (\beta_j - \beta_{0,j})^2 + (\mu - \mu_0)^2 + c^2 p \Vert \gamma - \gamma_0 \Vert_2^2 \right\} \\
			&= 4 \left( (\mu - \mu_0)^2 + c^2p \Vert \gamma - \gamma_0 \Vert_2^2 + \frac{1}{N}   \sum_{ i \neq j} \left\{ (\alpha_i - \alpha_{0,i})^2 + (\beta_j - \beta_{0,j})^2 \right\} \right) \\
			&= 4 \left( (\mu - \mu_0)^2 + c^2p \Vert \gamma - \gamma_0 \Vert_2^2 + \frac{1}{N}   (n-1) \Vert \vartheta - \vartheta_0 \Vert_2^2 \right) \\
			&= 4 \left( (\mu - \mu_0)^2 + c^2p \Vert \gamma - \gamma_0 \Vert_2^2 + \left\Vert \frac{1}{\sqrt{n}} (\vartheta - \vartheta_0) \right\Vert_2^2 \right) \\
			&= 4 \left( (\mu - \mu_0)^2 + c^2p \Vert \gamma - \gamma_0 \Vert_2^2 + \Vert  \bar{\vartheta} - \bar{\vartheta}_0 \Vert_2^2 \right) \\
			&\le 4(1 \vee (c^2p)) \Vert \bar{\theta} - \bar{\theta}_0 \Vert_2^2 \\
			&\le 4(1 \vee (c^2p)) \Vert \bar{\theta} - \bar{\theta}_0 \Vert_1^2 \\
			&\le 4(1 \vee (c^2p)) M^2.
		\end{split}
	\end{align}
	Notice that for the second summand on the right-hand side in \eqref{Eq: decomposition of variance directed}, we have
	\begin{align*}
		\frac{2}{N} \sum_{ i \neq j}  \E[\vert D_{ij}^T(\theta - \theta_0)\vert]^2 &= \frac{2}{N} \sum_{ i \neq j} (\alpha_i + \beta_j + \mu - \alpha_{0,i} - \beta_{0,j} - \mu_0 + (\gamma - \gamma_0)^T\E[Z_{ij}])^2 \\
		&=\frac{2}{N} \sum_{ i \neq j} (\alpha_i + \beta_j + \mu - \alpha_{0,i} - \beta_{0,j} - \mu_0)^2. 
	\end{align*}
	So that we may use the same steps as in \eqref{Eq: bound on sum D_ij theta -theta*} to conclude that
	\[
	\frac{2}{N} \sum_{ i \neq j}  \E[\vert D_{ij}^T(\theta - \theta_0)\vert]^2 \le 6(1 \vee (c^2p)) M^2.
	\]
	Such that in total,
	\begin{align*}
		\frac{1}{N} \sum_{ i \neq j}\text{Var}(g_{\bar{\theta}}(T_{ij})) \le 22(1 \vee (c^2p)) M^2 \coloneqq \tau_n^2.
	\end{align*}
	Applying Bousquet's concentration Theorem \ref{Thm: concentration theorem Bousquet} with $\eta_n, \tau_n$ defined above, we obtain for all $z > 0$
	\begin{align}\label{Eq: prob inequality I covariates}
		\begin{split}
			\exp&\left(-Nz^2\right) \ge P\left( Z_M \ge \mathbb{E}[Z_M] + z \sqrt{2(\tau_n^2 + 2\eta_n \mathbb{E}[Z_M])} + \frac{2z^2\eta_n}{3} \right) \\
			&=  P\left( Z_M \ge \mathbb{E}[Z_M] + z \sqrt{2(22(1 \vee (c^2p)) M^2 + 4(1 \vee c) \sqrt{n} M  \mathbb{E}[Z_M])} + \frac{4z^2(1 \vee c) \sqrt{n} M }{3} \right).
		\end{split}
	\end{align}
	From Lemma \ref{Lem: 14.20 covariates}, we know
	\[
	\mathbb{E}[Z_M] \le 8 M \sqrt{\frac{2\log(2(2n+p+1))}{N}}(1 \vee c) = 8Ma_n.
	\]
	Using this, we obtain from (\ref{Eq: prob inequality I covariates})
	\begin{align*}
		\exp\left(-Nz^2\right)& \ge P\left( Z_M \ge 8Ma_n + z \sqrt{2(22(1 \vee (c^2p)) M^2 + 32(1 \vee c) \sqrt{n} M^2 a_n)} + \frac{4z^2(1 \vee c) \sqrt{n} M }{3} \right)  \\
		&= P\left( Z_M \ge M \left( 8a_n + 2z \sqrt{11(1 \vee (c^2p) ) + 16(1 \vee c) \sqrt{n}a_n } + \frac{4z^2(1 \vee c) \sqrt{n}}{3} \right)\right).
	\end{align*}
	Now, pick $z = \sqrt{t/N}$ to get
	\[
	P\left( Z_M \ge M \left( 8a_n + 2 \sqrt{\frac{t}{N}( 11 (1 \vee (c^2p) ) + 16(1 \vee c) \sqrt{n} a_n  )} + \frac{4t(1 \vee c) \sqrt{n}}{3N} \right) \right) \le \exp(-t),
	\]
	which is the claim.
\end{proof}

\subsection{Putting it all together}

\begin{proof}[Proof of Theorem \ref{Thm: consistency directed}]
	Theorem \ref{Thm: consistency directed} now follows from Theorem \ref{Thm: Consistency covariates} and corollary \ref{Kor: Probability bound Z_M covariates}. Recall the definition of $K_n$ in \eqref{Eq: Def K_n}, which simplifies to
	\begin{align*}
		K_n &= 2 \frac{(1 + \exp(r_{n,0}))^2}{\exp(r_{n,0})} = 2 \frac{\left( 1 + \exp\left( - \text{logit}(\rho_{n,0})  \right) \right)^2}{\exp\left( - \text{logit}(\rho_{n,0}) \right)} = \frac{2}{\rho_{n,0}}.
	\end{align*}
	Thus, under the conditions of Theorem \ref{Thm: consistency directed}, we have with high probability by Theorem \ref{Thm: Consistency covariates} and Corollary \ref{Kor: Probability bound Z_M covariates},
	\begin{equation*}
		\mathcal{E}(\hat{\theta}) + \bar{ \lambda} \left( \frac{1}{\sqrt{n}} \Vert \hat{\vartheta} - \vartheta_0  \Vert_1 + \vert \hat{\mu} - \mu_0 \vert + \Vert \hat{\gamma} - \gamma_0 \Vert_1  \right) \le C  \frac{s_{0,+}\bar{ \lambda}^2}{\rho_{n,0}}.
	\end{equation*}
	with constant $C = 64/c_{\min}$.
\end{proof}

\section{Proof of Theorem \ref{Thm: inference directed}}

\subsection{Inverting population and sample Gram matrices}\label{Sec: inverting Gram matrices}

Note that the function $f(x) = x(1-x)$ is monotonically increasing in $x$ for $x \le 1/2$ and monotonically decreasing in $x$ for $x \ge 1/2$. Thus, by considering the cases $p_{ij} \le 1/2$ and $p_{ij} \ge 1/2$  separately and using that $\rho_n \le 1/2$, we may employ the following lower bound for all $i \neq j$: $p_{ij}(\theta_0)(1-p_{ij}(\theta_0)) \ge 1/2 \rho_n$. Also, recall that by assumption \ref{Assum: minimum EW directed}, the minimum eigenvalue  $\lambda_{\text{min}}$ of $\mathbb{E}[Z^TZ/N]$ stays uniformly bounded away from zero for all $n$. Then, for any $n$ and $v \in \R^{p+1} \backslash \{0\}$ with components $v = (v_1, v_R^T)^T, v_R \in \R^p$, we have 
\begin{align*}
	v^T \Sigma_\xi v &\ge \frac{1}{2} \rho_n v^T \frac{1}{N} \mathbb{E}[ D_\xi^T D_\xi] v = \frac{1}{2} \rho_n v^T \begin{pmatrix}
		1 & \textbf{0} \\
		\textbf{0} & \frac{1}{N} \E[ Z^TZ]
	\end{pmatrix} v \\
	& = \frac{1}{2} \rho_n \left(v_1^2 + v_R^T \frac{1}{N} \E[Z^TZ] v_R  \right) \\
	&\ge  \frac{1}{2} \rho_n (v_1^2 + \lambda_{\text{min}}\Vert v_R \Vert_2^2) \ge \frac{1}{2} \rho_n (1 \wedge c_{\min}) \Vert v \Vert_2^2 > 0.
\end{align*}
Hence, for finite $n$ all eigenvalues of $\Sigma_\xi$ are strictly positive and consequently this matrix is invertible. We now want to show that the same hold with high probability for the sample matrix $\hat{\Sigma}_\xi$.
Using the tools deployed in the proofs of Lemma \ref{Lem: Lemma 6 in KockTang} and \ref{Lem: Lemma 7 in KockTang} we can now show that with high probability the minimum eigenvalue of $D_\xi^TD_\xi / N$ is also strictly larger than zero, which means that $D_\xi^TD_\xi/N$ is invertible with high probability, from which the desired properties of $\hat{\Sigma}_\xi$ follow.
More precisely, recall the definition of $\kappa(A, m)$ for square matrices $A$ and dimensions $m$. We want to consider the expression $\kappa^2\left(\frac{1}{N} \mathbb{E}[D_\xi^TD_\xi], p+1\right)$ which simplifies to
\[
\kappa^2\left(\frac{1}{N} \mathbb{E}[D_\xi^TD_\xi], p+1 \right) \coloneqq \min_{v \in \R^{p+1} \backslash\{0\}} \frac{v^T\frac{1}{N} \mathbb{E}[D_\xi^TD_\xi] v}{\frac{1}{p+1} \Vert v \Vert_1^2}
\]
and compare it to $\kappa^2\left(\frac{1}{N}D_\xi^TD_\xi,p+1\right)$. By assumption \ref{Assum: minimum EW directed} and the argument above, we have
\[
\kappa^2\left(\frac{1}{N} \mathbb{E}[D_\xi^TD_\xi], p+1\right) \ge C > 0
\]
for a universal constant $C$ independent of $n$.
With $\delta = \max_{kl} \left\vert  \left(\frac{1}{N}D_\xi^TD_\xi\right)_{kl} - \left( \frac{1}{N}\E[D_\xi^TD_\xi]\right)_{kl} \right\vert$, by Lemma \ref{Lem: Lemma 6 in KockTang}, we have
\[
\kappa^2\left(\frac{1}{N}D_\xi^TD_\xi, p+1\right) \ge \kappa^2\left(\frac{1}{N}\E[D_\xi^TD_\xi], p+1\right) - 16\delta (p+1).
\]
By looking at the proof of Lemma \ref{Lem: Lemma 6 in KockTang}, we see that in this particular case we do not even need the factor $16(p+1)$ on the right hand side above, but this does not matter anyways, so we keep it.
By the exact same arguments we have used in the proof of Lemma \ref{Lem: Lemma 7 in KockTang} for the blocks \textcircled{5}, \textcircled{6}, \textcircled{8} and \textcircled{9}, we now get 
\[
\delta = O_P\left( N^{-1/2}  \right).
\]
Thus, for $n$ large enough, we have with high probability $\delta \le \frac{\lambda_{\min}}{32}$. Then, by Lemma \ref{Lem: Lemma 6 in KockTang}, with high probability and uniformly in $n$,
\[
\kappa^2\left(\frac{1}{N}D_\xi^TD_\xi, p+1\right) \ge \kappa^2\left(\frac{1}{N}\E[D_\xi^TD_\xi], p+1\right) - 16\delta (p+1) \ge \frac{\lambda_{\min} (p+1)}{2} \ge C > 0.
\]
Yet, if $\kappa^2\left(\frac{1}{N} D_\xi^TD_\xi, p+1\right) \ge C > 0$ uniformly in $n$, then for any $v \neq 0, v^T\frac{1}{N} D_\xi^TD_\xi v \ge C \Vert v \Vert_2^2$. But we also know that the minimum eigenvalue of $\frac{1}{N} D_\xi^TD_\xi$ is the largest possible $C$ such that this bound holds (it is actually tight with equality for the eigenvectors corresponding to the minimum eigenvalue). Therefore, with high probability, the minimum eigenvalue of $\frac{1}{N}D_\xi^TD_\xi$ stays uniformly bounded away from zero.
Thus, for any $v \in \R^{p+1}\backslash\{0\}$ and any finite $n$:
\[
\frac{1}{N}v^T D_\xi^T \hat{W}^2 D_\xi v \ge \min_{i \neq j} \{ p_{ij}(\hat{\theta}) (1 - p_{ij}(\hat{\theta}))\} \left(	v^T \frac{1}{N} D_\xi^T D_\xi v\right) \ge C \rho_n \Vert v \Vert_2^2 >0.
\]
Thus, $\text{mineval}\left(\frac{1}{N}D_\xi^T \hat{W}^2 D_\xi \right) \ge C \rho_n\text{mineval}
\left(\frac{1}{N}D_\xi^T D_\xi  \right) > 0$. That is, for every finite $n$, $\frac{1}{N} D_\xi^T \hat{W}^2 D_\xi$ is invertible with high probability.

\subsection{Goal and approach}\label{Sec:goal and approach}

Our strategy will be inverting the KKT conditions, similar to \cite{vandegeer2014}. 
Recall our discussion of the KKT conditions in Section \ref{Sec: model selection}.
By the same arguments, we find that $0$ has to be contained in the subdifferential of $\frac{1}{N}\mathcal{L}(\theta) + \lambda \Vert \beta \Vert_1$ at $\hat{\theta}$, where this time we consider the KKT conditions with respect to the original parameters $\theta$. That is, there exists a $\hat{z} \in \R^{2n+1+p}$ such that
\begin{equation*}
	0 = \frac{1}{N}\nabla \left.\mathcal{L}(\theta)\right|_{\theta = \hat{\theta}} + \lambda \hat{z},
\end{equation*}
where $\nabla\left.\mathcal{L}(\theta)\right|_{\theta = \hat{\theta}}$ is the gradient of $\mathcal{L}(\theta)$ evaluated at $\hat{\theta}$ and for $i=1, \dots, 2n, \hat{z}_i = 1$ if $\hat{\vartheta}_i > 0$ and $\hat{z}_i \in [-1,1]$ if $\hat{\vartheta}_i = 0$, and for $i = 2n+1, \dots, 2n+ 1+p, \hat{z}_i = 0$.

Denoting $\nabla_\xi \left.\mathcal{L}(\theta)\right|_{\theta = \hat{\theta}} \in \R^{p+1}$ the gradient of $\mathcal{L}$ with respect to the unpenalized parameters $\xi = (\mu, \gamma^T)^T$ only, evaluated at $\hat{\theta}$, we have
\begin{equation}\label{Eq: KKT gamma directed}
	0 = \nabla_\xi \left.\mathcal{L}(\theta)\right|_{\theta = \hat{\theta}}.
\end{equation}
\textbf{Goal:} We want to show that for $k=1, \dots, p+1$,
\[
\sqrt{N}\frac{\hat \xi_k - \xi_{0,k}}{\sqrt{\hat \Theta_{\xi,k,k}}} \rightarrow \mathcal{N}(0,1).
\]

\noindent\textbf{Approach:} Recall the definition of the "one-sample-version" of $\mathcal{L}$, i.e. $l_{\theta}: \{0,1\} \times \R^{2n+1 + p} \rightarrow \R$,
for $\theta = (\alpha^T, \beta^T, \mu, \gamma^T)^T \in \Theta$,
\begin{equation*}
	l_{\theta}(y,x) \coloneqq -y \theta^Tx + \log(1 + \exp(\theta^Tx).
\end{equation*}
Then, the negative log-likelihood is given by
\[
\mathcal{L}(\theta) = \sum_{ i \neq j} l_\theta(A_{ij}, D_{ij}^T)
\]
and
\begin{align*}
	\nabla\mathcal{L}(\theta) = \sum_{ i \neq j} \nabla l_\theta(A_{ij},  D_{ij}^T),  \quad H\mathcal{L}(\theta)= \sum_{ i \neq j} Hl_\theta(A_{ij},  D_{ij}^T),
\end{align*}
where $H$ denotes the Hessian with respect to $\theta$. 
Consider $l_\theta$ as a function in $\theta^T x$ and introduce:
\begin{equation}\label{Eq: Def l (a)}
	l(y,a) \coloneqq -ya + \log(1+\exp(a)),
\end{equation}
with second derivative: $\ddot l (y,a) = \partial_{a^2}l(y,a) = \frac{\exp(a)}{(1+\exp(a))^2}$. Note, that $\partial_{a^2}l(y,a)$ is Lipschitz continuous (it has bounded derivative $\vert \partial_{a^3}l(y,a) \vert \le 1/(6\sqrt{3})$; Lipschitz continuity then follows by the Mean Value Theorem). Doing a first order Taylor expansion in $a$ of $\dot l (y,a) = \partial_al(y,a)$ in the point $(A_{ij}, D_{ij}^T\theta_0)$ evaluated at $(A_{ij}, D_{ij}^T\hat \theta)$, we get
\begin{equation}\label{Eq: Taylor of l}
	\partial_a l(A_{ij}, D_{ij}\hat \theta) = \partial_a l(A_{ij}, D_{ij}^T\theta_0) + \partial_{a^2} l(A_{ij}, \alpha) D_{ij}^T(\hat \theta - \theta_0),
\end{equation}
for an $\alpha$ between $D_{ij}^T\hat \theta$ and $D_{ij}^T \theta_0$. By Lipschitz continuity of $\partial_{a^2}l$, we also find
\begin{align}\label{Eq: Lipschitz l a^2}
	\begin{split}
		\vert \partial_{a^2}l(A_{ij}, \alpha) D_{ij}^T(\hat \theta - \theta_0) - \partial_{a^2}l(A_{ij}, D_{ij}^T\hat \theta) D_{ij}^T(\hat \theta - \theta_0) \vert &\le \vert \alpha - D_{ij}^T\hat \theta \vert \vert D_{ij}^T(\hat \theta - \theta_0) \vert \\
		&\le \vert D_{ij}^T(\hat \theta - \theta_0) \vert^2,
	\end{split}
\end{align}
where the last inequality follows, because $\alpha$ is between $D_{ij}^T\hat \theta$ and $D_{ij}^T \theta_0$.

Consider the vector $P_n\nabla l_{\hat \theta}$: By equation (\ref{Eq: Taylor of l}), with $\alpha_{ij}$ between $D_{ij}^T\hat \theta$ and $D_{ij}^T \theta_0$,
\begin{align*}
	P_n \nabla l_{\hat \theta} &= \frac{1}{N} \sum_{ i \neq j} \left( \partial_{\theta_k} l (A_{ij}, D_{ij}^T \hat \theta) \right)_{k=1, \dots, 2n+1+p}, \quad \text{ as a }(2n+1+p)\times 1\text{-vector} \\
	&= \frac{1}{N} \sum_{ i \neq j} \dot l (A_{ij}, D_{ij}^T \hat \theta) D_{ij} \\
	&= \frac{1}{N} \sum_{ i \neq j} (\dot l (A_{ij}, D_{ij}^T \theta_0) + \ddot l(A_{ij}, \alpha_{ij}) D_{ij}^T (\hat \theta - \theta_0)) D_{ij} \\
	\shortintertext{which by (\ref{Eq: Lipschitz l a^2}) gives}
	& = P_n\nabla l_{\theta_0} + \frac{1}{N} \sum_{ i \neq j} D_{ij} \left\{\ddot l(A_{ij}, D_{ij}^T\hat \theta) D_{ij}^T(\hat \theta - \theta_0) + O(\vert D_{ij}^T(\hat \theta - \theta_0) \vert^2)\right\}. \\
	\shortintertext{Noticing that $\ddot l(A_{ij}, D_{ij}^T\hat \theta) = p_{ij}(\hat \theta) (1 - p_{ij}(\hat{\theta}))$ and thus $\sum_{ i \neq j}\ddot l(A_{ij}, D_{ij}^T\hat \theta)D_{ij} D_{ij}^T(\hat \theta - \theta_0)= D^T \hat W^2 D (\hat \theta - \theta_0)$:}
	&= P_n\nabla l_{\theta_0} + P_n Hl_{\hat{\theta}}(\hat \theta - \theta_0)  + O\left(\frac{1}{N}\sum_{i\neq j} D_{ij} \vert D_{ij}^T(\hat \theta - \theta_0) \vert^2 \right) \\
	&= P_n\nabla l_{\theta_0} + \frac{1}{N} D^T \hat W^2 D (\hat{ \theta} - \theta_0) + O\left(\frac{1}{N}\sum_{i \neq j} D_{ij} \vert D_{ij}^T(\hat \theta - \theta_0) \vert^2 \right),
\end{align*}
where the $O$ notation is to be understood componentwise.
Above, we have equality of two $((2n+1+p) \times 1)$-vectors. We are only interested in the portion relating to $\xi = (\mu, \gamma^T)^T$, that is, in the last $p+1$ entries. Introduce the $((2n+1+p) \times (2n+1+p))$-matrix
\begin{equation*}
	A = \begin{pmatrix}
		\textbf{0} & \textbf{0} \\
		\textbf{0} & \hat{ \Theta}_\xi
	\end{pmatrix},
\end{equation*}
where $\textbf{0}$ are zero-matrices of appropriate dimensions. Multiplying the above with $A$ on both sides gives:
\begin{equation}\label{Eq: Taylor P_n directed}
	AP_n \nabla l_{\hat \theta} = A P_n \nabla l_{\theta_0} + A \frac{1}{N} D^T \hat W^2 D (\hat{ \theta} - \theta_0) + A O\left(\frac{1}{N}\sum_{i \neq j} D_{ij} \vert D_{ij}^T(\hat \theta - \theta_0) \vert^2 \right).
\end{equation}
Let us consider these terms in turn: Multiplication by $A$ means that the first $n$ entries of any of the vectors above are zero. Hence we only need to consider the last $p+1$ entries.
The left-hand side of (\ref{Eq: Taylor P_n directed}) is equal to zero by (\ref{Eq: KKT gamma directed}). The last $p+1$ entries of the first term on the right-hand side are $\hat \Theta_\xi P_n \nabla_\xi l_{\theta_0}$. For the second term on the right hand side, notice that
\[
\frac{1}{N} D^T \hat W^2 D = \frac{1}{N}\begin{bmatrix}
	X^T\hat W^2X & X^T\hat W^2\textbf{1} & X^T\hat W^2Z \\
	\textbf{1}^T\hat W^2X & \textbf{1}^T\hat W^2\textbf{1} & \textbf{1}^T\hat W^2Z \\
	Z^T\hat W^2X & Z^T\hat W^2\textbf{1} & Z^T\hat W^2Z
\end{bmatrix}.
\]
$\hat \Theta_\xi$ is the exact inverse of $\hat \Sigma_\xi$ which is the lower-right $(p+1) \times (p+1)$ block of above matrix. Thus,
\[
A \frac{1}{N} D^T \hat W^2 D = \begin{bmatrix}
	\textbf{0} &\textbf{0} \\
	\hat \Theta_\xi \frac{1}{N}D_\xi^T\hat W^2X & I_{(p+1) \times (p+1)}
\end{bmatrix}.
\]
Then, for the last $p+1$ entries of $A \frac{1}{N} D^T \hat W^2 D (\hat{ \theta} - \theta_0)$
\[
\left(A \frac{1}{N} D^T \hat W^2 D (\hat{ \theta} - \theta_0) \right)_{\text{last }p+1\text{ entries}}= \hat{ \Theta}_\xi \frac{1}{N}D_\xi^T\hat W^2X (\hat{ \vartheta} - \vartheta_0) + \begin{pmatrix}
	\hat \mu - \mu_0 \\
	\hat \gamma - \gamma_0
\end{pmatrix}.
\]
Thus, (\ref{Eq: Taylor P_n directed}) implies
\begin{equation*}
	0 = \hat \Theta_\xi P_n \nabla_\gamma l_{\theta_0} +
	\hat{ \Theta}_\xi \frac{1}{N}D_\xi^T\hat W^2X (\hat{ \vartheta} - \vartheta_0)
	+ \begin{pmatrix}
		\hat \mu - \mu_0 \\
		\hat \gamma - \gamma_0
	\end{pmatrix}
	+ O\left( \hat \Theta_\xi \frac{1}{N}\sum_{i \neq j} \begin{pmatrix}
		1 \\
		Z_{ij}
	\end{pmatrix} \vert D_{ij}^T(\hat \theta - \theta_0) \vert^2 \right),
\end{equation*}
which is equivalent to
\begin{equation}\label{Eq: inference equation directed}
	\begin{pmatrix}
		\hat \mu - \mu_0 \\
		\hat \gamma - \gamma_0
	\end{pmatrix} = - \hat \Theta_\xi P_n \nabla_\xi l_{\theta_0} -
	\hat{ \Theta}_\xi \frac{1}{N}D_\xi^T\hat W^2X (\hat{ \vartheta} - \vartheta_0)
	+ O\left( \hat \Theta_\xi \frac{1}{N}\sum_{i \neq j} \begin{pmatrix}
		1 \\
		Z_{ij}
	\end{pmatrix} \vert D_{ij}^T(\hat \theta - \theta_0) \vert^2 \right).
\end{equation}
Our goal is now to show that for each component $k = 1, \dots, p+1$,
\[
\sqrt{N}\frac{\hat \xi_k - \xi_{0,k}}{\sqrt{\hat \Theta_{\xi,k,k}}} \overset{d}{\longrightarrow} \mathcal{N}(0,1).
\]
as described in the \textbf{Goal} section.
To that end, by equation (\ref{Eq: inference equation directed}), we now need to solve the following three problems: Writing $\hat{ \Theta}_{\xi,k}$ for the $k$-th row of $\hat{ \Theta}_\xi$,
\begin{enumerate}
	\item $\sqrt{N} \frac{ \hat \Theta_{\xi,k} P_n \nabla_\xi l_{\theta_0}}{\sqrt{\hat \Theta_{\xi,k,k}}} \overset{d}{\longrightarrow}  \mathcal{N}(0,1)$,
	\item 
	$
	\frac{1}{\sqrt{\hat \Theta_{\xi,k,k}}} \hat{ \Theta}_{\xi,k} \frac{1}{N}D_\xi^T\hat W^2X (\hat{ \vartheta} - \vartheta_0) = o_P\left( N^{-1/2} \right),
	$
	\item $
	O\left( 	\frac{1}{\sqrt{\hat \Theta_{\xi,k,k}}} \hat \Theta_{\xi,k} \frac{1}{N}\sum_{i \neq j} \begin{pmatrix}
		1 \\
		Z_{ij}
	\end{pmatrix} \vert D_{ij}^T(\hat \theta - \theta_0) \vert^2 \right) = o_P\left( N^{-1/2} \right).
	$
\end{enumerate}

\subsection{Bounding inverses}\label{Sec: Bounding inverses directed}

The problems (1) - (3) above suggest that it will be essential to bound the norm and the distance of $\hat \Theta_\xi$ and $\Theta_\xi$ in an appropriate manner. Notice that for any invertible matrices $A, B \in \R^{m \times m}$ we have
\[
A^{-1} - B^{-1} = A^{-1}(B - A)B^{-1}.
\]
Thus, for any sub-multiplicative matrix norm $\Vert \; . \; \Vert$, we get
\begin{equation}\label{Eq: difference between inverse matrices}
	\Vert A^{-1} - B^{-1} \Vert \le  \Vert A^{-1} \Vert \Vert B^{-1} \Vert \Vert B - A \Vert.
\end{equation}
We are particularly interested in the matrix $\infty$-norm, defined as
\[
\Vert A \Vert_\infty \coloneqq \sup \left\{ \frac{\Vert Ax \Vert_\infty}{\Vert x \Vert_\infty}, x \neq 0  \right\} =  \sup \left\{ \Vert Ax \Vert_\infty, \Vert x \Vert_\infty = 1 \right\}  = \max_{1 \le i \le m} \sum_{j = 1}^m \vert A_{i,j} \vert,
\]
i.e. $\Vert A \Vert_\infty$ is the maximal row $\ell_1$-norm of $A$. It is well-known, that any such matrix norm induced by a vector norm is sub-multiplicative ($\Vert AB \Vert_\infty \le \Vert A \Vert_\infty \Vert B \Vert_\infty$) and consistent with the inducing vector norm ($\Vert Ax \Vert_\infty \le \Vert A \Vert_\infty \Vert x \Vert_\infty$ for any vector $x$ of appropriate dimension). We first want to bound the matrix $\infty$-norm in terms of the largest eigenvalue.

\begin{lemma}\label{Lem: bound matrix norm by eval}
	For any symmetric, positive semi-definite $(m \times m)$-matrix $A$ with maximal eigenvalue $\lambda > 0$, we have $\Vert A \Vert_\infty \le \sqrt{m} \lambda$.
\end{lemma}
\begin{proof}
	\begin{align*}
		\Vert A \Vert_\infty &=  \sup \left\{ \Vert Ax \Vert_\infty, \Vert x \Vert_\infty = 1 \right\}  \\
		&\le \sup \left\{ \Vert Ax \Vert_2, \Vert x \Vert_\infty = 1 \right\}, \quad \Vert Ax \Vert_\infty \le \Vert Ax \Vert_2 \\
		&= \sup \left\{ \frac{\Vert Ax \Vert_2}{\Vert x \Vert_2} \Vert x \Vert_2, \Vert x \Vert_\infty = 1   \right\} \\
		&\le \sqrt{m} \sup \left\{ \frac{\Vert Ax \Vert_2}{\Vert x \Vert_2}, \Vert x \Vert_\infty = 1   \right\}, \quad \text{ if } \Vert x \Vert_\infty = 1, \text{ then } \Vert x \Vert_2 \le \sqrt{m}, \\
		&\le  \sqrt{m} \sup \left\{ \frac{\Vert Ax \Vert_2}{\Vert x \Vert_2}, x \neq 0  \right\} \\
		&= \sqrt{m} \Vert A \Vert_2 = \sqrt{m}\lambda,
	\end{align*}
	where $\Vert A \Vert_2$ is the spectral norm of the matrix $A$ and we have used that for symmetric matrices, the spectral norm is equal to the modulus of the largest eigenvalue of $A$.
\end{proof}
Also, recall that the inverse of a symmetric matrix $A$ is itself symmetric:
\[
I = A A^{-1} = A^T A^{-1} \overset{\text{transpose}}{\Longrightarrow} I = (A^{-1})^TA^T \overset{\text{symmetry}}{=} (A^{-1})^TA \overset{\text{uniqueness of inverse}}{\Longrightarrow} (A^{-1})^T = A^{-1}.
\]
Hence, $\hat \Theta_\xi$ and $\Theta_\xi$ are symmetric and we may apply Lemma \ref{Lem: bound matrix norm by eval}. Using that $\lambda_{\max}(\Sigma_\xi^{-1}) = \frac{1}{\lambda_{\min}(\Sigma_\xi)}$, we get
\begin{equation*}
	\Vert \Theta_\xi \Vert_\infty \le \sqrt{p} \lambda_{\max}(\Sigma_\xi^{-1}) \le C \frac{1}{\rho_n},
\end{equation*}
and with high probability
\begin{equation*}
	\Vert \hat \Theta_\xi \Vert_\infty \le \sqrt{p} \lambda_{\max}(\hat \Sigma_\xi^{-1}) \le C \frac{1}{\rho_n},
\end{equation*}
with some absolute constant $C$. 
Finally, by (\ref{Eq: difference between inverse matrices}),
\begin{equation*}
	\Vert \hat \Theta_\xi - \Theta_\xi \Vert_\infty \le \Vert \hat \Theta_\xi \Vert_\infty \Vert \Theta_\xi \Vert_\infty \Vert \hat \Sigma_\xi - \Sigma_\xi \Vert_\infty \le \frac{C}{\rho_n^2} \Vert \hat \Sigma_\xi - \Sigma_\xi \Vert_\infty.
\end{equation*}
It remains to control $\Vert \hat \Sigma_\xi - \Sigma_\xi \Vert_\infty$. We have
\begin{align*}
	\hat \Sigma_\xi - \Sigma_\xi  &= \frac{1}{N} \left( D_\xi^T\hat W^2 D_\xi - \E[D_\xi^T W_0^2 D_\xi]  \right) \\
	&= \underbrace{\frac{1}{N} \left( D_\xi^T (\hat W^2 - W_0^2) D_\xi   \right)}_{({I})} + \underbrace{\frac{1}{N} \left( D_\xi^TW_0^2 D_\xi - \E[D_\xi^TW_0^2 D_\xi]  \right)}_{(II)}.
\end{align*}
Recall that $\hat w_{ij}^2 = p_{ij}(\hat \theta) (1 - p_{ij}(\hat \theta)) = \frac{\exp(D_{ij}^T\hat \theta)}{(1 + \exp(D_{ij}^T\hat{ \theta}))^2} = \partial_{a^2}l(A_{ij}, D_{ij}^T\hat{ \theta})$, with the function $l$ defined in (\ref{Eq: Def l (a)}). Also recall that $\partial_{a^2}l$ is Lipschitz with constant one, by the Mean Value Theorem and the fact that it has derivative $\partial_{a^3}l$ bounded by one. Thus, considering the $(k,l)$-th element of $(I)$ above, we get:
\begin{align*}
	\left| \frac{1}{N} \left( D_\xi^T (\hat W^2 - W_0^2) D_\xi   \right)_{kl} \right| &= \left| \frac{1}{N} \sum_{ i \neq j} D_{ij,n+k}D_{ij,n+l} (\hat w_{ij}^2 - w_{0,ij}^2)\right| \\
	&\le C \frac{1}{N} \sum_{ i \neq j} \vert \hat w_{ij}^2 - w_{0,ij}^2\vert, \quad \text{ by unifrom boundedness of } Z_{ij} \\
	&\le C \frac{1}{N} \sum_{ i \neq j} \vert D_{ij}^T(\hat \theta - \theta_0) \vert, \quad \text{ by Lipschitz continuity} \\
	&\le \frac{C}{N} \sum_{ i \neq j} \left\{  \vert \hat \alpha_i - \alpha_{0,i} \vert + \vert \hat \beta_j - \beta_{0,j} \vert + \vert \hat \mu - \mu_0 \vert + \vert Z_{ij}^T(\hat \gamma - \gamma_0) \vert \right\} \\
	&\le \frac{C}{N} \underbrace{ \left\{\sum_{ i \neq j}  \vert \hat \alpha_i - \alpha_{0,i} \vert + \vert \hat \beta_j - \beta_{0,j} \vert \right\}}_{= (n-1) \Vert \hat \vartheta - \vartheta_0 \Vert_1}  + C \vert \hat \mu - \mu_0 \vert  + C \Vert \hat \gamma - \gamma_0 \Vert_1 \\
	&\le C \left\{  \frac{1}{n} \Vert \hat \vartheta - \vartheta_0 \Vert_1 + \vert \hat \mu - \mu_0 \vert + \Vert \hat \gamma - \gamma_0 \Vert_1  \right\} \\
	&= O_P\left( s^*_+ \sqrt{\frac{\log(n)}{N}} \rho_n^{-1}  \right), \, \text{ under the conditions of theorem \ref{Thm: consistency directed}}.
\end{align*} 
Since the dimension of $(I)$ is $(p+1) \times (p+1)$ and thus remains fixed, any row of $(I)$ has $\ell_1$ norm of order $O_P\left( s^*_+ \sqrt{\frac{\log(n)}{N}} \rho_n^{-1}  \right)$ and thus
\[
\Vert(I) \Vert_\infty = O_P\left( s^*_+ \sqrt{\frac{\log(n)}{N}} \rho_n^{-1}  \right).
\]
Taking a look at the $(k,l)$-th element in $(II)$:
\begin{align*}
	\left\vert \frac{1}{N} \left( D_\xi^TW_0^2 D_\xi - \E[D_\xi^TW_0^2 D_\xi]  \right)_{kl} \right\vert = \left\vert \frac{1}{N} \sum_{ i \neq j} \left\{D_{ij,n+k}D_{ij,n+l}w_{0,ij}^2 - \E[D_{ij,n+k}D_{ij,n+l} w_{0,ij}^2] \right\}\right\vert.
\end{align*}
Note that the random variables $D_{ij,n+k}D_{ij,n+l}w_{0,ij}^2$ are bounded uniformly in $i,j,k,l$. Thus, by Hoeffding's inequality, for any $t \ge 0$,
\begin{align*}
	P\left( \left\vert \frac{1}{N} \sum_{ i \neq j} \left\{D_{ij,n+k}D_{ij,n+l}w_{0,ij}^2 - \E[D_{ij,n+k}D_{ij,n+l} w_{0,ij}^2] \right\}\right\vert \ge t  \right) \le 2 \exp\left( - C N t^2 \right).
\end{align*}
This means, $\left\vert \frac{1}{N} \left( D_\xi^TW_0^2 D_\xi - \E[D_\xi^TW_0^2 D_\xi]  \right)_{kl} \right\vert = O_P\left( N^{-1/2}  \right)$. Again, since the dimension $p+1$ is fixed, we get by a simple union bound
\[
\Vert (II) \Vert_\infty = O_P\left(  N^{-1/2} \right).
\]
In total, we thus get
\begin{align*}
	\Vert \hat \Sigma_\xi - \Sigma_\xi \Vert_\infty 
	&=O_P\left(  s^*_+ \sqrt{\frac{\log(n)}{N}} \rho_n^{-1}  + \frac{1}{\sqrt{N}}  \right) = O_P\left( s^*_+ \sqrt{\frac{\log(n)}{N}} \rho_n^{-1} \right).
\end{align*}
We can now obtain a rate for $	\Vert \hat \Theta_\xi - \Theta_\xi \Vert_\infty$.
\begin{align*}
	\Vert \hat \Theta_\xi - \Theta_\xi \Vert_\infty \le \frac{C}{\rho_n^2} \Vert \hat \Sigma_\xi - \Sigma_\xi \Vert_\infty = O_P\left(  s^*_+ \sqrt{\frac{\log(n)}{N}} \rho_n^{-3} \right).
\end{align*}
By assumption \ref{Assum: new rate of s and rho_n directed}, we have $s_+^*\frac{\sqrt{\log(n)}}{\sqrt{n}\rho_n^2} \rightarrow 0, n \rightarrow \infty$, which in particular also implies that the above is $o_P(1)$. Notice in particular, that we have now managed to get for $k = 1, \dots, p+1,$
\begin{itemize}
	\item $\Vert \hat \Theta_{\xi,k}  - \Theta_{\xi,k} \Vert_1 = o_P(1)$,
	\item $\hat \Theta_{\xi,k,k} = \Theta_{\xi,k,k} + o_p(1)$.
\end{itemize}

\subsection{Problem 1}

We can now take a look at the problems (1) - (3) outlined above. For problem (1), we want to show:
\[
\sqrt{N} \frac{ \hat \Theta_{\xi,k} P_n \nabla_\xi l_{\theta_0}}{\sqrt{\hat \Theta_{\xi,k,k}}} \rightarrow \mathcal{N}(0,1).
\]
\textbf{Step 1:} Show that
\begin{equation}\label{Eq: root n}
	\hat{ \Theta}_{\xi,k} P_n \nabla_\xi l_{\theta_0} = \Theta_{\xi,k} P_n \nabla_\xi l_{\theta_0}+ o_P\left( N^{-1/2}\right).
\end{equation}
We have
\begin{align*}
	\vert (\hat \Theta_{\xi,k} - \Theta_{\xi,k} ) P_n \nabla_\xi l_{\theta_0} \vert &\le \Vert \hat \Theta_{\xi,k} - \Theta_{\xi,k} \Vert_1 \left\Vert \frac{1}{N} \sum_{ i \neq j} \begin{pmatrix}
		1 \\
		Z_{ij}
	\end{pmatrix} (p_{ij}(\theta_0) - A_{ij}) \right\Vert_\infty \\
	&\le \Vert \hat \Theta_{\xi} - \Theta_{\xi} \Vert_\infty \left\Vert \frac{1}{N} \sum_{ i \neq j} D_{\xi, ij} (p_{ij}(\theta_0) - A_{ij}) \right\Vert_\infty.
\end{align*}
Consider the vector $\sum_{ i \neq j} D_{\xi, ij} (p_{ij}(\theta_0) - A_{ij}) \in \R^{p+1}$. The $k$-th component of it has the form $\sum_{ i \neq j} (p_{ij}(\theta_0) - A_{ij})$ for $k=1$ and $\sum_{ i \neq j} Z_{ij,k-1} (p_{ij}(\theta_0) - A_{ij}), k = 2, \dots, p+1$. Notice that for these components are all centred:
\[
\mathbb{E}[ D_{\xi, ij,k} (p_{ij}(\theta_0) - A_{ij})] = \mathbb{E}[D_{\xi,ij,k} \mathbb{E}[ (p_{ij}(\theta_0) - A_{ij})| Z_{ij}]] = \mathbb{E}[D_{\xi, ij, k} \cdot0] = 0,
\]
as well as $\vert D_{\xi, ij,k} (p_{ij}(\theta_0) - A_{ij}) \vert \le c$, where $c > 1$ is a universal constant bounding $\vert Z_{ij,k} \vert$ for all $i,j,k$.
Thus, by Hoeffding's inequality, for any $t > 0$,
\begin{align*}
	P\left( \left|   \frac{1}{N} \sum_{ i \neq j} D_{\xi, ij, k} (p_{ij}(\theta_0) - A_{ij})  \right| \ge t   \right) \le 2 \exp \left(- 2 \frac{N t^2 }{c^2}  \right)
\end{align*}
and thus, 
\[
\frac{1}{N} \sum_{ i \neq j} D_{\xi, ij} (p_{ij}(\theta_0) - A_{ij}) = O_P\left( N^{-1/2} \right).
\]
Since we have $\Vert \hat \Theta_Z - \Theta_Z \Vert_\infty = o_P(1)$, by section \ref{Sec: Bounding inverses directed}, step 1 is now concluded. 

\noindent \textbf{Step 2:} Show that

\begin{equation*}
	\hat \Theta_{\xi,k,k} = \Theta_{\xi,k,k} + o_P(1).
\end{equation*}
Since $\Vert \hat \Theta_\xi - \Theta_\xi \Vert_\infty = o_P(1)$, by section \ref{Sec: Bounding inverses directed}, for all $k$
\[
\vert \hat \Theta_{\xi,k,k} - \Theta_{\xi,k,k} \vert \le \Vert \hat \Theta_\xi - \Theta_\xi \Vert_\infty = o_P(1)
\]
and step 2 is concluded.

\noindent\textbf{Step 3:} Show that
\[
\left\vert \frac{1}{\Theta_{\xi,k,k}} \right\vert \le C < \infty,
\]
for some universal constant $C > 0$. Then, we may conclude from step 1 and step 2 that
\[
\sqrt{N} \frac{ \hat \Theta_{\xi,k} P_n \nabla_\xi l_{\theta_0}}{\sqrt{\hat \Theta_{\xi,k,k}}} = \sqrt{N}\frac{\Theta_{\xi,k} P_n \nabla_\xi l_{\theta_0}}{\sqrt{\Theta_{\xi,k,k}}} + o_P(1). 
\]
To prove step 3, notice that $\Theta_\xi$ is symmetric and hence has only real eigenvalues. Therefore it is unitarily diagonalizable and for any $x \in \R^{p+1}$, we have $x^T \Theta_\xi x \ge \lambda_{\min}( \Theta_\xi) \Vert x \Vert_2^2$. We also know that
\[
\lambda_{\min}({  \Theta_\xi}) = \frac{1}{\lambda_{\max} ( \Sigma_\xi)}.
\]
Under assumption \ref{Assum: maximum EW directed} we can now deduce an upper bound on the maximum eigenvalue of $\Sigma_\xi$: For any $x \in \R^p$,
\[
x^T\Sigma_\xi x = x^T \frac{1}{N}\E[D_\xi^T W_0^2 D_\xi] x \le x^T \frac{1}{N}\E[D_\xi^T D_\xi] x \le (1 \vee \lambda_{\max}) \Vert x \Vert_2^2,
\]
where we have used that any entry in $W_0^2$ is bounded above by one. Since $x^T\Sigma_\xi x \le \lambda_{\max}(\Sigma_\xi) \Vert x \Vert_2^2$ and since this bound is tight (we have equality if $x$ is an eigenvector corresponding to $\lambda_{\max}$), we can conclude by assumption \ref{Assum: maximum EW directed} that $\lambda_{\max}(\Sigma_\xi) \le (1 \vee \lambda_{\max}) \le C < \infty$ for some universal constant $C > 0$.

In particular, since ${ \Theta}_{\xi,k,k} = e_k^T \Theta_\xi e_k$, we get
\[
{ \Theta}_{\xi,k,k} \ge \lambda_{\min}({ \Theta_\xi}) \Vert e_k \Vert_2^2 = \frac{1}{\lambda_{\max} ( \Sigma_\xi)} \ge C > 0,
\]
uniformly for all $n$. Consequently,
\[
0 <  \frac{1}{\Theta_{\xi,k,k}}  \le C < \infty.
\]
Step 3 is thus concluded.

\textbf{Step 4:} Finally, show that
\begin{equation*}
	\sqrt{N}\frac{\Theta_{\xi,k} P_n \nabla_\xi l_{\theta_0}}{\sqrt{\Theta_{\xi,k,k}}} \overset{d}{\longrightarrow} \mathcal{N}(0,1),
\end{equation*}
Such that by all the above
\begin{equation*}
	\sqrt{N} \frac{ \hat \Theta_{\xi,k} P_n \nabla_\xi l_{\theta_0}}{\sqrt{\hat \Theta_{\xi,k,k}}} \overset{d}{\longrightarrow} \mathcal{N}(0,1).
\end{equation*}
For brevity, we write $p_{ij}$ for the true link probabilities $p_{ij}(\theta_0)$. Also keep in mind that $\Theta_{\xi,k}$ denotes the $k$-th \textit{row} of $\Theta_\xi$, while $D_{\xi, ij}$ denote $((p+1) \times 1)$-\textit{column} vectors.
We want to apply the Lindeberg-Feller Central Limit Theorem. The random variables we study are the summands in 
\[
\sqrt{N}  \Theta_{\xi,k} P_n \nabla_\xi l_{\theta_0} = \sum_{ i \neq j} \left\{ \frac{1}{\sqrt{N}} \Theta_{\xi,k} D_{\xi, ij} (p_{ij} - A_{ij}) \right\}.
\]
First, notice that these random variables are centred:
\[
\mathbb{E}\left[ \frac{1}{\sqrt{N}} \Theta_{\xi,k}D_{\xi, ij} (p_{ij} - A_{ij}) \right] = \mathbb{E}\left[ \frac{1}{\sqrt{N}} \Theta_{\xi,k}D_{\xi, ij} \mathbb{E}[p_{ij} - A_{ij} | Z_{ij}]  \right] = \mathbb{E}\left[ \frac{1}{\sqrt{N}} \Theta_{\xi,k}D_{\xi, ij}  \cdot 0 \right] = 0.
\]
For the Lindeberg-Feller CLT we need to sum up the variances of these random variables. 
We claim that
\[
\sum_{ i \neq j} \text{Var}\left( \frac{1}{\sqrt{N}} \Theta_{\xi,k}D_{\xi, ij} (p_{ij} - A_{ij})   \right) = \Theta_{\xi,k,k}.
\]
Indeed, consider the vector-valued random variable $\sum_{ i \neq j} \left\{ \frac{1}{\sqrt{N}} D_{\xi, ij} (p_{ij} - A_{ij}) \right\} \in \R^{p+1}$. It has covariance matrix
\begin{align*}
	\mathbb{E}&\left[ \sum_{ i \neq j} \left\{ \frac{1}{\sqrt{N}} D_{\xi, ij} (p_{ij} - A_{ij}) \right\} \sum_{ i \neq j} \left\{ \frac{1}{\sqrt{N}} D_{\xi, ij} (p_{ij} - A_{ij}) \right\}^T   \right] \\
	&= \mathbb{E}\left[ \sum_{ i \neq j} \frac{1}{\sqrt{N}} D_{\xi, ij} (p_{ij} - A_{ij}) \frac{1}{\sqrt{N}} D_{\xi, ij}^T (p_{ij} - A_{ij})   \right], \quad \text{by independence accross } i,j \\
	&= \frac{1}{N} \sum_{i \neq j} \left[  \mathbb{E}[D_{\xi, ij,k}D_{\xi, ij,l}  (p_{ij} - A_{ij})^2] \right]_{k,l = 1, \dots, p+1}, \quad \text{ as a } ((p+1) \times (p+1))\text{-matrix} \\
	&= \frac{1}{N} \mathbb{E}[D_\xi^TW_0^2D_\xi] \\
	&= \Sigma_\xi.
\end{align*}
Thus, by independence across $i,j$,
\[
\sum_{ i \neq j} \text{Var}\left( \frac{1}{\sqrt{N}} \Theta_{\xi,k}D_{\xi,ij} (p_{ij} - A_{ij})   \right) = \text{Var}\left( \Theta_{\xi,k}  \sum_{ i \neq j}  \frac{1}{\sqrt{N}} D_{\xi, ij} (p_{ij} - A_{ij})  \right) = \Theta_{\xi,k} \Sigma_\xi \Theta_{\xi,k}^T = \Theta_{\xi,k,k},
\]
where for the last equality we have used that $\Theta_\xi$ is the inverse of $\Sigma_\xi$ and thus, $\Sigma_\xi\Theta_{\xi,k}^T= e_k$.
Now, we need to show that the Lindeberg condition holds. That is, we want that for any $\epsilon > 0$,
\begin{equation}\label{Eq: Lindeberg condition}
	\lim_{n\rightarrow \infty} \frac{1}{\Theta_{\xi,k,k}} \sum_{ i \neq j} \mathbb{E}\left[ \left\{  \frac{1}{\sqrt{N}} \Theta_{\xi,k}D_{\xi, ij} (p_{ij} - A_{ij}) \right\}^2 \mathbbm{1}\left( \vert \Theta_{\xi,k}D_{\xi, ij} (p_{ij} - A_{ij}) \vert > \epsilon \sqrt{N\Theta_{\xi,k,k} } \right) \right] = 0.
\end{equation}
We have
\[
\vert \Theta_{\xi,k}D_{\xi, ij} (p_{ij} - A_{ij}) \vert \le p \cdot c \cdot \Vert \Theta_{\xi,k} \Vert_1 \le C \Vert \Theta_\xi \Vert_\infty \le C \rho_n^{-1}.
\]
At the same time, we know from step 3 that $\Theta_{Z,k,k} \ge C > 0$ for some universal $C$. Then, as long as $\rho_n^{-1}$ goes to infinity at a rate slower than $n$, which is enforced by assumption \ref{Assum: new rate of s and rho_n directed}, we must have for $n$ large enough
\[
\vert \Theta_{\xi,k}D_{\xi, ij} (p_{ij} - A_{ij}) \vert < \epsilon \sqrt{N\Theta_{\xi,k,k} }
\]
uniformly in $i,j$. Thus, the indicator function and therefore each summand in (\ref{Eq: Lindeberg condition}) is equal to zero for $n$ large enough. Hence, (\ref{Eq: Lindeberg condition}) holds. Then, by the Lindeberg-Feller CLT,
\[
\sqrt{N}\frac{\Theta_{\xi,k} P_n \nabla_\xi l_{\theta_0}}{\sqrt{\Theta_{\xi,k,k}}} \overset{d}{\longrightarrow} \mathcal{N}(0,1).
\]
Now, by steps 1-4,
\begin{equation*}
	\sqrt{N} \frac{ \hat \Theta_{\xi,k} P_n \nabla_\xi l_{\theta_0}}{\sqrt{\hat \Theta_{\xi,k,k}}} \overset{d}{\longrightarrow} \mathcal{N}(0,1).
\end{equation*}
This concludes solving problem 1.

\subsection{Problem 2}\label{Sec: Problem 2 directed}
For problem 2 we must show
\[
\frac{1}{\sqrt{\hat \Theta_{\xi,k,k}}} \hat{ \Theta}_{\xi,k} \frac{1}{N}D_\xi^T\hat W^2X (\hat{ \vartheta} - \vartheta_0) = o_P\left( N^{-1/2} \right).
\]
Since we have $\Vert \hat \Theta_\xi - \Theta_\xi \Vert_\infty = o_P(1)$, we do not need to worry about $\frac{1}{\sqrt{\hat \Theta_{\xi,k,k}}}$, because $\hat{ \Theta}_{\xi,k,k} = \Theta_{\xi,k,k} + o_P(1)$ and $\frac{1}{\sqrt{ \Theta_{\xi,k,k}}} \le C < \infty$, i.e. $\frac{1}{\sqrt{ \hat \Theta_{\xi,k,k}}} = O_P(1)$ . By Theorem \ref{Thm: consistency directed} we also have a high-probability error bound on $\Vert \hat \vartheta - \vartheta_0 \Vert_1$. The problem will be bounding the corresponding matrix norms.
\[
\left\vert \hat{ \Theta}_{\xi,k} \frac{1}{N}D_\xi^T\hat W^2X (\hat{ \vartheta} - \vartheta_0) \right\vert \le \left\Vert \frac{1}{N} X^T\hat W^2 D_\xi \hat{ \Theta}_{\xi,k}^T\right\Vert_\infty \Vert \hat \vartheta - \vartheta_0 \Vert_1.
\]
Notice that in the display above we have the vector $\ell_\infty$-norm.
Also,
\[
\left\Vert\frac{1}{N} X^T\hat W^2 D_\xi \hat{ \Theta}_{\xi,k}^T \right\Vert_\infty \le \Vert \hat{ \Theta}_{\xi,k}^T \Vert_\infty \left\Vert \frac{1}{N}X^T\hat W^2D_\xi \right\Vert_\infty.
\]
Here we used the compatibility of the matrix $\ell_\infty$-norm with the vector $\ell_\infty$-norm. The first term is the vector norm, the second the matrix norm.
We know,
\[
\Vert \hat{ \Theta}_{\xi,k}^T \Vert_\infty \le \Vert \hat \Theta_\xi \Vert_\infty \le C \rho_n^{-1},
\]
where on the left hand side we have the vector norm and in the middle display the matrix norm.
Finally, $1/N \cdot X^T\hat W^2D_\xi$ is a $(n \times (p+1))$-matrix. The $(k,l)$-th element looks like
$1/N \cdot S_{k,l}$, where $S_{k,l}$ is the sum of $n-1$ terms of the form $D_{\xi, il,k} \hat{ w}_{il}^2$, summed over the appropriate indices $i,j$, all of which are uniformly bounded. Thus, 
\[
\left\vert \left( \frac{1}{N}X^T\hat W^2D_\xi \right)_{k,l} \right\vert \le \frac{1}{N} \cdot (n - 1)  \cdot c = \frac{C}{n}.
\]
Thus, the $\ell_1$-norm of any row of $\frac{1}{N}X^T\hat W^2D_\xi$ is bounded by $pC/n$ and thus
\[
\left\Vert \frac{1}{N}X^T\hat W^2D_\xi \right\Vert_\infty \le \frac{C}{n}.
\]
Recall that $\Vert \hat \vartheta - \vartheta_0 \Vert_1 = O_P\left( s^*_+ \frac{\sqrt{\log(n)}}{\sqrt{n}} \rho_n^{-1}  \right)$ by Theorem \ref{Thm: consistency directed}. Then,
\begin{align*}
	\left\vert \hat{ \Theta}_{\xi,k} \frac{1}{N}X^T\hat W^2D_\xi (\hat{ \vartheta} - \vartheta_0) \right\vert &\le \Vert \hat{ \Theta}_{\xi,k}^T \Vert_\infty \left\Vert \frac{1}{N}D_\xi^T\hat W^2X \right\Vert_\infty \Vert \hat \vartheta - \vartheta_0 \Vert_1 \\
	&=O_P\left( \frac{s^*_+}{\rho_n^2 \cdot n} \cdot  \frac{\sqrt{\log(n)}}{\sqrt{n}}  \right).
\end{align*}
Multiplying by $\sqrt{N} = O(n)$, gives
\begin{align*}
	\sqrt{N}\left\vert \hat{ \Theta}_{\vartheta,k} \frac{1}{N}D_\vartheta^T\hat W^2X (\hat{ \vartheta} - \vartheta_0) \right\vert &= O_P\left(  \frac{s^*_+}{\rho_n^2 } \cdot  \frac{\sqrt{\log(n)}}{\sqrt{n}}   \right),
\end{align*}
which is $o_P(1)$ under Assumption \ref{Assum: new rate of s and rho_n directed}.

\subsection{Problem 3}

Finally, we must show
\[
O\left( 	\frac{1}{\sqrt{\hat \Theta_{\xi,k,k}}} \hat \Theta_{\xi,k} \frac{1}{N}\sum_{i \neq j} \begin{pmatrix}
	1 \\
	Z_{ij}
\end{pmatrix} \vert D_{ij}^T(\hat \theta - \theta_0) \vert^2 \right) = o_P\left( N^{-1/2} \right).
\]
Again, since $\hat{ \Theta}_{\xi,k,k} = \Theta_{\xi,k,k} + o_P(1)$ and $\Theta_{\xi,k,k} \ge C > 0$ uniformly in $n$, we do not need to worry about the factor $	\frac{1}{\sqrt{\hat \Theta_{\xi,k,k}}}$ and it remains to show
\[
O\left( \hat \Theta_{\xi,k} \frac{1}{N}\sum_{i \neq j} D_{\xi, ij} \vert D_{ij}^\top(\hat \theta - \theta_0) \vert^2 \right) =  o_P\left( N^{-1/2} \right).
\]
We have
\begin{align*}
	\left\vert  \hat \Theta_{\xi,k} \frac{1}{N}\sum_{i \neq j} D_{\xi, ij} \vert D_{ij}^\top(\hat \theta - \theta_0) \vert^2  \right\vert &\le   \frac{1}{N}\sum_{i \neq j}  \vert \hat \Theta_{\xi,k}D_{\xi, ij} \vert \vert D_{ij}^T(\hat \theta - \theta_0) \vert^2 \\
	&\le c \Vert \hat \Theta_{\xi,k} \Vert_1 \frac{1}{N}\sum_{i \neq j} \vert D_{ij}^T(\hat \theta - \theta_0) \vert^2 \\
	&\le C \frac{1}{\rho_n} \frac{1}{N}\sum_{i \neq j} \vert D_{ij}^T(\hat \theta - \theta_0) \vert^2,
\end{align*}
where for the last inequality we have used that $\Vert \hat \Theta_{\xi,k} \Vert_1 \le \Vert \hat \Theta_\xi \Vert_\infty \le C \frac{1}{\rho_n}$. 
Now remember from \eqref{Eq: bound on sum D_ij theta -theta*} that
\[
\frac{1}{N}\sum_{i \neq j} \vert D_{ij}^T(\hat \theta - \theta_0) \vert^2 \le C \Vert \hat {\bar \theta} - \bar \theta_0 \Vert_1^2,
\]
where we make use of the fact that $\bar D \bar \theta = D \theta$.
From Theorem \ref{Thm: consistency directed} we know that under the assumptions of Theorem \ref{Thm: inference directed}, $\Vert \hat{ \bar \theta} - \bar \theta_0 \Vert_1 = O_P\left( s_{0,+} \sqrt{\frac{\log(n)}{N}} \rho_n^{-1} \right)$.
Thus,
\[
\sqrt{N}\left\vert  \hat \Theta_{\xi, k} \frac{1}{N}\sum_{i \neq j} D_{\xi, ij} \vert D_{ij}^T(\hat \theta - \theta_0) \vert^2  \right\vert = O_P\left( (s_{0,+})^2 \frac{\log(n)}{\sqrt{N}} \rho_n^{-3} \right).
\]
We see that this is $o_P(1)$ by applying assumption \ref{Assum: new rate of s and rho_n directed} twice. Problem 3 is solved.

\begin{proof}[of Theorem \ref{Thm: inference directed}]
	Theorem \ref{Thm: inference directed} now follows from the solved problems (1) - (3).
\end{proof}

\section{Proof of Theorem \ref{Thm: model selection consistency}}\label{Sec: model selection proofs}

\subsection{Proof of Lemmas}\label{Sec:auxiliary lemmas}

To make the representation cleaner, for the remainder of Section \ref{Sec: model selection proofs} we will simply write $S$ for $S_0$ and $S_+$ for $S_{0,+}$. Recall that we use $S_+^c$ to denote the complement of $S_+$ in $[2n+1+p]$, that is $S_+^c = [ 2n+1+p ] \backslash S_+$. We also use $S^c$ to refer to the complement of $S$ in $[2n]$ \textit{only}: $S^c = [2n] \backslash S$.

We begin by providing proofs of the lemmas in section \ref{Sec: model selection}.

\begin{proof}[of Lemma \ref{Lem: unique active set}]
	Since $\bar{\theta}^\dagger$ and $\hat{\bar{\theta}}$ both solve \eqref{Eq: Penalized llhd with covariates bar}, we must have
	\[
	\frac{1}{N} \bar{\mathcal{L}} ( \bar{\theta}^\dagger ) + \bar{\lambda} \Vert \bar{\vartheta}^\dagger \Vert_1 = \frac{1}{N} \bar{\mathcal{L}} (\hat{\bar \theta}) + \bar{\lambda} \Vert \hat{\bar{\vartheta}} \Vert_1.
	\]
	Denote by $\bar{z}^\dagger_{\vartheta}$ the first $2n$ components of $\bar{z}^{\dagger}$. Then, by \eqref{Eq: 22a} and \eqref{Eq: 22b}, $\langle \bar{z}^{\dagger }_{\vartheta} , \bar{\vartheta}^\dagger \rangle = \Vert \bar{\vartheta}^\dagger \Vert_1$. Thus,
	\begin{align*}
		\frac{1}{N} \bar{\mathcal{L}} ( \bar{\theta}^\dagger ) + \bar{\lambda} \langle \bar{z}^{\dagger }_{\vartheta} , \bar{\vartheta}^\dagger \rangle  = \frac{1}{N} \bar{\mathcal{L}} (\hat{\bar \theta}) + \bar{\lambda} \Vert \hat{\bar{\vartheta}} \Vert_1.
		\shortintertext{Hence, using that the last $p+1$ components of $\bar{z}^\dagger$ are zero,}
		\frac{1}{N} \bar{\mathcal{L}} ( \bar{\theta}^\dagger ) + \bar{\lambda} \langle \bar{z}^{\dagger } , \bar{\theta}^\dagger - \hat{\bar{\theta}} \rangle  = \frac{1}{N} \bar{\mathcal{L}} (\hat{\bar \theta}) + \bar{\lambda} \left(\Vert \hat{\bar{\vartheta}} \Vert_1 - \langle \bar{z}^\dagger, \hat{\bar{\theta}}\rangle \right).
		\shortintertext{But by \eqref{Eq: 21}, $\bar{\lambda }\bar{z}^\dagger = - 1/ N \cdot \nabla \bar{\mathcal{L}} (\bar{\theta}^\dagger) $ and therefore}
		\frac{1}{N} \bar{\mathcal{L}} ( \bar{\theta}^\dagger ) - \langle 1/ N \cdot \nabla \bar{\mathcal{L}} (\bar{\theta}^\dagger) , \bar{\theta}^\dagger - \hat{\bar{\theta}} \rangle - \frac{1}{N} \bar{\mathcal{L}} (\hat{\bar \theta}) = \bar{\lambda} \left(\Vert \hat{\bar{\vartheta}} \Vert_1 - \langle \bar{z}^\dagger, \hat{\bar{\theta}}\rangle \right).
	\end{align*}
	By the convexity of $\bar{\mathcal{L}}$, the left-hand side in the above display is negative. Therefore,
	\[
	\Vert \hat{\bar{\vartheta}} \Vert_1 \le \langle \bar{z}^\dagger, \hat{\bar{\theta}}\rangle = \langle \bar{z}^\dagger_\vartheta, \hat{\bar{\vartheta}} \rangle \le \Vert \bar{z}^\dagger_\vartheta \Vert_\infty \Vert \hat{\bar{\vartheta}} \Vert_1 \le \Vert \hat{\bar{\vartheta}} \Vert_1.
	\]
	Hence, $\langle \bar{z}^\dagger_\vartheta, \hat{\bar{\vartheta}} \rangle = \Vert \hat{\bar{\vartheta}} \Vert_1$. But since $\Vert \bar{z}^\dagger_{S^{\dagger c}} \Vert_\infty < 1$ by \eqref{Eq: correct exclusion}, this can only hold if $\hat{\bar{\vartheta}}_{S^{\dagger c}} = 0$. The claim follows.
\end{proof}

For the proof of Theorem \ref{Thm: model selection consistency} we need conditions similar to the ones in \cite{Ravikumar:2010}.
The first condition is the so-called \textit{dependency condition} which demands that the population Hessian of $\bar{\mathcal{L}}$ with respect to the variables contained in the active set $S$ is invertible. 
For our specific case, we let
\begin{equation}\label{Eq: Def of Q}
	Q \coloneqq \frac{1}{n-1} X^TW_0^2X = H_{\bar{\vartheta} \times \bar{\vartheta}} \bar{\mathcal{L}} (\bar{\theta})\in \R^{2n \times 2n},
\end{equation}
where $W_0 = \text{diag}(\sqrt{p_{ij}(\theta_0)(1-p_{ij}(\theta_0))}, i \neq j)$,
be the Hessian of $\bar{\mathcal{L}}$ with respect to  $\bar{\vartheta}$ only.
\begin{lemma}[Dependency condition]\label{Lem: dependency condition}
	For any $n$, the minimum eigenvalue of $Q_{S,S}$ satisfies
	\[
	\lambda_{min} \left( Q_{S,S} \right) \ge \frac{1}{2} \rho_n \cdot \left( 1 - \frac{\max\{s_\alpha, s_\beta\}}{n - 1}\right) > 0.
	\]
\end{lemma}

\begin{proof}[of Lemma \ref{Lem: dependency condition}]
	Notice that
	\[
	\frac{1}{n-1} X^TX = \begin{bmatrix}
		I_n & B \\
		B & I_n
	\end{bmatrix} \in \R^{2n \times 2n},
	\]
	where $I_n$ is the $(n \times n)$ identity matrix and $B$ is a matrix with zeros on the diagonal and $1/(n-1)$ everywhere else. Now consider the submatrix with only those rows and columns belonging to $S$
	\[
	P \coloneqq \frac{1}{n-1}( X^TX)_{S \times S} = \begin{bmatrix}
		I_{s_\alpha} & B_{S_\alpha, S_\beta} \\
		B_{S_\beta, S_\alpha} & I_{s_\beta}
	\end{bmatrix} \in \R^{s \times s}.
	\]
	This matrix $P$ is strictly diagonally dominant. Indeed,
	\begin{align*}
		\sum_{j \in S, j\neq i } P_{ij} &= \frac{s_\beta}{n-1} < 1 = P_{ii}, \quad i \in S_\alpha \\
		\sum_{j \in S, j\neq i } P_{ij} &= \frac{s_\alpha}{n-1} < 1 = P_{ii}, \quad i \in S_\beta,
	\end{align*}
	where the strict inequalities hold because $\min_i \{\alpha_{0,i}\} = \min_j \{\beta_{0,j}\} = 0$.
	Thus, $P$ is strictly positive definite. More, by the Gershgorin Circle Theorem, all the eigenvalues of $P$ must lie in one of the discs $D(P_{ii}, R_i)$, where $R_i = 	\sum_{j \in S, j\neq i } P_{ij}$ and $D(P_{ii}, R_i)$ is the disc with radius $R_i$ centred at $P_{ii}$. In particular,
	\[
	\text{mineval}(P) \ge 1 - \frac{\max \{s_\alpha, s_\beta\}}{n - 1}.
	\]
	But now, for any $v \in \R^s$,
	\begin{align*}
		v^T Q_{S,S} v \ge \frac{1}{2} \rho_n \cdot v^T P v \ge \frac{1}{2} \rho_n \left(  1 - \frac{\max \{s_\alpha, s_\beta\}}{n - 1} \right) \Vert v \Vert_2^2
	\end{align*}
	and the claim follows.
\end{proof}

\begin{lemma}[Incoherence condition]\label{Lem: incoherence condition}
	For any $n$,
	\[
	\Vert Q_{S^c,S} Q_{S,S}^{-1} \Vert_\infty \le \frac{1}{2} \rho_n^{-1} \cdot \frac{\max \{s_\alpha, s_\beta\}}{n - \max \{s_\alpha, s_\beta\}}.	
	\]
\end{lemma}
By Lemma \ref{Lem: dependency condition} the left-hand side of Lemma \ref{Lem: incoherence condition} is well-defined. Furthermore, under Assumption \ref{Assum: model selection assumption}, the right-hand side in Lemma \ref{Lem: incoherence condition} tends to zero as $n$ tends to infinity.

\begin{proof}[of Lemma \ref{Lem: incoherence condition}]
	We make use of the following bound of a the infinity norm of the inverse of a diagonally dominant matrix (see for example \cite{Varah:1975})
	\[
	\Vert Q_{S,S}^{-1} \Vert_\infty \le \max_{i \in S} \left\{\frac{1}{\vert q_{ii} \vert - R_i}\right\},
	\]
	where $q_{ii}$ is the $i$th diagonal entry of $Q_{S,S}$ and $R_i$ is the sum of the off-diagonal elements of the $i$th row of $Q_{S,S}$. That is, for $i \in S_\alpha$,
	\begin{equation*}
		q_{ii} - R_i = \frac{1}{n-1} \sum_{j=1, j\neq i}^n p_{ij}(1-p_{ij}) - \frac{1}{n-1} \sum_{j \in S_\beta } p_{ij}(1-p_ij) \ge \frac{1}{2(n-1)} \rho_n (n - s_\beta),
	\end{equation*}
	and analogously for $i \in S_\beta$,
	\begin{equation*}
		q_{ii} - R_i \ge  \frac{1}{2(n-1)} \rho_n (n - s_\alpha).
	\end{equation*}
	Thus,
	\begin{equation*}
		q_{ii} - R_i \ge  \frac{1}{2(n-1)} \rho_n (n - \max \{s_\alpha, s_\beta\})
	\end{equation*}
	and therefore,
	\begin{equation}\label{Eq: bound Q_SS inverse}
		\Vert Q_{S,S}^{-1} \Vert_\infty \le 2 \rho_n^{-1} \cdot \frac{n-1}{n - \max \{s_\alpha, s_\beta\}}.
	\end{equation}
	Furthermore, notice that any row of $Q_{S^c,S}$ has either $s_\alpha$ or $s_\beta$ non-zero entries, each of the form $1/(n-1) \cdot p_{ij}(1-p_{ij}) \le 1/(4(n-1))$. Hence,
	\[
	\Vert Q_{S^c,S} \Vert_\infty \le \frac{\max \{s_\alpha, s_\beta\}}{4(n-1)}.
	\]
	The claim follows by the submultiplicativity of the matrix infinity norm.
\end{proof}

\subsection{General strategy}

The proof of Theorem \ref{Thm: model selection consistency} hinges on the construction of $(\bar{\theta}^\dagger, \bar{z}^\dagger)$ succeeding with high probability and the challenge in proving this is proving that $(\bar{\theta}^\dagger, \bar{z}^\dagger)$ fulfils conditions \eqref{Eq: correct inclusion} and \eqref{Eq: correct exclusion}. Our proof relies on the following derivations. From \eqref{Eq: 21} we obtain
\begin{equation*}
	0 = \frac{1}{N}\nabla \bar{\mathcal{L}}(\constrtheta) + \bar{\lambda} \constrz - \frac{1}{N} \nabla \bar{\mathcal{L}}(\bar{\theta}_0) +  \frac{1}{N} \nabla \bar{\mathcal{L}}(\bar{\theta}_0).
\end{equation*}
Doing a Taylor expansion along the same lines as \eqref{Eq: Taylor of l} and \eqref{Eq: Lipschitz l a^2}, we obtain
\begin{align*}
	\frac{1}{N}\nabla \bar{\mathcal{L}}(\constrtheta) - \frac{1}{N} \nabla \bar{\mathcal{L}}(\bar{\theta}_0) =  \frac{1}{N} \bar{D}^T W_0^2 \bar{D} (\constrtheta - \bar{\theta}_0) + O\left(\frac{1}{N}\sum_{i \neq j} \bar{D}_{ij} \vert \bar{D}_{ij}^T( \constrtheta - \bar{\theta}_0) \vert^2 \right),
\end{align*}
where we have used the fact that we are taking derivatives with respect to $\bar{\theta}$ and used $\bar{D}_{ij}\bar{\theta}_0$ in \eqref{Eq: Lipschitz l a^2}, to obtain $W_0^2$ instead of $\hat{W}^2$ above. Combining the last two equations, we obtain
\begin{equation*}
	\frac{1}{N} \bar{D}^T W_0^2 \bar{D} (\constrtheta - \bar{\theta}_0) = - \bar{\lambda} \constrz -  \frac{1}{N} \nabla \bar{\mathcal{L}}(\bar{\theta}_0)+ O\left(\frac{1}{N}\sum_{i \neq j} \bar{D}_{ij} \vert \bar{D}_{ij}^T( \constrtheta - \bar{\theta}_0) \vert^2 \right).
\end{equation*}
Taking only the first $2n$ entries of that equation we obtain
\begin{align}\label{Eq: Taylor KKT}
	\begin{split}
		\frac{1}{N} \bar{X}^TW_0^2\bar{X} (\constrvartheta - \bar{\vartheta}_0) = &- \frac{1}{N} 	\nabla_{\bar{\vartheta}} \bar{\mathcal{L}}(\bar{\theta}_0) + \frac{1}{N} \bar{X}^TW_0^2\left[
		\begin{array}{c|c}
			\textbf{1} & Z
		\end{array}
		\right] (\xi^\dagger - \xi_0) - \bar{\lambda} \constrz_{1:2n} + \bar{R}
	\end{split}
\end{align}
where we use $\constrz_{1:2n}$ to refer to the first $2n$ components of $\constrz_{1:2n}$, use our shorthand notation $\xi = (\mu, \gamma^T)^T$ and let
\[
\bar{R} = O\left(\frac{1}{N}\sum_{i \neq j} 	\bar{X}_{ij} \vert \bar{D}_{ij}^T( \constrtheta - \bar{\theta}_0) \vert^2 \right).
\]
Notice that we the left-hand side in \eqref{Eq: Taylor KKT} is equal to
\[
Q (\constrvartheta - \bar{\vartheta}_0) = Q_{-,S} (\constrvartheta - \bar{\vartheta}_0)_S + Q_{-,S^c} \underbrace{(\constrvartheta - \bar{\vartheta}_0)_{S^c}}_{=0}.
\]
Plugging this into \eqref{Eq: Taylor KKT} and splitting up by rows, we get
\begin{subequations}
	\begin{align}
		Q_{S,S} (\constrvartheta - \bar{\vartheta}_0)_S = &- \frac{1}{N} 	\left(\nabla_{\bar{\vartheta}} \bar{\mathcal{L}}(\bar{\theta}_0)\right)_S+ \frac{1}{N} \bar{X}_S^TW_0^2\left[
		\begin{array}{c|c}
			\textbf{1} & Z
		\end{array}
		\right] (\xi^\dagger - \xi_0) - \bar{\lambda} \constrz_{1:2n, S} + \bar{R}_S \label{Eq: 31a} \\
		Q_{S^c,S} (\constrvartheta - \bar{\vartheta}_0)_{S} = &- \frac{1}{N} 	\left(\nabla_{\bar{\vartheta}} \bar{\mathcal{L}}(\bar{\theta}_0)\right)_{S^c} + \frac{1}{N} \bar{X}_{S^c}^TW_0^2\left[
		\begin{array}{c|c}
			\textbf{1} & Z
		\end{array}
		\right] (\xi^\dagger - \xi_0) - \bar{\lambda} \constrz_{1:2n, S^c} + \bar{R}_{S^c} \label{Eq: 31b},
	\end{align}
\end{subequations}
where it is important to remember that $S^c = [2n] \backslash S$ refers to the complement of $S$ in $[2n]$. We solve \eqref{Eq: 31a} for $(\bar{\vartheta}^\dagger - \bar{\vartheta}_0)_S$ and plug the result into \eqref{Eq: 31b}. Finally we rearrange for $- \bar{\lambda} \constrz_{1:2n, S^c}$,
\begin{align*}
	\begin{split}
		- \bar{\lambda} \constrz_{1:2n, S^c} = &Q_{S^c,S} Q_{S,S}^{-1} \bigg\{ - \frac{1}{N} 	\left(\nabla_{\bar{\vartheta}} \bar{\mathcal{L}}(\bar{\theta}_0)\right)_S+ \frac{1}{N} \bar{X}_S^TW_0^2\left[
		\begin{array}{c|c}
			\textbf{1} & Z
		\end{array}
		\right] (\xi^\dagger - \xi_0) - \bar{\lambda} \constrz_{1:2n, S} + \bar{R}_S  \bigg\} \\
		&+  \frac{1}{N} \left(\nabla_{\bar{\vartheta}} \bar{\mathcal{L}}(\bar{\theta}_0)\right)_{S^c} - \frac{1}{N} \bar{X}_{S^c}^TW_0^2\left[
		\begin{array}{c|c}
			\textbf{1} & Z
		\end{array}
		\right] (\xi^\dagger - \xi_0) - \bar{R}_{S^c}
	\end{split}
\end{align*}
Now, divide by $\bar{\lambda}$ and take the $\infty$-norm on both sides. Rearrange corresponding terms.
\begin{align*}
	\Vert \constrz_{1:2n, S^c} \Vert_\infty &\le \frac{1}{\bar{\lambda}} \left\{ \Vert Q_{S^c,S} Q_{S,S}^{-1} \Vert_\infty + 1  \right\} \left\Vert \frac{1}{N} 	\nabla_{\bar{\vartheta}} \bar{\mathcal{L}}(\bar{\theta}_0) \right\Vert_\infty & (I) \\
	& + \frac{1}{\bar{\lambda}} \left\{ \Vert Q_{S^c,S} Q_{S,S}^{-1} \Vert_\infty + 1 \right\} \left\Vert \bar{R} \right\Vert_\infty & (II) \\
	& + \frac{1}{\bar{\lambda}} \left\{ \Vert Q_{S^c,S} Q_{S,S}^{-1} \Vert_\infty + 1  \right\} \left\Vert \frac{1}{N} \bar{X}^TW_0^2\left[
	\begin{array}{c|c}
		\textbf{1} & Z
	\end{array}
	\right] (\xi^\dagger - \xi_0) \right\Vert_\infty& (III) \\
	& + \left\Vert Q_{S^c,S} Q_{S,S}^{-1} \right\Vert_\infty & (IV).
\end{align*}
By appropriately bounding the terms $(I) - (IV)$ on the right-hand side, we will proceed to show that for sufficiently large $n$, with high probability, $\Vert \constrz_{1:2n, S^c} \Vert_\infty < 1$, which is clearly equivalent to \eqref{Eq: correct exclusion}. Notice that we already may control term $(IV)$ as well as the terms $\Vert Q_{S^c,S} Q_{S,S}^{-1} \Vert_\infty + 1$ by the incoherence condition, Lemma \ref{Lem: incoherence condition}.

\subsection{Controlling term \texorpdfstring{$(I)$}{(I)}}

Notice that the $i$th component of $\frac{1}{N} \nabla_{\bar{\vartheta}} \bar{\mathcal{L}}(\bar{\theta}_0)$ is of the form
\[
\frac{1}{N} \sqrt{n} \sum_{j = 1, j \neq i} (A_{ij} - p_{ij}) = \frac{1}{\sqrt{n}} \cdot \frac{1}{n-1}  \sum_{j = 1, j \neq i} (A_{ij} - p_{ij}).
\]
In particular, each summand is a centred, bounded random variable. By Hoeffding's inequality, we have for every $t > 0$,
\[
P\left(  \left|    \frac{1}{n-1}  \sum_{j = 1, j \neq i} (A_{ij} - p_{ij})   \right|  \ge t \right)	 \le 2 \exp \left(   - \frac{n-1}{2} t^2\right).
\]
Thus, for any $\epsilon > 0$, picking $t = \epsilon \sqrt{n} \bar{\lambda}$, gives
\[
P\left(  \frac{1}{\bar{\lambda}}  \left|    \frac{1}{N}  \nabla_{\bar{\vartheta}} \bar{\mathcal{L}}(\bar{\theta}_0)_i  \right|  \ge \epsilon \right)	 \le 2 \exp \left(   - \frac{N \bar{\lambda}^2}{2} \epsilon^2\right).
\]
Taking a union bound over all $2n$ components of $\nabla_{\bar{\vartheta}} \bar{\mathcal{L}}(\bar{\theta}_0)$, leads to
\begin{equation}\label{Eq: bound gradient L wrt vartheta}
	P\left(  \frac{1}{\bar{\lambda}}  \left\Vert    \frac{1}{N}  \nabla_{\bar{\vartheta}} \bar{\mathcal{L}}(\bar{\theta}_0)  \right\Vert_\infty \ge \epsilon \right)	 \le 4n \cdot \exp \left(   - \frac{N \bar{\lambda}^2}{2} \epsilon^2\right) = 4 \cdot \exp \left(   - \frac{N \bar{\lambda}^2}{2} \epsilon^2 + \log(n) \right).
\end{equation}
In the next section, when controlling term $(II)$, we will also need a similar bound on the components of $ \frac{1}{N}  \nabla \bar{\mathcal{L}}(\bar{\theta}_0)$ corresponding to $\xi = (\mu, \gamma^T)^T$, which is why we derive the respective bounds now. Using analogous arguments to the above, we obtain
\begin{equation}\label{Eq: bound gradient L wrt xi}
	P\left(  \frac{1}{\bar{\lambda}}  \left\Vert    \frac{1}{N}  \nabla_{\xi} \bar{\mathcal{L}}(\bar{\theta}_0)  \right\Vert_\infty \ge \epsilon \right)	 \le 2(p+1) \cdot \exp \left(   - \frac{N \bar{\lambda}^2}{2 (1 \vee c^2)} \epsilon^2\right).
\end{equation}
Combining \eqref{Eq: bound gradient L wrt vartheta} and \eqref{Eq: bound gradient L wrt xi}, we obtain a bound on the infinity norm of the full gradient,
\begin{equation}\label{Eq: bound gradient}
	P\left(  \frac{1}{\bar{\lambda}}  \left\Vert    \frac{1}{N}  \nabla \bar{\mathcal{L}}(\bar{\theta}_0)  \right\Vert_\infty \ge \epsilon \right)	 \le 4 \cdot \exp \left(   - \frac{N \bar{\lambda}^2}{2} \epsilon^2 + \log(n) \right) + 2(p+1) \cdot \exp \left( - \frac{N \bar{\lambda}^2}{2 (1 \vee c^2)} \epsilon^2\right),
\end{equation}
which tends to zero, as long as $- \frac{N \bar{\lambda}^2}{2} \epsilon^2 + \log(n) \rightarrow \infty$, as $n$ tends to infinity.

\subsection{Controlling term \texorpdfstring{$(II)$}{(II)}}

Controlling term $(II)$ is by far the most involved step in controlling $\Vert \constrz_{1:2n, S^c} \Vert_\infty$. We start by controlling the $\ell_2$-error between our construction $\bar{\theta}^\dagger$ and the truth $\bar{\theta}_0$.

\begin{lemma}\label{Lem: bound construction error}
	Under assumptions \ref{Assum: minimum EW directed}, \ref{Assum: no approximation error},
	\ref{Assum: model selection assumption}, \ref{Assum: model selection lambda},
	for $n$ large enough, for any $\epsilon > 0$, with probability at least
	\begin{align*}
		\begin{split}
			1 &- 4 \cdot \exp \left(   - \frac{N \bar{\lambda}^2}{2} \epsilon^2 + \log(n) \right) - 2(p+1) \cdot \exp \left( - \frac{N \bar{\lambda}^2}{2 (1 \vee c^2)} \epsilon^2\right) \\
			&- p(p+3) \exp\left( - N \frac{c_{\min}^2}{2048 s^{2}_{+}\tilde{ c}}  \right),
		\end{split}
	\end{align*}
	which tends to one as long as $- \frac{N \bar{\lambda}^2}{2} \epsilon^2 + \log(n) \rightarrow - \infty$, as $n$ tends to infinity, we have
	\[
	\Vert \constrtheta - \bar{\theta}_0 \Vert_1 \le (1+\epsilon)\frac{9}{c_{\min} } \rho_n^{-1} s_+ \bar{\lambda}.
	\]
\end{lemma}

\begin{proof}
	Keep in mind that $\constrtheta - \bar{\theta}_0 = \constrtheta_{S_+} - \bar{\theta}_{0,S_+}$.
	Define a function $G: \R^{s + 1+ p} \rightarrow \R$,
	\[
	G(u) = \frac{1}{N} \left\{  \bar{\mathcal{L}} (\bar{\theta}_{0,S_+} + u) - \bar{\mathcal{L}} (\bar{\theta}_{0,S_+})  \right\} + \bar{\lambda} \left( \Vert \bar{\theta}_{0,S} + u_S \Vert_1 - \Vert \bar{\theta}_{0,S} \Vert_1 \right),
	\]
	where for the addition $\bar{\theta}_{0,S_+} + u$ to be well-defined, we use the canonical embedding of $\R^{s + 1 + p} \hookrightarrow \R^{2n + 1 + p}$, by setting the components not contained in $S$ to zero. In the following we will make use of that embedding without explicitly mentioning it if there is no chance of confusion. Also pay close attention to the distinction between $S_+$ and $S$ in above display. Clearly, $G(0) = 0$ and $G$ is minimized at $\bar{u}^\dagger = \constrtheta_{S_+} - \bar{\theta}_{0,S_+}$, which implies that $G(\bar{u}^\dagger) \le 0$. Also, $G$ is convex.
	
	Now suppose we manage to find some $B \in \R, B > 0$, such that for all $u \in \R^{s+1+p}$ with $\Vert u \Vert_1 = B$ it holds $G(u) > 0$. We claim that in that case it must hold $\Vert \bar{u}^\dagger \Vert_1 \le B$. Indeed, if $\Vert \bar{u}^\dagger \Vert_1 > B$, then there exists a $t \in (0,1)$ such that for $\tilde{u} = t \bar{u}^\dagger$ we have $\Vert \tilde{u} \Vert_1 = B$. But then, by convexity of $G$, $G(\tilde{u}) \le t G(\bar{u}^\dagger) + (1-t) G(0) = tG(\bar{u}^\dagger)\le 0$. A contradiction.
	
	Thus, we need to find an appropriate $B$. Let $B > 0$, the correct form to be determined later. Now, pick any $u \in \R^{s+1+p}$ with $\Vert u \Vert_1 = B$. We do a first order Taylor expansion of $\bar{\mathcal{L}}$ in the point $\bar{\theta}_{0,S_+}$, evaluated at $\bar{\theta}_{0,S_+} + u$. This yields
	\begin{align*}
		G(u) =& \frac{1}{N} \left\{  \nabla_{S_+} \bar{\mathcal{L}} (\bar{\theta}_{0,S_+})^T(\bar{\theta}_{0,S_+} + u - \bar{\theta}_{0,S_+}  ) + \frac{1}{2} \cdot u^TH_{S_+, S_+} \bar{\mathcal{L}}(\bar{\theta}_{0,S_+} + u \alpha) u \right\} \\
		&+ \bar{\lambda} \left( \Vert \bar{\theta}_{0,S} + u_S \Vert_1 - \Vert \bar{\theta}_{0,S} \Vert_1 \right),
	\end{align*}
	for some $\alpha \in [0,1]$.
	Now, using \eqref{Eq: bound gradient}, we know that with high-probability,
	\begin{equation}\label{Eq: 52}
		\left|  \frac{1}{N} \nabla_{S_+} \bar{\mathcal{L}} (\bar{\theta}_{0,S_+})^T u  \right| \le \left\Vert \frac{1}{N} \nabla_{S_+} \bar{\mathcal{L}} (\bar{\theta}_{0,S_+}) \right\Vert_\infty \Vert u \Vert_1 \le \epsilon \bar{\lambda} B
	\end{equation}
	with the $\epsilon$ from  \eqref{Eq: bound gradient}. Furthermore, by using the triangle inequality, we obtain
	\begin{equation}\label{Eq: 53}
		\bar{\lambda} \left( \Vert \bar{\theta}_{0,S} + u_S \Vert_1 - \Vert \bar{\theta}_{0,S} \Vert_1 \right) \ge - \bar{\lambda} \Vert u \Vert_1 = - \bar{\lambda} B
	\end{equation}
	Clearly, the canonical embedding of $u$ into $\R^{2n+1+p}$ fulfils the condition of the empirical compatibility condition, Proposition \ref{Prop: compatibility condition Sigma}. Also, keep in mind that assumptions \ref{Assum: model selection assumption} and \ref{Assum: model selection lambda} together imply $n^{-1/2} \rho_n^{-1}s_+ \rightarrow 0$, which in particular implies $s_+ = o(\sqrt{n})$. Thus, Proposition \ref{Prop: compatibility condition Sigma} is applicable and with high probability as prescribed in Proposition \ref{Prop: compatibility condition Sigma}, we have
	\begin{align}\label{Eq: 54}
		\begin{split}
			\frac{1}{2} \cdot u^TH_{S_+, S_+} \bar{\mathcal{L}}(\bar{\theta}_{0,S_+} + u \alpha) u &\ge \frac{1}{4} \rho_n u^T \left\{ \frac{1}{N} \bar{D}^T\bar{D} \right\}_{S_+,S_+} u \\
			&= \frac{1}{4} \rho_n u^T \Sigma u \\
			&\ge \frac{1}{8} \rho_n \frac{c_{\min}}{s_+} \Vert u \Vert_1^2 \\
			& = \frac{1}{8} \rho_n \frac{c_{\min}}{s_+} B^2
		\end{split}
	\end{align}
	Combining \eqref{Eq: 52}, \eqref{Eq: 53}, \eqref{Eq: 54}, we find
	\[
	G(u) \ge - \epsilon \bar{\lambda} B - \bar{\lambda} B + \frac{1}{8} \rho_n \frac{c_{\min}}{s_+} B^2.
	\]
	The right-hand side of this equation is strictly larger zero, whenever
	\[
	B > (1+\epsilon)\frac{8}{c_{\min} } \rho_n^{-1} s_+ \bar{\lambda}.
	\]
	Thus, the claim follows from picking
	\[
	B = (1+\epsilon)\frac{9}{c_{\min} } \rho_n^{-1} s_+ \bar{\lambda}.
	\]
\end{proof}

\begin{lemma}\label{Lem: bound on bar R}
	Under assumptions \ref{Assum: minimum EW directed}, \ref{Assum: no approximation error}, 
	\ref{Assum: model selection assumption}, \ref{Assum: model selection lambda},
	for $n$ large enough, for any $\epsilon > 0$, with probability at least
	\begin{align*}
		\begin{split}
			1 &- 4 \cdot \exp \left(   - \frac{N \bar{\lambda}^2}{2} \epsilon^2 + \log(n) \right) - 2(p+1) \cdot \exp \left( - \frac{N \bar{\lambda}^2}{2 (1 \vee c^2)} \epsilon^2\right) \\
			&- p(p+3) \exp\left( - N \frac{c_{\min}^2}{2048 s^{2}_+\tilde{ c}}  \right),
		\end{split}
	\end{align*}
	which tends to one as long as $- \frac{N \bar{\lambda}^2}{2} \epsilon^2 + \log(n) \rightarrow - \infty$, as $n$ tends to infinity, we have
	\[
	\frac{1}{\bar{\lambda}} \Vert \bar{R} \Vert_\infty \le 
	\frac{324 (1 \vee (c^2p)) (1+\epsilon)^2 }{c_{\min}^2 } \cdot \sqrt{n} \rho_n^{-2} s_+^2 \bar{\lambda}.
	\]
\end{lemma}

\begin{proof}
	Consider the $i$th component of $\bar{R}$, for $i \in S_\alpha$. Similar to \eqref{Eq: bound on sum D_ij theta -theta*} we obtain,
	\begin{align*}
		\bar{R}_i &= \frac{1}{N} \sum_{j = 1, j \neq i}^n 	\bar{X}_{ij} \vert \bar{D}_{ij}^T( \constrtheta - \bar{\theta}_0) \vert^2 \\
		&=  \frac{1}{\sqrt{n}} \cdot \frac{1}{n-1} \sum_{j = 1, j \neq i}^n \vert \bar{D}_{ij}^T( \constrtheta - \bar{\theta}_0) \vert^2 \\
		&= \frac{1}{\sqrt{n}} \cdot \frac{1}{n-1} \sum_{j = 1, j \neq i}^n \left( \alpha_i^\dagger - \alpha_{0,i} + \beta_j^\dagger - \beta_{0,j}^\dagger + \mu^\dagger - \mu_0 + Z_{ij}^T(\gamma^\dagger - \gamma_0) \right)^2 \\
		&\le \frac{4}{\sqrt{n}} \cdot \frac{1}{n-1} \sum_{j = 1, j \neq i}^n \left( ( \alpha_i^\dagger - \alpha_{0,i})^2 + (\beta_j^\dagger - \beta_{0,j})^2 +( \mu^\dagger - \mu_0)^2 + c^2p  \Vert \gamma^\dagger - \gamma_0\Vert_2^2 \right) \\
		&= \frac{4}{\sqrt{n}} \left\{ ( \alpha_i^\dagger - \alpha_{0,i})^2 +  ( \mu^\dagger - \mu_0)^2 +    c^2p  \Vert \gamma^\dagger - \gamma_0\Vert_2^2 \right\} + \frac{4}{\sqrt{n}} \cdot \frac{1}{n-1} \sum_{j = 1, j \neq i}^n (\beta_j^\dagger - \beta_{0,j})^2 \\
		&\le \frac{4}{\sqrt{n}} (1 \vee (c^2p)) \left\{ ( \alpha_i^\dagger - \alpha_{0,i})^2 +  ( \mu^\dagger - \mu_0)^2 +  \Vert \gamma^\dagger - \gamma_0\Vert_2^2 \right\} + \frac{\sqrt{n}}{n-1} \Vert \bar{\beta}^\dagger - \bar{\beta}_0 \Vert_2^2.
	\end{align*}
	We have
	\[
	( \alpha_i^\dagger - \alpha_{0,i})^2 = n ( \bar{\alpha}_i^\dagger - \bar{\alpha}_{0,i})^2 \le n \Vert \bar{\alpha}^\dagger - \bar{\alpha}_0 \Vert_2^2.
	\]
	Thus, by Lemma \ref{Lem: bound construction error}, with at least the prescribed probability and for all $i \in S_\alpha$,
	\begin{align*}
		\frac{R_i}{\bar{\lambda}} \le& 4 (1 \vee (c^2p)) \sqrt{n} \Vert \bar{\theta}^\dagger - \bar{\theta}_0 \Vert_2^2 \le 4 (1 \vee (c^2p)) \sqrt{n} \Vert \bar{\theta}^\dagger - \bar{\theta}_0 \Vert_1^2 \\
		&\le  \frac{324 (1 \vee (c^2p)) (1+\epsilon)^2 }{c_{\min}^2 } \cdot \sqrt{n} \rho_n^{-2} s_+^2 \bar{\lambda}.
	\end{align*}
	The same bound is found for all $i \in S_\beta$ using the exact same steps. Since the right-hand side above does not depend on $i$ the claim follows.
\end{proof}

\subsection{Controlling term \texorpdfstring{$(III)$}{(III)}}

\begin{lemma}\label{Lem: bound on model selection bias}
	Under assumptions \ref{Assum: minimum EW directed}, \ref{Assum: no approximation error}, 
	\ref{Assum: model selection assumption}, \ref{Assum: model selection lambda}
	for $n$ large enough, for any $\epsilon > 0$, with probability at least
	\begin{align*}
		\begin{split}
			1 &- 4 \cdot \exp \left(   - \frac{N \bar{\lambda}^2}{2} \epsilon^2 + \log(n) \right) - 2(p+1) \cdot \exp \left( - \frac{N \bar{\lambda}^2}{2 (1 \vee c^2)} \epsilon^2\right) \\
			&- p(p+3) \exp\left( - N \frac{c_{\min}^2}{2048 s^{2}_+\tilde{ c}}  \right),
		\end{split}
	\end{align*}
	which tends to one as long as $- \frac{N \bar{\lambda}^2}{2} \epsilon^2 + \log(n) \rightarrow - \infty$, as $n$ tends to infinity, we have
	\[
	\frac{1}{\bar{\lambda}} \left\Vert \frac{1}{N} \bar{X}^TW_0^2\left[
	\begin{array}{c|c}
		\textbf{1} & Z
	\end{array}
	\right] (\xi^\dagger - \xi_0) \right\Vert_\infty \le \frac{9 (1 \vee c) (1+\epsilon) (p+1)}{4c_{\min} } \cdot \frac{1}{\sqrt{n}} \rho_n^{-1} s_+.
	\]
\end{lemma}

\begin{proof}
	We have
	\begin{align*}
		\left\Vert \frac{1}{N} \bar{X}^TW_0^2\left[
		\begin{array}{c|c}
			\textbf{1} & Z
		\end{array}
		\right] (\xi^\dagger - \xi_0) \right\Vert_\infty &\le  \left\Vert \frac{1}{N} \bar{X}^TW_0^2\left[
		\begin{array}{c|c}
			\textbf{1} & Z
		\end{array}
		\right] \right\Vert_\infty \Vert  (\xi^\dagger - \xi_0) \Vert_\infty \\
		&\le \left\Vert \frac{1}{N} \bar{X}^TW_0^2\left[
		\begin{array}{c|c}
			\textbf{1} & Z
		\end{array}
		\right] \right\Vert_\infty \Vert \bar{\theta}^\dagger - \bar{\theta}_0 \Vert_1.
	\end{align*}
	Consider the $i$th row of the matrix $\frac{1}{N} \bar{X}^TW_0^2\left[
	\begin{array}{c|c}
		\textbf{1} & Z
	\end{array}
	\right]$,
	\[
	\left\Vert \left( \frac{1}{N} \bar{X}_{-,i}^TW_0^2\left[
	\begin{array}{c|c}
		\textbf{1} & Z
	\end{array}
	\right] \right)^T \right\Vert_1 \le \frac{1}{N} \sqrt{n} (n-1) \cdot \frac{1}{4} (1 \vee c) (p+1) = \frac{p+1}{4} (1 \vee c) \frac{1}{\sqrt{n}},
	\]
	where we have used that the $i$th column of $\bar{X}$ has exactly $(n-1)$ non-zero entries, each with value $\sqrt{n}$, each entry of $W_0^2$ is upper bounded by $1/4$ and any row of $\left[
	\begin{array}{c|c}
		\textbf{1} & Z
	\end{array}
	\right]$ has $p+1$ entries, each of which is upper bounded by $1 \vee c$. Thus, by Lemma \ref{Lem: bound construction error}, with the prescribed probability,
	\[
	\frac{1}{\bar{\lambda}} \left\Vert \frac{1}{N} \bar{X}^TW_0^2\left[
	\begin{array}{c|c}
		\textbf{1} & Z
	\end{array}
	\right] (\xi^\dagger - \xi_0) \right\Vert_\infty \le \frac{9 (1 \vee c) (1+\epsilon) (p+1)}{4c_{\min} } \cdot \frac{1}{\sqrt{n}} \rho_n^{-1} s_+.
	\]
\end{proof}

\subsection{Condition \ref{Eq: correct exclusion}}

\begin{lemma}\label{Lem: Condition 10b}
	Under assumptions \ref{Assum: minimum EW directed}, \ref{Assum: no approximation error},
	\ref{Assum: model selection assumption}, \ref{Assum: model selection lambda},
	for $n$ large enough, with probability at least
	\begin{align*}
		\begin{split}
			1 &- 4 \cdot \exp \left(   - \frac{N \bar{\lambda}^2}{18}  + \log(n) \right) - 2(p+1) \cdot \exp \left( - \frac{N \bar{\lambda}^2}{18 (1 \vee c^2)}\right) - p(p+3) \exp\left( - N \frac{c_{\min}^2}{2048 s^{2}_+\tilde{ c}}  \right),
		\end{split}
	\end{align*}
	which tends to one as long as $- \frac{N \bar{\lambda}^2}{18} + \log(n) \rightarrow - \infty$, as $n$ tends to infinity, we have
	\begin{align*}
		\Vert \constrz_{1:2n, S^c} \Vert_\infty < 1.
	\end{align*}
\end{lemma}

\begin{proof}
	By equation \eqref{Eq: bound gradient L wrt vartheta}, Lemmas \ref{Lem: incoherence condition}, \ref{Lem: bound on bar R}, \ref{Lem: bound on model selection bias}, with the probability given in those Lemmas, for any $\epsilon > 0$,
	\begin{align*}
		\Vert \constrz_{1:2n, S^c} \Vert_\infty \le &\left\{ \Vert Q_{S^c,S} Q_{S,S}^{-1} \Vert_\infty + 1  \right\} \epsilon \\
		& + \left\{ \Vert Q_{S^c,S} Q_{S,S}^{-1} \Vert_\infty + 1  \right\} 	 \frac{324 (1 \vee (c^2p)) (1+\epsilon)^2 }{c_{\min}^2 } \cdot \sqrt{n} \rho_n^{-2} s_+^2 \bar{\lambda} \\
		& + \left\{ \Vert Q_{S^c,S} Q_{S,S}^{-1} \Vert_\infty + 1  \right\}  \frac{9 (1 \vee c) (1+\epsilon) (p+1)}{4c_{\min} } \cdot \frac{1}{\sqrt{n}} \rho_n^{-1} s_+ \\
		&+ \frac{1}{2} 	\Vert Q_{S^c,S} Q_{S,S}^{-1} \Vert_\infty.
	\end{align*}
	By Lemma \ref{Lem: incoherence condition}, for $n$ sufficiently large, we have $ \Vert Q_{S^c,S} Q_{S,S}^{-1} \Vert_\infty < 1/2$. Thus, by equation \eqref{Eq: bound gradient L wrt vartheta}, Lemmas \ref{Lem: bound on bar R} and \ref{Lem: bound on model selection bias}, for $n$ sufficiently large, with the prescribed probability,
	\begin{align*}
		\Vert \constrz_{1:2n, S^c} \Vert_\infty \le &\frac{3}{2} \epsilon + \frac{1}{4}\\
		& + \frac{486 (1 \vee (c^2p)) (1+\epsilon)^2 }{c_{\min}^2 } \cdot \sqrt{n} \rho_n^{-2} s_+^2 \bar{\lambda} \\
		& +  \frac{27 (1 \vee c) (1+\epsilon) (p+1)}{8c_{\min} } \cdot \frac{1}{\sqrt{n}} \rho_n^{-1} s_+.
	\end{align*}
	Pick $\epsilon = 1/3$, to obtain
	\begin{align*}
		\Vert \constrz_{1:2n, S^c} \Vert_\infty \le &\frac{3}{4} + \frac{486 (1 \vee (c^2p)) (4/3)^2 }{c_{\min}^2 } \cdot \sqrt{n} \rho_n^{-2} s_+^2 \bar{\lambda} +  \frac{27 (1 \vee c) (4/3) (p+1)}{8c_{\min} } \cdot \frac{1}{\sqrt{n}} \rho_n^{-1} s_+.
	\end{align*}
	The second and third term go to zero as $n$ tends to infinity by assumption \ref{Assum: model selection assumption}. Indeed, the second term is assumption \ref{Assum: model selection assumption} exactly. For the third term note that
	by assumption \ref{Assum: model selection assumption}, $\sqrt{n} s_+^2 \bar{\lambda} \rho_n^{-2} = n^{-1/2}\rho_n^{-1}s_+ \cdot n\rho_ns_+\bar{\lambda} \rightarrow 0$ as, $n \rightarrow \infty$. On the other hand, by assumption \ref{Assum: model selection lambda}, $n\rho_ns_+\bar{\lambda} \ge C\rho_n^{-1}s_+\log(n) \rightarrow \infty$. Therefore it must hold that $n^{-1/2}\rho_n^{-1}s_+ \rightarrow 0$. The claim follows.
\end{proof}

\subsection{Proof of Theorem \ref{Thm: model selection consistency}}

\begin{proof}[of Theorem \ref{Thm: model selection consistency}]
	By Lemma \ref{Lem: Condition 10b}, we know that with probability at least as large as 
	\begin{align*}
		\begin{split}
			1 &- 4 \cdot \exp \left(   - \frac{N \bar{\lambda}^2}{18}  + \log(n) \right) - 2(p+1) \cdot \exp \left( - \frac{N \bar{\lambda}^2}{18 (1 \vee c^2)}\right) - p(p+3) \exp\left( - N \frac{c_{\min}^2}{2048 s^{2}_+\tilde{ c}}  \right),
		\end{split}
	\end{align*}
	property \eqref{Eq: correct exclusion} holds for the construction $(\constrtheta, \constrz)$. Thus, by Lemma \ref{Lem: unique active set}, $\hat{S} = S^\dagger$ and in particular $\hat{S} \cap S^c = \emptyset$. 
	
	For the second part of Theorem \ref{Thm: model selection consistency}, recall that 
	by equation \eqref{Eq: 31a},
	\begin{equation}\label{Eq: 31a - again}
		\constrvartheta_S = \bar{\vartheta}_{0,S} + Q_{S,S}^{-1} \left\{ - \frac{1}{N} 	\left(\nabla_{\bar{\vartheta}} \bar{\mathcal{L}}(\bar{\theta}_0)\right)_S+ \frac{1}{N} \bar{X}_S^TW_0^2\left[
		\begin{array}{c|c}
			\textbf{1} & Z
		\end{array}
		\right] (\xi^\dagger - \xi_0) - \bar{\lambda} \constrz_{1:2n, S} + \bar{R}_S \right\}.
	\end{equation}
	Thus, $S^\dagger$ contains all those indices $i$ with
	\[
	\left\Vert  Q_{S,S}^{-1} \left\{ - \frac{1}{N} 	\left(\nabla_{\bar{\vartheta}} \bar{\mathcal{L}}(\bar{\theta}_0)\right)_S+ \frac{1}{N} \bar{X}_S^TW_0^2\left[
	\begin{array}{c|c}
		\textbf{1} & Z
	\end{array}
	\right] (\xi^\dagger - \xi_0) - \bar{\lambda} \constrz_{1:2n, S} + \bar{R}_S \right\} \right\Vert_\infty < \bar{\vartheta}_{0,i}.
	\]
	Hence, consider
	\begin{align*}
		&\left\Vert  Q_{S,S}^{-1} \left\{ - \frac{1}{N} 	\left(\nabla_{\bar{\vartheta}} \bar{\mathcal{L}}(\bar{\theta}_0)\right)_S+ \frac{1}{N} \bar{X}_S^TW_0^2\left[
		\begin{array}{c|c}
			\textbf{1} & Z
		\end{array}
		\right] (\xi^\dagger - \xi_0) - \bar{\lambda} \constrz_{1:2n, S} + \bar{R}_S \right\} \right\Vert_\infty \\ & \quad \quad \le \Vert Q_{S,S}^{-1}  \Vert_\infty \left\Vert  - \frac{1}{N} 	\left(\nabla_{\bar{\vartheta}} \bar{\mathcal{L}}(\bar{\theta}_0)\right)_S+ \frac{1}{N} \bar{X}_S^TW_0^2\left[
		\begin{array}{c|c}
			\textbf{1} & Z
		\end{array}
		\right] (\xi^\dagger - \xi_0) - \bar{\lambda} \constrz_{1:2n, S} + \bar{R}_S  \right\Vert_\infty \\
		& \quad \quad \le 2 \rho_n^{-1} \cdot \frac{n-1}{n - \max \{s_\alpha, s_\beta\}} \\
		& \quad \quad \quad \bigg\{ \epsilon\bar{\lambda} + \bar{\lambda} \\
		& \quad \quad \quad + \frac{9 (1 \vee c) (1+\epsilon) (p+1)}{4c_{\min} } \cdot \frac{1}{\sqrt{n}} \rho_n^{-1} s_+ \bar{\lambda} \\
		&\quad \quad \quad +  \frac{324 (1 \vee (c^2p)) (1+\epsilon)^2 }{c_{\min}^2 } \cdot \sqrt{n} \rho_n^{-2} s_+^2 \bar{\lambda}^2 \bigg\}
	\end{align*}
	where we used \eqref{Eq: bound Q_SS inverse}, \eqref{Eq: bound gradient} and Lemmas \ref{Lem: bound on bar R} and \ref{Lem: bound on model selection bias}. 
	
	By assumption $\bar{\lambda} \le C \cdot \sqrt{\log(n)/N}$ for some $C > 0$, thus the first two terms in the bracket may be upper bound by $C \cdot \sqrt{\log(n)/N}$, for a possibly different $C$. The third term is $o(1) \cdot 1 /n$ by assumption \ref{Assum: model selection assumption} and the last term is $o(1) \cdot \sqrt{\log(n)} / n$ by assumption \ref{Assum: model selection assumption}. Since $(n-1) / (n - \max \{s_\alpha, s_\beta\}) = O(1)$, the entire right-hand side is less or equal
	$C \rho_n^{-1} \frac{\sqrt{\log(n)}}{n}.$ 
	Multiply \eqref{Eq: 31a - again} by $\sqrt{n}$ to transition to the unscaled parameters $\vartheta^\dagger_S$ and the claim follows.
	In particular, $	C \rho_n^{-1} \frac{\sqrt{\log(n)}}{n}$ 
	goes to zero as $n$ tends to infinity, which implies that for $n$ large enough, with at least the prescribed probability the construction fulfils \eqref{Eq: correct inclusion} and thus
$	\hat{S} = S^\dagger = S.$
\end{proof}

\section{Sparse $\beta$-models and the power law}\label{sec:power law}

In this appendix we show that the degrees in sparse $\beta$-models can exhibit power law distributions. For that we leverage the results in \cite{Britton:etal:2006} that show that the $\beta$-model can generate node degrees asymptotically following a power law  and the empirical degree distribution converging in probability to the same power law if $\beta_{i}$ are randomly generated in a suitable way. Recall that in the $\beta$-model, each node is associated with a degree heterogeneity parameter $\beta_i$ and links are formed independently with probabilities
\begin{equation}
	\label{eq:betamodel}
	P(A_{ij}=1)=p_{ij}=\frac{ e^{\beta_i + \beta_j} }{ 1 + e^{\beta_i + \beta_j} }.
\end{equation}
We show the result for the model introduced in \cite{Chen:etal:19}, which is an undirected version of out model without covariates. Recall that in this model, links are made independently with probabilities
\begin{equation}\label{eq:sbm}
	P(A_{ij}=1)=p_{ij}=\frac{ e^{\mu+\beta_i+\beta_j} }{ 1 + e^{\mu+\beta_i+\beta_j} }.
\end{equation}

\begin{proposition}[S$\beta$M and power law]
	\label{prop:powerlaw}
	Let $\{ W_{i} \}_{i=1}^{\infty}$ be i.i.d.  random variables supported in $[1,\infty)$ with $P(W_{1} > w) \sim cw^{-\rho}$ as $w \to \infty$ for some $c > 0$ and $\rho \in (0,1)$. 
	For the S$\beta$M in \eqref{eq:sbm}, 
	suppose that $\mu = -\rho^{-1} \log n$ and $\beta_{i}$'s  are generated as $\beta_{i} = \log W_{i}$. 
	Then the limiting distribution of  each node degree $d_i$ as $n \to \infty$ is a power law with exponent $\tau = 2$, that is, 
	\[ 
	p_{k} := \lim_{n \to \infty} P(d_i=k) \sim k^{-2}, \ k \to \infty.
	\]
	In addition, for $N_{k} := | \{ i \in \{ 1,\dots, n \} : d_{i} = k \}|$, we have $N_{k}/n \stackrel{P}{\to} p_{k}$ as $n \to \infty$. 
\end{proposition}

In the above proposition, we do not impose sparsity on $\beta$, but it is possible to do so by assuming a mass at $1$ to the distribution of $W_{1}$ since the assumption only requires the tail of the distribution of $W_{1}$ to behave like $cw^{-\rho}$. The proposition follows from the results of \cite{Britton:etal:2006}, which are restated in the following.
We first recall the definition of a mixed Poisson distribution. 

\begin{definition}[Mixed Poisson distribution]
	Let $F$ be a distribution function supported in $\R_{+}$. A random variable $X$ taking values in the nonnegative integers follows the  \textit{mixed Poisson distribution} with mixing  distribution $F$ if
	\[
	P(X=k)=\int_{[0,\infty)} e^{-w}\frac{w^k}{k!} dF(w), \quad k =0,1,2,\dots.
	\] 
	If $W \sim F$, then we also say that $X$ follows the mixed Poisson distribution with parameter $W$. 
\end{definition}

The next lemma shows that the tail behavior of a mixed Poisson distribution is determined solely by that of the mixing distribution. 
\begin{lemma}\label{lemma:mixed-Poisson}
	Let $F$ be a distribution function supported in $\R_{+}$ such that $c_1 x^{1-\tau} \le 1-F(x) \le c_2 x^{1-\tau}$ for large $x$ for some $0 < c_1 < c_2 < \infty$. Then
	there exist $0 < c'_1 < c'_2 < \infty$ such that the distribution function $G$ of a mixed Poisson distribution with mixing distribution $F$ satisfies $c'_1x^{1-\tau} \le 1-G(x) \le c'_2x^{1-\tau}$ for large $x$. 
\end{lemma}

\begin{proof}
	See \cite{Hofstad:2016} Exercise 6.12.
\end{proof}

The following results are taken from Theorem 3.2 and Proposition 3.1 in \cite{Britton:etal:2006}. In the following, the variables $d_{1},\dots,d_{n}$ are indeed a triangular sequence and hence should be indexed by $n$, but this is suppressed for the notational convenience. 

\begin{theorem}[$\beta$-model and mixed Poisson distribution]
	\label{thm:mixed-Poisson} 
	Let $\{ W_{i} \}_{i=1}^{\infty}$ be i.i.d. positive random variables with $P(W_{1} > w) \sim cw^{-\rho}$ as $w \to \infty$ for some $c > 0$ and $\rho \in (0,1)$. 
	For the $\beta$-model in \eqref{eq:betamodel}, suppose that $\beta_{1},\dots,\beta_{n}$ are generated as $\beta_{i} = \log W_{i} - (\log n)/(2\rho)$ for $i=1,\dots,n$ for each $n=1,2,\dots$. 
	Then:
	\begin{itemize}
		\item[(i)] The limiting
		distribution of  each node degree $d_{i}$ as $n \to \infty$ is the  mixed Poisson distribution with parameter
		$\varrho W_{1}^{\rho}$, where $\varrho = c\int_{0}^{\infty}(1+x)^{-2}x^{-\rho}dx$. 
		\item[(ii)] For  $N_k := |\{ i \in \{1,\dots,n \}: d_i=k \}|$, we have $N_{k}/n \stackrel{P}{\to} P(\varrho W_{1}^{\rho} = k)$
		as $n \to \infty$, where $\varrho$ appears in (i). 
	\end{itemize}
\end{theorem}

Combined with Lemma \ref{lemma:mixed-Poisson}, we know that 
\[
\begin{split}
	P(d_i \geq y) &\approx P(\varrho W_{1}^\rho \geq y) \\
	&= P (W_1 \ge (y/\varrho)^{1/\rho}) \\
	&\sim c\varrho y^{-1},
\end{split}
\]
so that the limiting distribution of $d_i$ is  a power law with exponent $\tau=2$ in the sense that 
\[ 
\lim_{n \to \infty} P(d_i=k) \sim k^{-2}, \ k \to \infty. 
\] 
Likewise,  the empirical degree distribution  converges in probability to the same power law, which proves Proposition \ref{prop:powerlaw}.

\bibliography{bib}

\end{document}